\newcommand{\E}[0]{{\sf E}}

\newcommand{\qq}[0]{d}

\documentclass[letterpaper,11pt]{article}

\usepackage{amsthm}
\usepackage{amsmath}

\usepackage{amsfonts}
\usepackage{amscd,amsthm,amsfonts,amsopn,amssymb,verbatim}

\newtheorem{thm}{Theorem}[section]

\newtheorem{lemma}[thm]{Lemma}
\newtheorem{prop}[thm]{Proposition}
\newtheorem{cor}[thm]{Corollary}

\newtheorem{conj}[thm]{Conjecture}

\newcommand{\beq}[1]{\begin{equation}\label{#1}}
\newcommand{\enq}[0]{\end{equation}}

\newcommand{\bn}[0]{\bigskip\noindent}
\newcommand{\mn}[0]{\medskip\noindent}
\newcommand{\nin}[0]{\noindent}

\newcommand{\sub}[0]{\subseteq}
\newcommand{\sm}[0]{\setminus}
\renewcommand{\dots}[0]{,\ldots,}

\newcommand{\ov}[0]{\overline}

\newcommand{\A}[0]{{\cal A}}
\newcommand{\B}[0]{{\cal B}}
\newcommand{\cee}[0]{{\cal C}}

\newcommand{\eee}[0]{{\cal E}}
\newcommand{\f}[0]{{\cal F}}
\newcommand{\g}[0]{{\cal G}}
\newcommand{\h}[0]{{\cal H}}
\newcommand{\I}[0]{{\cal I}}

\newcommand{\K}[0]{{\cal K}}

\newcommand{\Q}[0]{{\cal Q}}
\newcommand{\Qq}[0]{{\cal Q}^*}
\newcommand{\R}[0]{{\cal R}}
\newcommand{\sss}[0]{{\cal S}}
\newcommand{\T}[0]{{\cal T}}

\newcommand{\ra}[0]{\rightarrow}
\newcommand{\Ra}[0]{\Rightarrow}

\newcommand{\inte}[0]{{\rm int}}
\newcommand{\crit}[0]{{\rm crit}}
\newcommand{\ext}[0]{{\rm ext}}
\newcommand{\deff}[0]{{\rm def}}

\newcommand{\bb}[0]{b}
\newcommand{\cc}[0]{\mbox{{\sf c}}}

\newcommand{\ttt}[0]{t}

\newcommand{\aaa}[0]{\mbox{{\sf a}}}
\newcommand{\bbb}[0]{\mbox{{\sf b}}}
\newcommand{\ww}{\mbox{{\sf w}}}
\newcommand{\zz}[0]{\mbox{{\sf z}}}

\newcommand{\xxx}[0]{B}
\newcommand{\PPP}[0]{\mathbb P}

\newcommand{\0}[0]{\emptyset}

\renewcommand{\qed}[0]{\begin{flushright} \rule{2mm}{3mm} \end{flushright}}

\def\qqqed{\null\nobreak\hfill\hbox{\rule{2mm}{3mm} }\par\smallskip}

\newcommand{\C}[2]{{{#1}\choose{{#2}}}}
\newcommand{\Cc}[0]{\tbinom}
\newcommand{\ga}[0]{\alpha }
\newcommand{\gb}[0]{\beta }
\newcommand{\gc}[0]{\gamma }
\newcommand{\gd}[0]{\delta }
\newcommand{\gD}[0]{\Delta }
\newcommand{\gG}[0]{\Gamma }
\newcommand{\gk}[0]{\kappa }
\newcommand{\gl}[0]{\lambda }
\newcommand{\gL}[0]{\Lambda}
\newcommand{\go}[0]{\omega}
\newcommand{\gO}[0]{\Omega}

\newcommand{\gs}[0]{\sigma}
\newcommand{\gS}[0]{\Sigma}
\newcommand{\gz}[0]{\zeta}
\newcommand{\eps}[0]{\varepsilon }
\newcommand{\vt}[0]{\vartheta}
\newcommand{\vs}[0]{\varsigma}
\newcommand{\vr}[0]{\varrho}
\newcommand{\vp}[0]{\varphi}

\newcommand{\glr}[0]{\gL_r(n,p)}
\newcommand{\glq}[0]{\gL_r(n,q)}

\newcommand{\comments}[1]{}
\usepackage[usenames]{color}

\begin{document}

\renewcommand{\thefootnote}{\fnsymbol{footnote}}

\title{Tur\'an's Theorem for random graphs}
\author{B. DeMarco\footnotemark $~$ and J. Kahn\footnotemark
}
\date{}
\footnotetext{ * supported by the U.S.
Department of Homeland Security under Grant Award Number 2007-ST-104-000006.}
\footnotetext{ $\dag$ Supported by NSF grants DMS0701175 and DMS1201337.}

\date{}

\maketitle

\begin{abstract}
For a graph $G$, denote by $t_r(G)$ (resp. $b_r(G)$)
the maximum size of a $K_r$-free (resp. $(r-1)$-partite)
subgraph of $G$.
Of course $t_r(G) \geq b_r(G)$ for any $G$,
and Tur\'an's Theorem
says that equality holds for complete graphs.
With $G_{n,p}$ the usual (``binomial" or
``Erd\H{o}s-R\'enyi") random graph, we show:

\mn
{\bf Theorem}
{\em For each fixed r there is a C such that if
\[
p=p(n) > Cn^{-\tfrac{2}{r+1}}\log^{\tfrac{2}{(r+1)(r-2)}}n,
\]
then $\Pr(t_r(G_{n,p})=b_r(G_{n,p}))\ra 1$ as $n\ra\infty$.}

\mn
This is best possible (apart from the value of $C$) and settles a question
first considered by Babai, Simonovits and Spencer about 25 years ago.

\end{abstract}


\section{Introduction}\label{Intro}

Write $\ttt_r(G)$
for the maximum size
of a $K_r$-free
subgraph of a graph $G$
(where a graph is $K_r$-{\em free} if it contains no copy
of the complete graph $K_r$ and size means number of edges),
and $b_r(G)$ for the maximum size of an
$(r-1)$-partite subgraph of $G$.
Of course $\ttt_r(G)\geq b_r(G)$ for any $G$,
while the classic theorem of
Tur\'an \cite{Turan}---commonly held to have
initiated extremal graph theory---says that equality holds
when $G$ is the complete graph $K_n$.

Here we are interested in understanding, for a given $r$,
when
equality is likely
to hold for the usual (``binomial" or ``Erd\H{o}s-R\'enyi") random graph
$G=G_{n,p}$---that is, for what $p=p(n)$
one has
\beq{tttbb}
\ttt_r(G_{n,p})=\bb_r(G_{n,p}) ~~~\mbox{{\em w.h.p.} }
\enq
(An event holds {\em with high probability}
(w.h.p.) if its probability tends to 1 as $n\ra\infty$.
Note \eqref{tttbb}
holds for sufficiently small $p$ for the silly reason
that $G$ is itself likely to be $(r-1)$-partite, but we are
thinking of more interesting values of $p$.)

First results on this problem were given by
Babai, Simonovits and Spencer \cite{BSS}
(apparently in response to a conjecture of Paul Erd\H{o}s
\cite{Pach}).
They showed that
for $r=3$---in which case Tur\'an's Theorem is actually
Mantel's \cite{Mantel}---\eqref{tttbb} holds when
$p>1/2$
(more precisely, when $p>1/2-\eps$ for some
fixed $\eps >0$),
and asked whether their result could be extended to $p > n^{-c}$
for some fixed positive $c$.
This was finally accomplished (with $c = 1/250$)
in an ingenious paper of
Brightwell, Panagiotou and Steger \cite{BPS},
which actually proved a similar statement for every (fixed)
$r$:
\begin{thm}[\cite{BPS}]\label{Tbps}
For each r there is a $c>0$ such that if $p> n^{-c}$
then w.h.p. every largest
$K_r$-free subgraph of $G_{n,p}$ is
$(r-1)$-partite.
\end{thm}

\nin
(Actually \cite{BSS} considers the problem with a
general forbidden graph $H$ in place of a clique---though
the discussion there is mostly
confined to $H$'s of chromatic number three---and
\cite{BPS} also suggests that
Theorem~\ref{Tbps} may hold for more than cliques;
see Section \ref{Remarks} for a little more on this.)

\medskip
It was also suggested in \cite{BPS} that when $r=3$,
$p> n^{-1/2+\eps}$ might suffice for \eqref{tttbb}, and
the precise answer in this case---\eqref{tttbb} holds for
$p>Cn^{-1/2}\log^{1/2} n$---was proved in
\cite{DKMantel}.
(The more conservative suggestion in \cite{BPS} seems due
to an excess of caution on the part of the authors,
who surely realized that $\Theta (n^{-1/2}\log^{1/2} n)$
is the natural guess \cite{StegerPC}.)

Here we settle the problem for every $r$:
\begin{thm}\label{MT}
For each r there is a C such that if
\beq{p}
p > Cn^{-\tfrac{2}{r+1}}\log^{\tfrac{2}{(r+1)(r-2)}}n,
\enq
then w.h.p. every largest $K_r$-free subgraph of $G_{n,p}$ is $(r-1)$-partite.
\end{thm}

\mn
This
is best possible (apart from the value of $C$),
basically because (formal proof omitted)
for smaller $p$
there are usually edges of $G:=G_{n,p}$ not lying in $K_r$'s;
and while these
automatically belong to all largest $K_r$-free
subgraphs of $G$, there's no reason to expect that
they are all contained
in every largest $(r-1)$-partite subgraph
(and if they are not, then $t_r(G)>b_r(G)$).

\mn
{\bf Context.}
One of the most interesting combinatorial directions of the
last few decades has been the study of
``sparse random" versions of classical results
(e.g. the theorems of Ramsey, Szemer\'edi and
Tur\'an)---that is,
of the extent to which these results remain true in a
random setting.  These developments, initiated by
Frankl and R\"odl
\cite{FR86} and the aforementioned Babai {\em et al.} \cite{BSS}
and given additional impetus by the ideas of R\"odl and
Ruci\'nski \cite{RR1,RR2} and
Kohayakawa, \L uczak and R\"odl \cite{KLRAP,KLR},
led in more recent years to a number of major
results, beginning with the breakthroughs of
Conlon and Gowers \cite{Conlon-Gowers}
and Schacht \cite{Schacht}.
The following are special cases, the second of which
will be needed for the
proof of Theorem~\ref{MT}.
\begin{thm}[\cite{Conlon-Gowers,Schacht}]\label{Tdensity}
For each $\vartheta>0$ there is a $K$ such that
if
\[
p > Kn^{-2/(r+1)}
\]
then w.h.p. $t_r(G_{n,p}) < (1-\frac{1}{r-1}+\vt)|G_{n,p}|$.
\end{thm}
\begin{thm}[\cite{Conlon-Gowers}]\label{8.34}
For each $\vartheta>0$ there is a $K$ such that
if
\beq{p'}
p > Kn^{-2/(r+1)}
\enq
then w.h.p. each $K_r$-free
subgraph of $G=G_{n,p}$ of size at least $(1-\frac{1}{r-1})|G|$
can be made $(r-1)$-partite by deletion of at most
$\vartheta n^2p$ edges.
\end{thm}
\nin
These may be considered sparse random versions of
Tur\'an's Theorem and the ``Erd\H{o}s-Simonovits
Stability Theorem" \cite{Erdos,Simonovits} respectively.
Both were conjectured by Kohayakawa {\em et al.}
\cite{KLR}, who proved Theorem~\ref{8.34} for $r=3$,
the weaker Theorem  \ref{Tdensity} for $r=3$ having been proved earlier
by Frankl and R\"odl \cite{FR86}.
(See also \cite{GSS,Gerke} for further progress preceding
Theorems~\ref{Tdensity} and \ref{8.34},
and \cite{Samotij} for a common generalization of
\cite{Conlon-Gowers} and \cite{Schacht}.)

Even more recently, related (but independent) papers of
Balogh, Morris and Samotij \cite{BMS}
and Saxton and Thomason \cite{ST} prove remarkable ``container" theorems---more
asymptotic counting than
probabilistic methods---which, {\em once established,}
yield surprisingly simple proofs of many of the
very difficult results mentioned above.
See also \cite{Rodl-Schacht} for a survey of these
and related developments.

Though it does finally establish
the ``true" random analogue
of Tur\'an's Theorem,
one cannot really say
that Theorem~\ref{MT}
is the culmination of some of this earlier work.
First,
it does not quite imply Theorem~\ref{Tdensity},
whose conclusion holds for $p$ in a somewhat larger range, and
its conclusion is not comparable to that of Theorem~\ref{8.34}.
(Of course it 
is much stronger than Theorem~\ref{Tdensity} in the range where
it does apply.)
Second,
apart
from a black-box application of Theorem~\ref{8.34},
the problem addressed by Theorem~\ref{MT}
seems immune to the powerful ideas developed to prove
the aforementioned results.
(Conversations with several interested parties support
this opinion and
suggest
that the paucity of
results in the direction of Theorem~\ref{MT}
is not due to lack of effort.)

\mn
{\bf Plan.}
We prove Theorem~\ref{MT} only for $r\geq 4$;
the proof could presumably be adapted to $r=3$, but
this seems pointless given that we already have
the far simpler argument of \cite{DKMantel}.

We begin with terminology and such in Section \ref{Usage},
but defer further preliminaries
in order to give an early idea
of where we are headed.
Thus Section~\ref{Skeleton}
just states the main points---Lemmas \ref{L2.2} and
\ref{L2.3}---underlying
Theorem~\ref{MT}
and
shows how they imply the theorem.

Section~\ref{PB} then
collects machinery needed for the arguments to come.
One new item here
is Lemma~\ref{TRW}, an extension of the recent
Riordan-Warnke generalization
\cite{RW} of the Janson Inequalities
\cite{Janson}, that seems likely to be useful elsewhere.

We next, in Sections \ref{S2.2} and \ref{S2.3},
outline the proofs of
Lemmas~\ref{L2.2} and \ref{L2.3},
again meaning we state main points and derive the lemmas
from them.  The assertions underlying
Lemma~\ref{L2.2} are proved in
Sections \ref{PLxKLS}-\ref{PL2.2B}
and those
underlying
Lemma~\ref{L2.2} in
Sections \ref{RandC}-\ref{PL2.3'}.
(The two parts both require Lemma~\ref{xKLS} and
the material of Section~\ref{PB}, but are otherwise independent.)
Finally, Section \ref{Remarks} mentions a few
related questions.

\mn
{\bf Discussion.}
The basic structure of the argument---deriving Theorem~\ref{MT} from
Lemmas~\ref{L2.2} and \ref{L2.3}---seems natural
and is similar to that in
\cite{Mantel}.
(See the remark following Lemma~\ref{L2.3}.
The reader familiar with \cite{Mantel} may notice that the
rather {\em ad hoc} conditions around the analogue of
$Q(\Pi)$---here defined in
the second paragraph after \eqref{po1}---have now
disappeared.)
It should, however, be stressed
that the nature and
difficulty of the problem undergo a drastic change
when we move from $r=3$ to $r\geq 4$, and that most of the ideas
of \cite{Mantel} are pretty clearly useless for present
purposes.  (This feels
akin to the familiar jump in difficulty
when one moves from graphs to hypergraphs.)
In the event, most of the
key ideas in what follows are without much in the way of antecedents,
the most notable exception being that the uses of
Harris' Inequality in Section~\ref{RandC} were inspired by a
related use in \cite{BPS}.

We will try to say a little more about various aspects
of the argument when we are in a position to do so intelligibly.
The most interesting points
are the proof of
Lemma~\ref{L2.2B} (the last of the lemmas supporting
Lemma~\ref{L2.2};
what's most interesting here is how tricky
this innocent-looking statement was to prove)
and, especially, the several ideas developed in
Sections~\ref{RandC}-\ref{PL2.3'} to deal with
Lemma~\ref{L2.3}.

\section{Usage and definitions}\label{Usage}

For integers $a\leq b$, we use
$[a,b]$ for $\{a\dots b\}$ and $[b]$
for $[1,b]$ (assuming $b\geq 1$).
As usual, $2^X$ and $\C{X}{k}$ are
the collections of subsets and $k$-subsets of the set $X$.
We write
$\ga=(1\pm \gd)\gb$ for $(1-\gd)\gb<\ga<(1+\gd )\gb$
and $\log $ for natural logarithm.
Following a common abuse, we usually omit superfluous
floor and ceiling symbols.

We use
$B(n,p)$ for a random variable with the binomial
distribution ${\rm Bin}(n,p)$.
In line with recent practice, we occasionally use
$X_p$ for the ``binomial" random subset of $X$
given by
\beq{Xp}
\Pr(X_p=A) = p^{|A|}(1-p)^{|X\sm A|} ~~ ~~( A\sub X).
\enq

Throughout the paper $V=[n]$ 
is our default
vertex set.
The random graphs $G_{n,p}$
$G_{n,M}$
are defined as usual; see e.g. \cite{JLR}.
{\em We will usually use $G$ as an abbreviation for
$G_{n,p}$,}
so for the present discussion use $H$ for a general graph (on $V$).

We use $|H|$ for the {\em size} (i.e. number of edges) of
$H$,
$N_H(x)$ for the set of neighbors
of $x$ in $H$, $d_H(x)$
for the degree of $x$ in $H$ (i.e. $|N_H(x)|$),
$d_H(x,y)$ for $|N_H(x)\cap N_H(y)|$ and so on.
When the identity of $H$ is clear---usually
meaning $H=G=G_{n,p}$---we will sometimes drop the
subscript (thus $N(x)$ or $N_x$, $d(x)$ etc.) and may then,
a little abusively,
use, for example, $N_B(x)$ for the set of neighbors
of $x$ in $B\sub V$ or $N_L(x)$ for the set of vertices
joined to $x$ by members of $L\sub \C{V}{2}$.
We use $\gD_H$ for the maximum degree of $H$.

As usual, $H[A]$ is the subgraph of $H$ induced by $A\sub V$.
For disjoint $A_1\dots A_k\sub V$, we use $\nabla(A_1\dots A_k)$
for the set of pairs $\{x,y\}$ meeting two distinct
$A_i$'s, and often
write $\nabla_H(A_1\dots A_k)$
for $H\cap \nabla(A_1\dots A_k)$.
%
We will tend to use $xy$ ($=yx$), rather than
$\{x,y\}$, for an element of $\C{V}{2}$.
Unless stated otherwise, $V(L)$ is the set of vertices belonging to
members of $L\sub\C{V}{2}$.

\medskip
A {\em cut} is an ordered $(r-1)$-partition $\Pi=(A_1\dots A_{r-1})$
of $V$.  (The order of $A_2\dots A_{r-1}$ isn't important,
but $A_1$ will play a special role.)
{\em Throughout the paper $\Pi$ will denote
a cut.}
We say $\Pi$ is {\em balanced} if each of its blocks has size
$(1\pm \gd)n/(r-1)$, where $\gd$ is a small (positive) constant
(see the discussion at the end of this section).

For $\Pi$ as above we sometimes use $\ext(\Pi)$
for $\nabla(A_1\dots A_{r-1})$ and $\inte(\Pi)$ for
$\C{V}{2}\sm \ext(\Pi)$
(and give $\ext_H(\Pi)$,
$\inte_H(\Pi)$ their obvious meanings).
We will also use $|\Pi|$ for
$|\nabla(A_1\dots A_{r-1})|$ and
$|\Pi_H|$ for
$|\nabla_H(A_1\dots A_{r-1})|$;
thus
$b_r(H) = \max_\Pi|\Pi_H|$.
The {\em defect of $\Pi$ with respect to $H$} is
\[\deff_H(\Pi) = b_r(H) -|\Pi_H|,\]
and the {\em defect} of $\Pi$ is its defect with respect to $G=G_{n,p}$.

\medskip
Though it may take some getting used to,
the following notation will
be quite helpful.
Suppose that for $i=1\dots s$, $X_i$ is a collection of
$a_i$-subsets of $V$
(we will usually have $a_i\leq 2$)
and that $\sum a_i\leq r$.
We then write
$\kappa_H(X_1\dots X_s)$ for the number of ways to choose
disjoint $Y_1\in X_1\dots Y_s\in X_s$ and
an $(r-\sum a_i$)-subset $Z$ of $V\sm \cup Y_i$ so that
\[
\C{Y_1\cup\cdots \cup Y_s\cup Z}{2}
\sm \bigcup_{i=1}^s\C{Y_i}{2}\sub H.
\]
When $X_i$ consists of a single set, say $\{x_1\dots x_{a_i}\}$,
we omit set brackets and commas in the specification;
for example:  (i) $\kappa_H(xy)$ counts choices of $Z\in \C{V}{r-2}$
with all pairs from $Z\cup \{x,y\}$ other than (possibly)
$\{x,y\}$ belonging to $H$,
and
(ii) $\kappa_H(x_1x_2x_3,T)$, with $T\sub \C{V}{2}$,
counts choices of $\{x_4,x_5\}\in T$ with
$\{x_4,x_5\}\cap \{x_1,x_2,x_3\}=\0$
and
$\{x_6\dots x_r\}\sub V\sm \{x_1\dots x_5\}$
(with $x_4\neq x_5$ and $x_6\dots x_r$ distinct)
such that all members
of $\C{\{x_1\dots x_r\}}{2}$
other than those in
$\C{\{x_1\dots x_3\}}{2}\cup \{\{x_4,x_5\}\}$
lie in $H$.

In one special case, when $s=r-1$, $a_1=2$, $X_2\dots X_s$
are disjoint and no pair from $X_1$ meets $X_2\cup\cdots\cup X_s$,
we will on a few occasions use $\K_H(X_1\dots X_s)$ for the collection
counted by $\kappa_H(X_1\dots X_s)$
(members of which may be thought of as copies of $K_r^-$,
the graph obtained from $K_r$ by deleting an edge).

When $H=G=G_{n,p}$, we will tend to
drop subscripts and write simply $\kappa(\cdots)$
and $\K(\cdots)$.

\medskip
The quantity
\beq{gLr}
\gL_r(n,p) := n^{r-2}p^{\Cc{r}{2}-1},
\enq
which, up to scalar, is
the expectation of $\kappa(xy)$ (for given $x,y$),
will appear frequently (so we give it a name).

\mn
{\bf Constants.}
There will be quite a few of these, but not so many that are more than local.
The most important are $\gd$ (see the above definition of a balanced cut);
$\gc$ (used in the definition of a ``bad" pair for a given cut following \eqref{po1});
$\ga$ (see ``rigidity" in Section~\ref{RandC});
and $C$ (in \eqref{p}).
The few constants that are given explicitly will, superfluously, be subscripted by $r$.
For the main constants, apart from an explicit constraint on
$\gc$ in Section~\ref{S2.3}
(see \eqref{gcsat}), we will not bother with actual values,
but the hierarchy is (of course) important:  we assume
$C^{-1} \prec \gd \prec \ga,\gc$
(where, just for the present discussion, ``$a\prec b$" means
$a$ is small enough relative to $b$ to support our arguments),
with $\ga$ (and $\gc$, but this will follow from \eqref{gcsat})
small relative to $r$.
(We could take $\ga=\gc$,
but prefer to distinguish them to emphasize their separate roles.)
In particular the constant $C$---and $n$---are always assumed to be large enough to
support our various assertions.

\section{Skeleton}\label{Skeleton}

In this section we state the main points underlying
Theorem~\ref{MT} and derive the theorem from these
(with a small assist from one of the standard large
deviation assertions of Section \ref{PB}).

We fix $r\geq 4$ and
assume
$p$ is as in \eqref{p}
with
$C$ a suitable constant
(and, as always,
$n$ large enough to support our arguments).
Though not really necessary, it will also be
convenient to assume, as we may by
Theorem~\ref{Tbps}, that
\beq{po1}
p=o(1).
\enq

Fix $\gc>0$ to be specified below.
(The specification will make more sense in
Section \ref{S2.3}, where we outline the proof of
Lemma~\ref{L2.3}, so we postpone it until then.)

For disjoint
$S_1\dots S_{r-2}\sub V$, a pair
$\{x,y\}$ is
{\em bad for $(S_1\dots S_{r-2})$ in} $G$ if
$\kappa_G(xy,S_1\dots S_{r-2})< \gc \glr$.
For $\Pi=(A_1\dots A_{r-1})$, we write
$Q_G(\Pi)$ for the set of pairs from $A_1$ that are
bad for $(A_2\dots A_{r-1})$ in $G$, and for
$F\sub \C{V}{2}$ let
\[\vp(F,\Pi) = (r-1)|F[A_1]| +|F\cap \ext(\Pi)|. \]

We now write $G$ for $G_{n,p}$.
The next two statements are our main points.

\begin{lemma}\label{L2.2}
There is an $\eta>0$ such that
w.h.p.
\[\vp(F,\Pi) < |\Pi_G|\]
whenever $\Pi=(A_1\dots A_{r-1})$
is balanced
and $F\sub G$ is
$K_r$-free and satisfies $F\not\supseteq \Pi_G$, $F\cap Q_G(\Pi)=\0$,
\beq{HP2}
|F[A_1]|<\eta n^2p,
\enq
and
\beq{HP3}
|N_F(x)\cap A_1|=\min\{|N_F(x)\cap A_i|:i\in [r-1]\}
~~\forall x\in A_1.
\enq
\end{lemma}

\begin{lemma}\label{L2.3}
W.h.p. $\deff_G(\Pi) \geq r |Q|$ whenever
$\Pi=(A_1\dots A_{r-1})$
is balanced, $Q\sub Q_G(\Pi)$
and
\beq{dQx'}
d_Q(x) \leq \min\{|N_G(x)\cap A_i|:2\leq i\leq r-1\}
~~\forall x\in A_1.
\enq
\end{lemma}
\nin
(The lemma holds with any constant in place of
$r$ in the defect bound, but $r$ (actually $r-1$) is
what's needed in
the final inequality \eqref{final2} below.)

\mn
{\em Remark.}
Technicalities aside, the dichotomy embodied in
Lemmas \ref{L2.2} and \ref{L2.3} is quite natural.
If $\Pi=(A_1\dots A_{r-1})$ is a cut and
$xy\in G[A_1]$ (say), then any
$K_r$-free $F$ $\sub G$ containing $xy$ must miss
at least one edge of each member of $ \K_G(xy,A_2\dots A_{r-1})$.
For a typical $xy$ there are (by our choice of $p$) many
of these, and one may hope that this forces
$\ext_G(\Pi)\sm F$ to be (much)
larger than the number of such $xy$'s in $F$,
which gives $|F|< |\Pi_G|$ ($\leq b_r(G)$) provided
a decent fraction of the edges of
$F\cap {\rm int}(\Pi)$ are ``typical" edges of $G[A_1]$.
Something of this sort is shown in Lemma~\ref{L2.2}.

Of course for a general $\Pi $ and
$xy$ as above,
$ \kappa_G(xy,A_2\dots A_{r-1})$ need not be large, or even
positive.
This more interesting situation---in which
membership of $xy$ in $F$ says
less about $\ext_G(\Pi)\sm F$---is handled by
Lemma~\ref{L2.3}, which says, roughly, that the defect of $\Pi$
is large relative to the number of
pairs $x,y$---or edges, but adjacency of $x,y$
is irrelevant here---from $A_1$ for which
$ \kappa_G(xy,A_2\dots A_{r-1})$ is small.

Notice, for example, that $t_r(G)>b_r(G)$ whenever there are a maximum
cut $(A_1\dots A_{r-1})$ and $xy\in G[A_1]$ with
$ \kappa_G(xy,A_2\dots A_{r-1})=0$; thus a baby requirement for
Theorem~\ref{MT} is that this situation be unlikely, and
in fact we don't know how to show even this much without some
portion of the machinery of Sections \ref{RandC}-\ref{PL2.3'}.

\medskip
At any rate,
given Lemmas \ref{L2.2} and \ref{L2.3} we finish easily,
as follows.
Let $F_0$ be a largest $K_r$-free subgraph of $G$ and
$\Pi=(A_1\dots A_{r-1})$
a cut maximizing $|F_0\cap \ext(\Pi)|$, with
$|F_0[A_1]|=\max_i |F_0[A_i]|$.
Then \eqref{HP3} holds
with $F_0$ in place of $F$ (if it did not,
we could move a violating
$x$ from $A_1$ to some other $A_i$
to increase $|F_0\cap \ext(\Pi)|$),
and Theorem~\ref{8.34} implies that w.h.p.
$F_0$ also satisfies \eqref{HP2} (actually with $o(1)$ in place of $\eta$;
here we use the standard observation that
$b_r(H)\geq (r-2)|H|/(r-1)$ for any $H$).
Moreover $\Pi$ is balanced (actually with $o(1)$
in place of $\gd$ in the definition of balance)
w.h.p., since (w.h.p.)

\begin{eqnarray}
|\nabla(A_1\dots A_{r-1})|p +o(n^2p)
&>&
|\Pi_G| \label{bal1}\\
&\geq &|F_0\cap \ext(\Pi)|\nonumber\\
&>&|F_0| -o(n^2p)\label{bal2} \\
&> &(r-2)|G|/(r-1)-o(n^2p)\label{bal3}\\
&> &(r-2)n^2p/[2(r-1)]-o(n^2p)\label{bal4},
\end{eqnarray}
so that
\[
|\nabla(A_1\dots A_{r-1})| >(r-2)n^2/[2(r-1)]-o(n^2),
\]
which easily gives $|A_i|= (1\pm o(1))n/(r-1)$ $\forall i$.
Here \eqref{bal1} and \eqref{bal4} are easy applications of
``Chernoff's Inequality" (Theorem~\ref{Chern});
\eqref{bal2} follows from Theorem~\ref{8.34};
and \eqref{bal3} is the ``standard observation" mentioned above.

\medskip
Let
$F_1 =F_0[A_1]\cup (F_0\cap \ext(\Pi))$ and $F=F_1\sm Q_G(\Pi)$.
Noting that these modifications introduce no $K_r$'s
and preserve
\eqref{HP2} and \eqref{HP3},
we have, w.h.p.,
\begin{eqnarray}
\ttt_r(G) &=&|F_0|\nonumber\\
& \leq & \vp(F_1,\Pi)\nonumber\\
& =& \vp(F,\Pi) + (r-1)|F_1\cap Q_G(\Pi)|\nonumber\\
&\leq & |\Pi_G| + (r-1)|F_1\cap Q_G(\Pi)|\label{final1}\\
&\leq &\bb_r(G),\label{final2}
\end{eqnarray}
where \eqref{final1} is given by Lemma~\ref{L2.2}
(note that if $F\supseteq \Pi_G$ then $F\cap Q_G(\Pi)=\0$
implies $F=\Pi_G$)
and \eqref{final2} by Lemma~\ref{L2.3}
(applied with $Q=F_1\cap Q_G(\Pi)$, noting that
\eqref{dQx'} follows from the fact that \eqref{HP3}
holds for $F_1$).

This gives
\eqref{tttbb}.
For the slightly stronger assertion of Theorem~\ref{MT},
notice that we have strict inequality in
\eqref{final1} unless $F=\Pi_G$
and in \eqref{final2} unless $F_1\cap Q_G(\Pi)=\0$.
Thus $|F_0|=b_r(G)$ implies
$F_0[A_1]=F[A_1]\cup (F_1\cap Q_G(\Pi))=\0$,
so also $F_0[A_i]=\0$ for $i\geq 2$
(since we assume $|F_0[A_1]|=\max_i |F_0[A_i]|$).
\qed

\section{Preliminaries}\label{PB}

The following version of
Chernoff's Inequality may be found, for example, in
\cite[Theorem 2.1]{JLR}.
\begin{thm}\label{Chern}
For $\xi =B(m,p)$, $\mu=mp$ and any $\gl\geq 0$,
\begin{eqnarray*}
\Pr(\xi >\mu+\gl)& < &\exp [- \tfrac{\gl^2}{2(\mu+\gl/3)}],\\
\Pr(\xi < \mu-\gl) &<& \exp [- \tfrac{\gl^2}{2\mu}].
\end{eqnarray*}
\end{thm}

\medskip
We will also need the following inequality for weighted sums of Bernoullis,
which can be derived from, for instance, \cite[Lemma 8.2]{Beck-Chen}.

\begin{lemma}\label{Chernwts}
Suppose $w_1\dots w_m \in [0,z]$.
Let $\xi_1\dots \xi_m$ be independent Bernoullis,
$\xi = \sum \xi_i w_i$, and $\E\xi \leq\psi$.
Then for any $\eta \in [0,1]$ and $\gl\geq \eta \psi$,
\[\Pr(\xi \geq \psi+\gl) \leq \exp [- \eta\gl/(4z)].\]
\end{lemma}
\nin

\medskip
We
record a few easy consequences of
Theorem~\ref{Chern},
in which we again take $G=G_{n,p}$
(with $p$ as in \eqref{p}, which is more than is needed here).

\begin{prop}\label{vdegree}
W.h.p. for all $x,y\in V$,
\beq{Deg}
d(x)=(1\pm o(1))np ~~~\mbox{and}
~~~d(x,y)=(1\pm o(1))np^2.
\enq
\end{prop}

\begin{prop}\label{propG[X]}
{\rm (a)}
For each $\eps>0$ there is a K such that w.h.p.
\[
||G[X]|-|X|^2p/2| \leq \max\{\eps |X|^2p,K|X|\log n\}
~~\forall X\sub V.
\]

\mn
{\rm (b)}  There is a fixed $\eps>0$ such that w.h.p.
\[
|G[X]| <|X|\log n ~~\forall X\sub V ~\mbox{with}~
|X| < \eps p^{-1}\log n.
\]
\end{prop}

\begin{prop}\label{nablaST}
For all $\eps >0$ and c there is a K such that
w.h.p.
\[
|\nabla_G(S,T)| > (1-\eps)|S||T|p
\]
for all disjoint $S,T\sub V$
with $|S|>cn$ and $|T|>K/p$.
\end{prop}

\nin
We omit the straightforward proofs.

\subsection{Polynomial concentration}\label{SKV}

\mn
We will need two instances of the ``polynomial
concentration" machinery of J.H. Kim and V. Vu
\cite{KimVu,Vu,Vu2}.
Here we omit the polynomial
language and just recall what we actually use, for which
we assume the following setup.
Let $\h$ be a collection of $d$-subsets of $X=[N]$,
$\ww:\h\ra \Re^+$,
$Y=X_p$
(see
\eqref{Xp} for $X_p$)
and
\beq{xipoly}\mbox{$\xi = \sum\{\ww_A:A\in \h,A\sub Y\}$.}
\enq
For $L\sub X$ let
\[
\mbox{$\E_L =\sum\{\ww_Ap^{|A\sm L|}:L\sub A\in \h\}$}
\]
and
$\E_l =\max_{|L|=l}\E_L$ for $0\leq l<d$
(e.g. $\E_0=\E\xi$).

We will need the following particular consequences of \cite{KimVu,Vu2,Vu}.
(The first---as observed in \cite[Cor. 5.5]{JKV}, here slightly rephrased---follows
easily from results of \cite{KimVu} and \cite{Vu2}, and the second
is contained in
\cite[Cor. 2.6]{Vu}.)

\begin{lemma}\label{VuL1}
For each fixed $d$, $\eps>0$, b and $M$ there is a $J$
such that
if
$\max \ww_A\leq b$ and
$B\geq\max \{(1+\eps)\E_0,J \log N, N^\eps \max_{0<j<d}\E_j\}$,
then
\[\Pr(\xi >B)< N^{-M}.\]
\end{lemma}

\begin{lemma}\label{VuL2}
For each fixed $d$, $\eps>0$, b and $M$ there is a $J$
such that
if $\max \ww_A\leq b$ and
$\max_{0\leq l<d}\E_l < N^{-\eps}$, then
\[\Pr(\xi >J)< N^{-M}.\]
\end{lemma}

We will apply these results in the following setting.
(There is nothing surprising here---e.g. similar applications
of the above machinery appear in \cite{JKV}---but,
lacking a reference, we include a few details.)

Define a {\em rooted graph} to be a graph $H=(V(H),E(H))$
with members of some $R=R(H)\subset V(H)$ designated ``roots."
In what follows it will be convenient to fix
some ordering
of $R$ and speak of the {\em root sequence}, $(u_1\dots u_s)$,
of $H$.

Though
we allow edges between the roots, they play no role here
and we set $E'(H)=\{e\in E(H):e\not\sub R\}$,
$v_H=|V(H)\sm R|$, $e_H=|E'(H)|$
and $\rho(H) = e_H/v_H$, this last quantity
called the {\em density} of $H$.
(For typographical reasons we will sometimes use $v(H)$ and $e(H)$
in place of $v_H$ and $e_H$.)

For the purposes of this limited discussion a
{\em subgraph} of a rooted $H$ is a subgraph
(in the ordinary sense)
with the same roots.
We say $H$ is {\em balanced} if $\rho(H')\leq \rho(H)$
for all subgraphs $H'$ of $H$ with $v_{H'}\neq 0$
and
{\em strictly balanced} if the inequality is strict
whenever
$E'(H')\neq E'(H)$.

A {\em copy} of a rooted graph $H$ in a graph $G$
is an injection $\vp:V(H)\ra V(G)$ such that
$\vp(u)\vp(v)\in E(G)$ $\forall uv\in E'(H)$.
(Note that
here, for once, we are not assuming $G=G_{n,p}$;
and that when we do assume this below (that is, until the end of this section)
we are not placing any restriction on $p$.)

We use $\Phi(H,G)$
for the set of copies $\vp$ of $H$ in $G$.
If $(u_1\dots u_s)$ is the root sequence of $H$ and
$x_1\dots x_s$ are
distinct vertices of $G$, then we set
\[
\Phi(H,G;x_1\dots x_s)=\{\vp\in \Phi(H,G):
\vp(u_i)=x_i
~\forall i\in [s]\}
\]
and
\[
N(H,G;x_1\dots x_s)=|\Phi(H,G;x_1\dots x_s)|.
\]
(If $x_1\dots x_n$ are not all distinct then we set
$
N(H,G;x_1\dots x_s)=0.
$)

\medskip
We now take $G=G_{n,p}$.  Then
$
N(H,G;x_1\dots x_s)
$
is 
a random variable of the type treated in
Lemmas \ref{VuL1} and \ref{VuL2};
namely, with
$X=\C{[n]}{2}\sm \C{\{x_1\dots x_s\}}{2}$ (so $Y= G_{n,p}\cap X$), $d=e_H$ and
\[\h = \{\vp(E'):\vp\in \Phi(H,K_n;x_1\dots x_s)\},\]
we have
\beq{xipoly'}
\xi := N(H,G;x_1\dots x_s)
=\ww|\{A\in \h,A\sub Y\}|,
\enq
where $\ww$ is the number of {\em automorphisms
of the rooted graph} $H$ (that is,
permutations of $V(H)$ that are automorphisms in the usual sense
and fix all roots).
Of course if we set $\ww_A=\ww$ $\forall A\in\h$, then
\eqref{xipoly'} agrees with \eqref{xipoly}.

Notice that with these definitions we have
\beq{E0}
\E_0=\E\xi ~\left\{\begin{array}{c}
<\\\sim
\end{array}\right\}~
n^{v_H}p^{e_H}
\enq
and, for any $L\sub X$,
\beq{EL}
\E_L < (v_H)_{v_L}n^{v_H-v_L}p^{e_H-e_L},
\enq
with $v_L=|V(L)|$ (where $V(L)$ is the
set of vertices of $[n]\sm \{x_1\dots x_s\}$
incident with edges of $L$),
$e_L=|L|$ and, as usual,
$(j)_i=j(j-1)\cdots (j-i+1)$.
(The unnecessarily precise $(v_H)_{v_L}$ bounds the number of possibilities for a
bijection
from some $v_L$-subset of $V(H)\sm R$ to $V(L)$.)

Notice that
$v_L=v(H_L)$ and $e_L=e(H_L)$, where $H_L$ is the rooted
graph with vertex set $\{x_1\dots x_s\}\cup V(L)$,
edge set $L$ and root sequence $(x_1\dots x_s)$,
and that
$\E_L= 0$ unless
\beq{HL}
\mbox{$H_L$ is isomorphic (as a rooted graph)
to some subgraph of $H$.}
\enq
On the other hand, for $L$ with $e_L<e_H$ satisfying \eqref{HL},
\beq{vL}
v_L\geq \left\{\begin{array}{ll}
v_He_L/e_H &\mbox{if $H$ is balanced,}\\
v_He_L/e_H+\varrho_H &\mbox{if $H$ is strictly balanced,}
\end{array}\right.
\enq
where $\varrho_H$ is some positive constant (depending on $H$);
thus, recalling \eqref{EL} and
writing $\zz$ for the constant $(v_H)_{v_L}$ appearing there,
we have (for such $L$)
\beq{EL''}
E_L <
\left\{\begin{array}{ll}
\zz (n^{v_H}p^{e_H})^{1-e_L/e_H}&\mbox{if $H$ is balanced,}\\
\zz (n^{v_H}p^{e_H})^{1-e_L/e_H}n^{-\varrho{_H}}
&\mbox{if $H$ is strictly balanced.}
\end{array}\right.
\enq

\medskip
In the next two propositions we assume the above setup
and mainly aim for statements that suffice for our purposes.
The quantity $\E_0$ is, of course, the same
for all choices of $x_1\dots x_s$ and we now use it for
this common value.
The propositions and the corollary that follows them
are trivial when the $x_i$'s in question are not all distinct,
and the proofs accordingly ignore this possibility.

\begin{prop}\label{KVcor1}
If H is balanced, then for each $\theta>0$ there is a K
such that if $S>K\log n$ and
$\E_0 < n^{-\theta}S$,
then w.h.p.
\[
N(H,G;x_1\dots x_s)
< \theta S/\log n   ~~~~\forall x_1\dots x_s\in [n].
\]
\end{prop}

\nin
{\em Proof.}
It is enough to show that, for any $M$, we can choose $K$
so that (for any $x_1\dots x_s$)
\beq{KVcor1toshow}
\Pr(N(H,G;x_1\dots x_s)
\geq  \theta S/\log n) < n^{-M}.
\enq
To see this, we fix $x_1\dots x_s$, set
$\xi=N(H,G;x_1\dots x_s)$ and follow the notation introduced above.
We first claim that \eqref{KVcor1toshow} will follow if we show
\beq{ETSeps}
\mbox{there is a fixed $\eps>0$ such that}
~~\E_{l} < n^{-\eps}S ~~~\forall 0\leq l< e_H .
\enq
To see that this is enough,
suppose first that $S\geq n^{\eps/2}$.
We may then apply Lemma~\ref{VuL1} with, for example,
$B = n^{-\eps/4}S$
to say that with probability at least $1-n^{-M}$
for any fixed $M$,
\[\xi \leq B =o(S/\log n).\]
If, on the other hand, $S< n^{\eps/2}$, then
Lemma~\ref{VuL2} gives
\[
\Pr(\xi >N) < n^{-M}
\]
for a suitable $N$, and taking $K=N/\theta$ gives
\eqref{KVcor1toshow}.

\medskip
Finally, for the proof of \eqref{ETSeps}, we have, using
\eqref{EL''}, \eqref{E0}
and our hypotheses
(with $\zz$ as in \eqref{EL''} and $L$ of size less than $e_H$),
\begin{eqnarray*}
E_L
&< & (1+o(1))\zz (\E_0)^{1-e_L/e_H}\\
&<& (1+o(1))\zz (n^{-\theta}S)^{1-e_L/e_H}
< n^{-\theta(1-e_L/e_H)}S,
\end{eqnarray*}
which gives \eqref{ETSeps} with $\eps = \theta/e_H$.\qed

\begin{prop}\label{PExt}
If $H$ is strictly balanced, then for any
$\gb>0$ and M there is a K such that
for any $x_1\dots x_s\in [n]$, with probability at least
$1-n^{-M}$,
\beq{NHG}
N(H,G;x_1\dots x_s) < \left\{\begin{array}{ll}
K & \mbox{if $\E_0< n^{-\gb}$,}\\
\max\{(1+\gb)\E_0,K\log n\} &\mbox{otherwise.}
\end{array}\right.
\enq
In particular, w.h.p. \eqref{NHG} holds
for all $x_1\dots x_s\in [n]$.
\end{prop}

\nin
{\em Proof.}
We will apply one of
Lemmas \ref{VuL1}, \ref{VuL2} with
$\xi=N(H,G;x_1\dots x_s)$, $N=\C{n}{2}-\C{s}{2}$,
$d=e_H$, $b=\ww$ ($\ww$ as in \eqref{xipoly'}) and
$\eps =\frac13\min\{\gb,\vr_H\}$.

Suppose first that $\E_0\leq 1$ and let $K$ be the $J$ of
Lemma~\ref{VuL1}.
By \eqref{E0} and \eqref{EL''} we have
\beq{Elleqn}
\E_l\leq (1+o(1))\zz n^{-\vr_H}\leq N^{-\eps}
\enq
for all $ 0<l<d$;
so may take $B$ in Lemma~\ref{VuL1} to be $K \log n$,
and then the lemma gives $\Pr(\xi> K\log n)<N^{-M}$ ($ < n^{-M}$)
as desired.
If, {\em a fortiori}, $\E_0<n^{-\gb}$, then we also have
\eqref{Elleqn} for $l=0$, so with
$K$ equal to the $J$ of Lemma~\ref{VuL2}, that lemma gives
$\Pr(\xi>K)< N^{-M}$.

Finally, if $\E_0>1$, then \eqref{E0} and \eqref{EL''} give
$\E_l< (1+o(1))\zz n^{-\vr_H}\E_0\leq N^{-\eps}\E_0$ for $l\in [d-1]$;
so, with $K$ again the $J$ of Lemma~\ref{VuL1} and
$B=\max\{(1+\gb)\E_0,K\log n\}$, the lemma again
gives the relevant bound in \eqref{NHG}.\qed

In fact all our applications of
Proposition~\ref{PExt} will be instances of the next assertion
(so we really only use the proposition with $H=K_r$).

\begin{cor}\label{KrCor}
For all $s<r$,
$\gb>0$ and M there is a K such that,
with
$Z=n^{r-s}p^{\C{r}{2}-\C{s}{2}}$:
for any $x_1\dots x_s\in [n]$, with probability
at least $1-n^{-M}$,
\beq{NHG'}
\kappa(x_1\ldots x_s) < \left\{\begin{array}{ll}
K & \mbox{if $Z< n^{-\gb}$,}\\
\max\{\tfrac{(1+\gb)}{(r-s)!}Z,K\log n\}
&\mbox{otherwise}
\end{array}\right.
\enq
(where $\kappa =\kappa_G$).
In particular, w.h.p. \eqref{NHG'} holds
for all $x_1\dots x_s\in [n]$.
%
%
\end{cor}

\nin
{\em Proof.}
It is easy to see
that all rooted versions of $K_r$ are strictly balanced.
Note also that, again
taking $\xi =N(K_r,G;x_1\dots x_s)$
(for some specified choice of roots $u_1\dots u_s$ for $K_r$), we have
$\kappa(x_1\cdots x_s) = \xi/(r-s)!$
and $\E_0:=\E\xi< Z$
(see \eqref{E0}).
The assertion thus follows from
Proposition~\ref{PExt}.\qed

\subsection{Harris}\label{SecHarris}

Before continuing we quickly recall the seminal
correlation inequality of T.E. Harris \cite{Harris}.
Fix a set $I$ and set $\gO = \{0,1\}^I\equiv 2^I$
(where we make the usual identification of a set with its
indicator).
For $f:\gO\ra \Re$,
recall that
$f$ is {\em increasing in} $J\sub I$
if $f(x)\geq f(y)$ whenever $i\in J$, $x_i\geq y_i$
and $x_j=y_j$ for $j\neq i$
({\em decreasing in J} is defined similarly),
and is {\em determined by J} if $f(x)=f(y)$
whenever $x_i=y_i ~\forall i\in J$.
An event (i.e. subset of $\gO$)
$F$ is increasing in $J$ if its indicator is, and similarly
for ``decreasing in" and ``determined by."

Harris' Inequality (for Bernoullis)
says that, with expectations taken with respect to some product measure on $\gO$,
if $f$ and $g$ are increasing (i.e. in I),
then $f$ and $g$ are positively correlated
(that is, $\E fg\geq \E f \E g$), while
if one of $f,g$ is increasing and the other is
decreasing then they are negatively correlated.
Though this will be used in the proof of Theorem~\ref{TRW},
it is familiar enough that a formal statement seems unnecessary;
but we do record the following, perhaps
less familiar variant, for use in the crucial
applications of Harris' Inequality in the proof of Lemma~\ref{L2.3}
(see Section \ref{RandC}).

\begin{thm}\label{Harris}
Suppose $\xi_i$, $i\in I$, are independent Bernoullis
and $f,g:\gO\ra\Re$
with $f$ decreasing in and determined by
$J\sub I$ and
$g $ increasing in $J$.
Then $f$ and $g$ are negatively correlated.
\end{thm}

\nin
To get this from Harris' Inequality as given above,
set $\xi =(\xi_i:i\in J)$ and $\gl =(\xi_i:i\in I\sm J)$,
write $f(\xi)$ for the common value of $f(\xi,\gl)$
and set $g_\gl(\xi)=g(\xi,\gl) $.
Then
\[
\E fg = \E_\gl\E_\xi f(\xi)g_\gl(\xi)
\leq \E_\gl [\E_\xi f(\xi)\E_\xi g_\gl(\xi)]
= \E_\xi f(\xi)\E_\gl  \E_\xi g_\gl(\xi)
=\E f \E g,
\]
where the inequality follows from Harris since,
given $\gl$, the functions $f(\xi)$ and $ g_{\gl}(\xi)$
are decreasing and increasing
(respectively).

\subsection{Lower tails}\label{SJRW}

\mn
We will make substantial use of
the following fundamental lower tail bound of
Svante Janson (\cite{Janson} or
\cite[Theorem 2.14]{JLR}), for which
we need a little notation.
Suppose $A_1\dots A_m$ are subsets of the
finite set $\gG$.
For $j\in [m]$, let $I_j$ be the indicator of the event
$\{\gG_p\supseteq A_j\}$, and set $X=\sum I_j$,
$\mu = \E X =\sum_j\E I_j$
and
\beq{Delta}
\mbox{$\ov{\gD} = \sum\sum\{\E I_iI_j: A_i\cap A_j\neq\0\}.$}
\enq
(Note this includes
diagonal terms.)

\begin{thm}\label{TJanson}
With notation as above, for any $t\in [0,\mu]$,
\[\Pr(X\leq \mu -t) \leq \exp[-t^2/(2\ov{\gD})].\]
\end{thm}

\medskip
A surprising recent result of O. Riordan and L. Warnke
\cite{RW} shows that Theorem~\ref{TJanson} continues
to hold when
the events $\{\gG_p\supseteq A_i\}$ are replaced by
members of
some union- and
intersection-closed
family $\I$ of events (in some probability space)
satisfying
\beq{harris}
\Pr(B\cap C)\geq \Pr(B)\Pr(C) ~~\forall B,C\in \I,
\enq
and ``$A_i\cap A_j\neq \0$" in \eqref{Delta} is replaced
by dependence of the corresponding events.
(For Theorem~\ref{TJanson}---which
really applies to general product measures on $2^\gG$---
$\I$ is the family of increasing events and
\eqref{harris} is Harris' Inequality.)
One crucial ingredient in the proof of Lemma~\ref{L2.2}
(see Section \ref{PL2.2B}) will be an
application of a further generalization,
which we state only for the Harris context (but see Remark 2 below).

Consider some product probability measure on $2^\gG$, and
suppose $\xxx_{ij} \sub 2^\gG$ are increasing and
$\xxx_i=\cup_j\xxx_{ij}$.
Write $(i,j)\sim (k,l)$ if $\xxx_{ij}$ and $\xxx_{kl}$
are dependent.
(Note that, unlike \cite{RW},
we take $(i,j)\sim (i,j)$.)
Let $I_{ij}$ and $I_i$ be
the indicators of $\xxx_{ij}$ and $\xxx_i$ and
set $X=\sum I_i$,
\[\mbox{$\mu = \sum_{i,j}\E I_{ij}$,}\]
\[\mbox{$\ov{\Theta}=
\sum_{i,j}\sum_k\Pr(B_{ij}\cap(\cup_l\{\xxx_{kl}:
(k,l)\sim (i,j)\}))$,}\]
\[\mbox{$\ov{\gD} = \sum\sum\{\E I_{ij}I_{kl}:
(i,j)\sim (k,l)\}~~~$ ($\geq \ov{\Theta}$)} \]
and
\[\mbox{$\gc =\sum_i\sum_{\{j,k\}}\E I_{ij}I_{ik}$,}\]
with the inner sum over (unordered) pairs with $j\neq k$.

Specializing the next statement to
the case when there is just one $j$ for each $i$
yields the result of \cite{RW}.

\begin{thm}\label{TRW}
With notation as above, for any $t\in [\gc,\mu]$,
\begin{eqnarray}\label{JbdA}
\Pr(X\leq \mu -t)
&\leq &\exp[-(t-\gc)^2/(2\ov{\Theta})] \nonumber\\
&\leq &\exp[-(t-\gc)^2/(2\ov{\gD})].
\end{eqnarray}
\end{thm}

\nin
{\em Remarks.}
1.  We could, of course, replace $\mu$ in \eqref{JbdA}
by $\E X$, yielding a more natural, if slightly weaker
statement.
We will find the theorem useful
when (roughly speaking)
$\Pr(\xxx_i)\approx \sum_j\Pr(\xxx_{ij})$; that is,
when the probability of seeing at least two $\xxx_{ij}$'s
for a given $i$ is small relative to the probability
of seeing just one.  In this case there is not much difference
between $\ov{\Theta}$ and $\ov{\gD}$, and in fact the main
reason for bothering with $\ov{\Theta}$ here is that it is needed
in the proof.

\mn
2.  As noted above, Theorem~\ref{TRW}
is actually valid in the same generality as
\cite{RW}---that is, with $B_{ij}$'s from some $\I$
as in the paragraph containing \eqref{harris}---this extension
requiring only formal changes in the proof
(\eqref{harris} in place of Harris and use of
a nice observation from \cite{RW} to give the independence of
$I_{ij}$ and $Z_{ij}$ below).
As in Theorem~\ref{TJanson} and \cite{RW},
the bound in the first line of
\eqref{JbdA} may be replaced
by the slightly smaller
$\exp[-\varphi((\gc-t)/\mu)\mu^2/\ov{\Theta}]$,
where $\varphi(x) = (1+x)\log (1+x)-x$.

\mn
{\em Proof of Theorem~\ref{TRW}.}
This is mostly as in
\cite{Janson} (again, see \cite{JLR}) and \cite{RW},
so we aim to be brief.
(We are basically copying the proof of Theorem 2.14
on pp. 32-33 of \cite{JLR}, adding one nice idea
(\eqref{RWobs} below) from \cite{RW} and
taking account of the extra terms corresponding to $\gc$.)

\medskip
Let
$I_{ijk}$ and $J_{ijk}$ be
the indicators of
$
\cup\{\xxx_{kl}:(k,l)\sim (i,j)\}$ and
$
\cup\{\xxx_{kl}:(k,l)\not\sim (i,j)\}$
(so $I_k\leq I_{ijk}+J_{ijk}$ for any $i,j,k$),
and set
$Y_{ij} = \sum_kI_{ijk}$ and
$Z_{ij} = \sum_kJ_{ijk}$.
Note that $I_{ij}$ and $Z_{ij}$ are independent
(since increasing events are independent---that is,
Harris' Inequality holds with equality---iff
they depend on disjoint subsets of $\gG$).

Set $\Psi(s) = \E e^{-sX}$ ($s\geq 0$).  The main point
is to give a lower bound on
\[
-\Psi'(s) =\E Xe^{-sX} = \sum\E I_ie^{-sX}.
\]
Using $I_i\geq \sum_j I_{ij} - \sum_{\{j,k\}}I_{ij}I_{ik}$
and $X\leq Y_{ij}+Z_{ij}$ (for any $i,j$), we have
\[
\E I_ie^{-sX}
\geq \sum_j \E [I_{ij}e^{-sY_{ij}}e^{-sZ_{ij}} ]
 - \sum_{\{j,k\}}\E [I_{ij}I_{ik}e^{-sX}].
\]

The key observation from \cite{RW}
(adapted to our setting) is
\begin{eqnarray}
E [I_{ij}e^{-sY_{ij}}e^{-sZ_{ij}} ]
&=&
\E [I_{ij}e^{-sZ_{ij}}] -
\E [e^{-sZ_{ij}}(1-e^{-sY_{ij}})I_{ij}]\nonumber\\
&\geq &
\E I_{ij}\E e^{-sZ_{ij}} -
\E e^{-sZ_{ij}}\E[(1-e^{-sY_{ij}})I_{ij}]\nonumber\\
&=&
\E [I_{ij}e^{-sY_{ij}}]\E e^{-sZ_{ij}}\nonumber\\
&\geq &
\E [I_{ij}e^{-sY_{ij}}]\E e^{-sX},
\label{RWobs}
\end{eqnarray}
where the first inequality follows from the independence
of $I_{ij}$ and $Z_{ij}$
(which gives
$\E [I_{ij}e^{-sZ_{ij}}]=\E I_{ij}\E e^{-sZ_{ij}}$)
together with
Harris' Inequality (and the observation that
$f:= e^{-sZ_{ij}}$ and $g:=(1-e^{-sY_{ij}})I_{ij}$ are, respectively,
decreasing and increasing).

On the other hand, again using Harris, we have
\[\E [I_{ij}I_{ik}e^{-sX}]\leq
\E [I_{ij}I_{ik}]\E e^{-sX}.\]

Combining the preceding observations gives
\begin{eqnarray}
-(\log \Psi(s))' =\frac{-\Psi'(s)}{\Psi(s)}
&\geq &
\sum_{i,j}\E [I_{ij}e^{-sY_{ij}}]
- \sum_i\sum_{\{j,k\}}\E [I_{ij}I_{ik}]\label{logderiv} \\
&\geq &\mu e^{-s\ov{\Theta}/\mu} -\gc.\nonumber
\end{eqnarray}
The lower bound $\mu\exp[-s\ov{\Theta}/\mu]$
on the first sum in \eqref{logderiv}
is obtained {\em via}
two applications of Jensen's Inequality as
in the last four lines of \cite[p. 32]{JLR}.

We then have
\[
-\log\Psi(s) \geq \int_0^s(\mu e^{-u\ov{\Theta}/\mu}-\gc)du
= \frac{\mu^2}{\ov{\Theta}}(1-e^{-s\ov{\Theta}/\mu}) -s\gc,
\]
yielding (by Markov's Inequality)
\[
\log\Pr(X\leq \mu-t)\leq \log \E e^{-sX} +s(\mu-t)
\leq
-\frac{\mu^2}{\ov{\Theta}}(1-e^{-s\ov{\Theta}/\mu}) +s(\mu-(t-\gc)),
\]
and applying this with
$s=-\log (1-(t-\gc)/\mu)\mu/\ov{\Theta}$ gives
\eqref{JbdA} (actually the slightly better bound mentioned in Remark 2
above; again, {\em cf.} \cite[p. 33]{JLR}).\qed

\subsection{A calculation}

The following observation will be needed
twice below
(see the proofs of Lemmas \ref{Rlemma} and \ref{basic'}), so
we include it in these preliminaries.

\begin{lemma}\label{newHlemma}
For each $\xi>0$ there is a $\vt>0$ so that the following is true
(as usual, provided $p$ is as in \eqref{p} with a large enough $C$).
Let $R\sub \C{V}{2}$ satisfy
\beq{Rdeg}
\gD_R < \vt np/\log n
\enq
and let $\h$ consist of all
sets of the form $K(xy,Z) := \C{\{x,y\}\cup Z}{2}\sm \{xy\}$
with $xy\in R$ and $Z\in \C{V\sm \{x,y\}}{r-2}$.
Let $I_K$ be the indicator of $\{K\sub G\}$
and
\beq{ovgD}
\ov{\gD}=\sum\sum \{\E I_KI_L:K,L\in \h, K\cap L\neq \0\}.
\enq
Then
\beq{RLdelta}
\ov{\gD} < \xi|R|\glr^2/\log n.
\enq
\end{lemma}

\mn
{\em Remark.}  The $\h$'s in our applications
will be subsets of the one here, which, of course, only
shrinks $\ov{\gD}$.

\mn
{\em Proof.}
For $K=K(xy,Z)$ as in the lemma, let $V(K) = \{x,y\}\cup Z$
and $e_K=\{x,y\}$.
We organize
$\g:=\{(K,L):K,L\in \h, K\cap L\neq \0\}$
as follows.
For $(K,L)\in \g$ set $a(K,L) = |e_K\cap e_L| $ ($\in \{0,1,2\}$)
and
$b(K,L)= |V(K)\cap V(L)|$
($\in \{2\dots r\}$).
Notice that if $a(K,L) = 2$ (that is, if $e_K=e_L$),
then $K\cap L\neq\0$ implies $b(K,L)\geq 3$.

Let
\[N(a,b)= |\{(K,L)\in \g:
a(K,L) =a, b(K,L) =b\}|.\]
Then, since $|V(K)\cup V(L)|= 2r-b(K,L)$,
there is a fixe $B=B_r$ (e.g., crudely, $B= r!$) such that
\[
N(0,b) < B|R|^2n^{2r-4-b},
\]
\[\mbox{$N(1,b) < B\sum_xd^2_R(x)n^{2r-3-b}\leq B|R|\gD_{_R} n^{2r-3-b}$}\]
and
\[N(2,b) < |R|n^{2r-2-b}.\]

\mn
On the other hand,
\[
\E I_KI_L = p^{r^2-r-2-|K\cap L|}
\]
and
\[
|K\cap L| ~\left\{\begin{array}{ll}
= \Cc{b(K,L)}{2}-1&\mbox{if $a(K,L)=2$,}\\
\leq \Cc{b(K,L)}{2}&\mbox{otherwise.}
\end{array}\right.
\]
Combining these observations we have

\begin{eqnarray*}
\ov{\gD} &\leq &
\sum_{b=3}^rN(2,b)p^{r^2-r-1-\Cc{b}{2}} +
B \sum_{a=0}^1\sum_{b=2}^rN(a,b)p^{r^2-r-2-\Cc{b}{2}} \\
&<&
|R|n^{2r-4}p^{r^2-r-2}\left[\sum_{b=3}^rn^{2-b}p^{1-\Cc{b}{2}} +
B(|R|+\gD_{_R} n)\sum_{b=2}^rn^{-b}p^{-\Cc{b}{2}}\right],
\end{eqnarray*}
so for \eqref{RLdelta} would like to say that the
expression in square brackets is less than
$\xi \log^{-1} n$,
which a little checking---using \eqref{Rdeg}
with a small enough $\vt$ (something like $\xi/B$ will do)
and our lower bound on $p$---shows
to be true.
(The largest contributions are (i) $n^{2-r}p^{1-\C{r}{2}}=\glr^{-1}$
corresponding to $b=r$ in the first sum, and
(ii)  B$\gD_{_R}n^{-1}p^{-1}$, which is the $(b=2)$-term from
the second sum multiplied by $B\gD_{_R} n$.)\qqqed

\subsection{Miscellaneous}

In closing these preliminaries we mention two last, easy points.
First, we recall just one detail
(borrowed from \cite[Lemma 4.1]{JKV}), of the connection between
$G_{n,p}$ and $G_{n,M}$:
\begin{lemma} \label{lemma:twomodels}
Let $n^{\Omega(1)}= M \le {n \choose 2}$ be an integer
and $p= M/{n \choose 2}$.
If an event $\eee$
holds with probability at least $1-\eps$ in $G_{n,p}$,
then it  holds with probability at least $1- O(n\eps)$ in $G_{n,M}$.
\end{lemma}
\nin
(The extra $n$ in the conclusion won't
be a problem, since
our exceptional probabilities will be much smaller than $1/n$.
We will also want something in the other direction, but
defer the trivial statement until needed; see \eqref{transfer2'}.)

\medskip
We will also find a couple of uses for the following observation
\cite{deWerra,McDiarmid}, in which we call a coloring
{\em equitable} if the sizes of the color
classes differ by at most one.

\begin{prop}\label{PCol}
For any $m\geq \gD+1$,
the edges of any simple graph of maximum degree at most
$\gD$ can be equitably colored with $m$ colors.
\end{prop}
\nin

\section{Main points for the proof of Lemma~\ref{L2.2}}\label{S2.2}

Here we derive Lemma~\ref{L2.2} from the following three
assertions, which are proved in Sections \ref{PLxKLS}-\ref{PL2.2B}.
%
We use
$\tau(A_1\dots A_{r-1})$ for the number of choices of distinct
$x_1\dots x_{r-1}$ with $x_i\in A_i$ and
all pairs from $\{x_1\dots x_{r-1}\}$ belonging to $G$,
and also
write
$\tau(A,B,C)$ for $\tau(A,B,C\dots C)$ (with $r-3$ copies of $C$).

\begin{lemma} \label{xKLS}
For fixed $\theta, \vr>0$,
w.h.p.
\beq{tRST}
\tau(S_1\dots S_{r-1}) >
(1-\vr) |S_1|\cdots |S_{r-1}|p^{\Cc{r-1}{2}}
\enq
whenever $v\in V$ and $S_1\dots S_{r-1}$ are disjoint subsets
of $N_v$ with each of $|S_2|\dots |S_{r-1}|$
at least $\theta np$ and
\beq{tlb}
|S_1|> \vr^{-2}6r\log r\cdot2^r
\max\left\{\frac{1}{\theta p}~,\frac{np}{\theta^{r-2}\glr}\right\}.
\enq
\end{lemma}

\mn
{\em Remarks.}
For the proof of Lemma~\ref{L2.2} we could
replace
the lower bound in
\eqref{tlb} by $\theta np$,
but the present stronger version
(with the weaker lower bound on $|S_1|$)
will be needed in the
proof of Lemma~\ref{L2.3}.
The constants in \eqref{tlb} (i.e. the $\theta$'s and the expression preceding the ``max") are unimportant.

\begin{lemma}\label{xKLS'}
For fixed $\pi\geq\eps >0$,
w.h.p.
\beq{TRST}
\tau(S,T,R) \leq 8\pi |T| \glr
\enq
whenever $v\in V$; $S$ and $T$ are disjoint subsets
of $N_v$ with $|T| < \eps np$ and
$|S|< \pi |T|/\eps$; and $R\sub N_v\sm (S\cup T)$.
\end{lemma}

\nin
(Of course there is no change in content if we say
``$R= N_v\sm (S\cup T)$," but
the stated version will be convenient.)

\begin{lemma}\label{L2.2B}
For each $\pi>0$ there is an $\eps>0$ such that w.h.p.
\beq{kappaST}
\kappa(S,T) < \pi|S|\glr
\enq
whenever
$S\sub \Cc{V}{2} $ and $T\sub G$ satisfy
$\gD_S\leq 2np,$ $V(S)\cap V(T) = \0$
and
$|T|\leq |S|<\eps n^2p$.
\end{lemma}

\mn
{\em Remarks.}
Each of Lemmas~\ref{xKLS'} and \ref{L2.2B} bounds
the quantity in question by something like its natural value;
namely, the r.h.s. of \eqref{TRST} is, up to scalar, the
natural value of the l.h.s. when $|T|\approx \eps np$ and
$|S|\approx \pi np$
(and $R= N_v\sm (S\cup T)$),
while
\eqref{kappaST} says that for a small
$T\sub G$, $\kappa(S,T)$ can account for only a small
fraction of $\kappa(S)\asymp|S|\glr$.
(It is easy to see that it is not enough to bound $T$ without
reference to $|S|$.)

The proof of Lemma~\ref{L2.2B} turns out to be much less straightforward than one might expect,
and a small puzzle may  be worth mentioning.
The bound on $\gD_S$, which will eventually come for free
because we will have $S\sub G$, happens to be just
what's needed for the current proof of the lemma,
but we don't know that it is really necessary.
When $r=4$ the lemma can fail with $\gD_S$ as small as $n^2p^3$
(note that by \eqref{p}, $\gD_S\approx n^2p^3$ requires $r\leq 4$ since $\gD_S< n$): for some $x\in V$,
take $T=G[N_G(x)]$ and let $S$ consist of $n^2p^3$ pairs
containing $x$ and avoiding $N_G(x)$;
then $|S|\ll n^2p$ and
(typically) $|T| \approx n^2p^3/2 <|S|$, while $\kappa(S,T)\approx |S||T|p^2\approx |S|\glr/2$.
But we don't know whether
$\gD_S\ll n^2p^3$ suffices when $r=4$ or
whether {\em any} bound (on $\gD_S$) is needed for larger $r$.
It would be interesting to understand what's going on here,
and so see whether this seemingly innocuous lemma can be proved less circuitously.

\medskip
We will also need the fact, contained in Proposition~\ref{propG[X]}(a),
that for some fixed M, w.h.p.
\beq{G[X]}
|G[X]| < \max\{|X|^2p,M|X|\log n\} ~~\forall X\sub V.
\enq

\mn
{\em Proof of Lemma} \ref{L2.2}.
Recalling $\gc$ from the definition of $Q_G(\Pi)$ preceding Lemma~\ref{L2.2},
choose constants $\vr,\gz,\eta>0$ with $\gc^{r-2}\gg \vr\gg\gz$
and $\eta\ll \gz^2$
small enough so that the conclusion
of Lemma~\ref{L2.2B} holds when $\pi =\gz$ and $\eps = (r-1)\eta$.
We assert that Lemma~\ref{L2.2} holds with this value of $\eta$.

\medskip
What we actually show is that the ``w.h.p." statement in
Lemma~\ref{L2.2} is true provided \eqref{Deg},
\eqref{G[X]}
and the conclusions of
Lemmas \ref{xKLS}-\ref{L2.2B} hold for suitable
values of the parameters.
To say this properly, define properties:

\mn
(A)  \eqref{tRST} holds for $\theta = (2\vr/(1-\vr))^{1/(r-2)}$
and all $v,S_1\dots S_{r-1}$
as in Lemma~\ref{xKLS};

\mn
(B)  \eqref{TRST} holds whenever
$(\eps,\pi)\in \{(\gz,\vr/(8(r-2))),
(\theta,(\gc-2\vr)/(8(r-2)))\}$ (with $\theta$ as in (A))
and $v,S,T,R$ are
as in Lemma~\ref{xKLS'};

\mn
(C)  \eqref{kappaST} holds for $(\pi,\eps) = (\gz,(r-1)\eta)$
and $S,T$
as in Lemma~\ref{L2.2B};

\mn
(D)  \eqref{G[X]} holds (for some fixed $M$, whose value will
be irrelevant here).

\mn
By Lemmas \ref{xKLS}-\ref{L2.2B} and Propositions \ref{vdegree} and \ref{propG[X]}(a),
it is enough to show that
the conclusion of Lemma~\ref{L2.2} holds whenever \eqref{Deg}
(we just need degrees bounded by $2np$) and
(A)-(D) are satisfied, which we now assume.

\medskip
Fix $F$ and $\Pi$ as in Lemma~\ref{L2.2} and set
$I=F[A_1]$ and $L = \Pi_G\sm F$.
We should show (provided $I\neq \0$)
\beq{L3I}
|L|>(r-1)|I|,
\enq
so assume \eqref{L3I} fails and aim for a
contradiction.

Set
$X =\{x: d_I(x)>\gz np\}$
($\sub A_1$) and $Y=A_1\sm X$.
We begin with the observation that not many edges of $I$
lie in $X$:
noting that $|I|> |X|\gz np/2$ and
(consequently) $|X|< 2(\eta/\gz)n$,
and using (D),
we have
\beq{GX}
|G[X]|
< \tfrac{2}{\gz}\max\{\tfrac{|X|}{n},\tfrac{M\log n}{np}\}|I|
< 4 \eta \gz^{-2}|I|.
\enq

\medskip
We now use $K$ for members of $\K(I,A_2\dots A_{r-1})$
(recall this was defined near the end of Section~\ref{Usage}).
Say such a $K$
is {\em covered at} $v\in A_1$ if it contains
an edge of $L(v):=\{e\in L:v\in e\}$
(so in particular contains $v$),
is {\em covered at} $W\sub A_1$ if it is covered at some $v\in W$,
and is
{\em covered off} $A_1$ if
it contains an edge of $L\cap \nabla(A_2\dots A_{r-1})$.
Let $I_1= \cup_{y\in Y}I_1(y)$, where (for $y\in Y$)
\[
I_1(y) =\{yw\in I: \mbox{at least
$\vr \glr$ of the $K$'s on $yw$
are covered at $y$}\};\]
\[I_2=\{e\in I: \mbox{at least $\vr \glr$ of the $K$'s on $e$
are covered off $A_1$}\};\]
and
$I_3= I\sm (G[X]\cup I_1\cup I_2)$.
Note that
each $e\in I_3$ is of the form $xy$ with $x\in X$, $y\in Y$
and at least $(\gc-2\vr)$ $K$'s from $\K(e,A_2\dots A_{r-1})$
covered at $x$.

In what follows we show that each of $I_1$, $I_2$, $I_3$ is
small compared to $I$, which with \eqref{GX} gives
the desired contradiction.
We first assert that
\beq{I1I}
|I_1| \leq 8(r-1)(r-2)\gz |I|/\vr,
\enq
which follows from
\beq{I1I2}
|L(y)| \geq \vr |I_1(y)|/(8(r-2)\gz) ~~\forall y\in Y
\enq
and our assumption that \eqref{L3I} fails.
For \eqref{I1I2}, notice that
if $|N_L(y)\cap A_i|< \vr |I_1(y)|/(8(r-2)\gz)$
for $i=2\dots r-1$, then (B),
applied, for each $i\in [2,r-1]$,
with $(\eps,\pi)=(\gz,\vr/(8(r-2))$, $v=y$
and
\[(S,T,R) = ( N_L(y)\cap A_i,N_{I_1}(y),N_y\sm (A_1\cup A_i))\]
(note $|I_1(y)|<\gz np$)
bounds
the number of $K$'s from $\K(I_1(y),A_2\dots A_{r-1})$
that are covered at $y$ by $\vr |I_1(y)|\glr $,
contradicting the assumption that each $e\in I_1(y)$ lies in more than
$\vr \glr $ such $K$'s.

\medskip
We next show that $I_2$ is small.
Set $J=L\cap \nabla(A_2\dots A_{r-1})$.
By the definition of $I_2$ we have
$\kappa(I_2,J) \geq \vr |I_2|\glr $.
But, as we will show in a moment, (C) together with
($|J|\leq $) $|L|\leq$ $(r-1)|I|$ gives
\beq{kappaI2J}
\kappa(I_2,J) \leq \kappa(I,J)< (r-1)\gz|I|\glr ,
\enq
so that \beq{I2I}
|I_2|< (r-1)\gz|I|/\vr.
\enq
For \eqref{kappaI2J} we use (C) with
$(S,T)=(I_2,J)$ if $|J|\leq |I_2|$ and
$(S,T)=(J,I_2)$ otherwise.
In either case we have $\gD_S\leq 2np$ (by \eqref{Deg}) and
$|S|\leq (r-1)|I| < (r-1)\eta n^2p$,
so that (C) gives
$\kappa(S,T) < \gz |S|\glr  \leq (r-1)\gz|I|\glr $.

\medskip
Finally, we show $I_3$ is small. Set
$\theta=(2\vr/(1-\vr))^{1/(r-2)}$ (as in (A)).
We first observe that
\beq{gagd}
|I_3(x)| \leq \theta np ~~\forall x\in X.
\enq

\mn
To see this, suppose instead that
$|I_3(x)| =d > \theta np$
for some $x\in X$.
By \eqref{HP3} we have
\[
d\leq |N_F(x)\cap A_1|\leq \min\{|N_F(x)\cap A_i|:2\leq i\leq r-1\};
\]
so according to (A)---note $\theta np$ is larger than
the bound in \eqref{tlb}---there are at least
$(1-\vr)d^{r-1}p^{\C{r-1}{2}}$
$K$'s containing $x$ and one vertex from each of $N_{I_3}(x),
N_F(x)\cap A_2\dots N_F(x)\cap A_{r-1}$---that is,
\beq{kappax}
\kappa(x,N_{I_3}(x),
N_F(x)\cap A_2\dots N_F(x)\cap A_{r-1})\geq
(1-\vr)d^{r-1}p^{\C{r-1}{2}}
\enq
---each of which must be covered either at $Y$ or off $A_1$.
But since an edge of $I_3(x)$
is contained in at most $2\vr \glr $ $K$'s that are covered
in one of these ways, the l.h.s. of \eqref{kappax} is
at most $2d\vr \glr $, which is less than the r.h.s.
Thus \eqref{gagd} does hold.

\medskip
We may now proceed as we did in bounding $|I_1|$:  for $x\in X$,
each $e\in I_3(x)$ lies in at least $(\gc-2\vr)\glr $
$K$'s from $\K(e,A_2\dots A_{r-1})$ that are covered at $x$, which
by (B) (with
$(\eps,\pi)=(\theta,(\gc-2\vr)/8(r-2))$)
implies
\[|L(x)|\geq \max\{|N_L(x)\cap A_i|:i\in [2,r-1]\}
\geq  \tfrac{\gc -2\vr}{8(r-2)\theta}|I_3(x)|  ~~\forall x\in X
\]
and (again using failure of \eqref{L3I})
\beq{lastI}
|I_3|\leq 8(r-1)(r-2)\theta|I|/(\gc -2\vr)
\enq
(a small multiple of $|I|$ because of our choice of $\theta$.)

Finally, combining
\eqref{GX}, \eqref{I1I},
\eqref{I2I}
and \eqref{lastI}, we have the contradiction
$|I| \leq |G[X]|+|I_1|+|I_2|+|I_3|<|I|$.
\qqqed

\section{Main points for the proof of Lemma~\ref{L2.3}}\label{S2.3}

\mn
Here we just state the two
main assertions underlying Lemma~\ref{L2.3} and
show that they suffice.
The assertions themselves are proved in Sections~\ref{PXdefect} and \ref{PL2.3'},
with both arguments rooted in the observations of Section~\ref{RandC}.

\medskip
Let
\beq{abphi}
\mbox{$\aaa_r=\frac{r-4}{2(r-3)}$, $~\bbb_r=\frac{r(r-3)}{2(r-1)^2}~$
and
$\cc_r= (\aaa_r+\bbb_r)/2$.}
\enq
We can now, finally, say something about $\gc$
(the parameter in the definition
of a bad pair in the passage preceding Lemma~\ref{L2.2}).
Here it is not necessary (or desirable) to specify an actual value;
we just stipulate that
\beq{gcsat}
0<\gc < \frac{1}{2}
\left(\frac{\cc_r-\aaa_r}{4r^2+6}\right)^{r-2}.
\enq

\medskip
For $x\in V$ and disjoint $A_1\dots A_{r-1}\sub V$, set
\beq{D}
D(x;A_1\dots A_{r-1}) =
\sum\{d_{A_i}(x)d_{A_j}(x):1\leq i<j\leq r-1\}.
\enq
For a cut
$\Pi=(A_1\dots A_{r-1})$, we also write
$D_\Pi(x) $ for $D(x;A_1\dots A_{r-1})$.
We say $x$ is {\em bad for} $\Pi$ if $x\in A_1$ and
$D_\Pi(x)< \cc_r n^2p^2$.

\mn
{\em Remark.}
When $d(x)\approx np$ and
$\Pi=(A_1\dots A_{r-1})$,
$\bbb_r n^2p^2$ is essentially the minimum possible value of
$D_\Pi(x)$ if at least $r-2$ of the $d_{A_i}(x)$'s are at least
$np/(r-1)$, and $\aaa_r n^2p^2$ is (essentially) the maximum possible
value if at least two of the $d_{A_i}(x)$'s are zero.
While the inequality $\aaa_r<\bbb_r$
is easily verified,
we don't see any
intuitive reason why it should be true;
yet our proof (of Theorem~\ref{MT}) collapses if it is not.

We will use $\cc_r>\aaa_r$ in the present section
(see \eqref{Zdeg}) and $\cc_r< \bbb_r$ twice in
Section \ref{PXdefect}
(see the proof of Proposition~\ref{propt} and
\eqref{deffG}).

%
\begin{lemma}\label{Xdefect}
There is a fixed $\nu>0$ such that
w.h.p.: for every $t>0$, every balanced cut with at least t bad
vertices has defect at least
$\nu t n^{3/2}p^2$.
\end{lemma}

\medskip
Set $\gz=(2\gc)^{1/(r-2)}$---thus \eqref{gcsat} is
\beq{gzsat}
0<\gz <
\frac{\cc_r-\aaa_r}{4r^2+6}
\enq
---and
\begin{eqnarray}
\gS &= &
24r\log r\cdot2^r
\max \{(\gz p)^{-1},np(\gz^{r-2}\glr)^{-1}\}.
\nonumber\\
&=&\left\{\begin{array}{ll}
\Theta(np/\gL_r(n,p))&\mbox{if $p<n^{-2/(r+2)}$,}\\
\Theta(p^{-1})&\mbox{otherwise.}
\end{array}\right.\label{s}
\end{eqnarray}

\begin{lemma}\label{L2.3'}
W.h.p. $\deff_G(\Pi)\geq 2r^2|Q|$ whenever
$\Pi=(A_1\dots A_{r-1})$ is a balanced cut and
$Q\sub Q_G(\Pi)$
satisfies
\beq{dQx}
d_Q(x)< \gS ~~\forall x.
\enq
\end{lemma}

\mn
{\em Remarks.}
As in Lemma~\ref{L2.3} the factor $2r^2$ in the defect bound
is what's needed below, but
could actually be any constant.
The relatively severe constraint on $\gc$ in
\eqref{gcsat} is needed for the derivation of Lemma~\ref{L2.3}, not for
Lemma~\ref{L2.3'} itself.

\mn
{\em Proof of Lemma} \ref{L2.3}.
We show that the conclusion of Lemma~\ref{L2.3}
holds whenever we have:
the conclusion of Lemma~\ref{xKLS} for
$\theta =\gz$ and $\vr=1/2$;
the conclusions of Lemmas \ref{Xdefect} and \ref{L2.3'}; and
\beq{deg1}
d(x) < (1+o(1))np ~~\forall x.
\enq
This is enough since (by the lemmas just mentioned and
Proposition~\ref{vdegree}) these assumptions hold w.h.p.

Suppose we have the stated conditions and
$\Pi$, $Q$ are as in Lemma~\ref{L2.3}.
We may of course assume $Q\neq \0$.
We first show that,
for any $x\in A_1$,
\beq{dBxdCx}
\min\{d_{A_i}(x):2\leq i\leq r-1\}>\gz np
~\Ra~ d_Q(x)< \gS
\enq
and
\beq{dQ}
d_Q(x) \leq \gz np .
\enq
{\em Proof.}
Notice that $\gS$ is simply the r.h.s. of \eqref{tlb}
for $\theta =\gz$ and $\vr=1/2$.
So if \eqref{dBxdCx} fails---that is, if
$d_{A_i}(x)>\gz np$ for $i\in \{2\dots r-1\}$
and $d_Q(x)\geq \gS$---then
the conclusion of Lemma~\ref{xKLS}, applied
with $v=x$, $S_1 =N_Q(x)$ and $S_i =N_{A_i}(x)$ for $i\in [2,r-1]$,
gives
\begin{eqnarray*}
\sum_{y\in N_Q(x)}\kappa(xy,A_2\dots A_{r-1})&=&
\tau(N_Q(x),N_{A_2}(x)\dots N_{A_{r-1}}(x))\\
&\geq &(1/2)d_Q(x)(\gz np)^{r-2}p^{\Cc{r-1}{2}} \\
&=&d_Q(x)\gc \gL_r(n,p),
\end{eqnarray*}
a contradiction since
$xy\in Q$ implies $\kappa(xy,A_2\dots A_{r-1})< \gc \glr$.
This gives \eqref{dBxdCx}, and \eqref{dQ} follows easily:
if \eqref{dQ} fails then
\eqref{dQx'} implies
$\min\{d_{A_i}(x):2\leq i\leq r-1\}>\gz np$,
and then \eqref{dBxdCx} gives $d_Q(x)< \gS$;
but
$\gS <\gz np$,
so we again have a contradiction.\qed

\medskip
Let $X$ be the set of vertices that are bad for $\Pi$ (so $X\sub A_1$)
and set
\[
Z_i=\{x\in A_1:d_{A_i}(x)\leq \gz np\}  ~~~~~
2\leq i\leq r-1
\]
and $Y = A_1\sm (X\cup Z_2\cup\cdots \cup Z_{r-1})$.

Let $Q_v$ be the set of pairs from $Q$ meeting $X$.
Then $|Q_v|\leq |X|\gz np$ (by \eqref{dQ}),
which with the conclusion of
Lemma~\ref{Xdefect} gives
\[\deff_G(\Pi) \geq \nu |X|n^{3/2}p^2\geq (\nu/ \gz)n^{1/2}p|Q_v|.\]
So we may assume $|Q_v| < r\gz|Q|/(\nu n^{1/2}p)$ ($\ll |Q|$).

We may further assume that $|Q[Y]|< (2r)^{-1}|Q|$;
otherwise---in view of \eqref{dBxdCx}, which implies
$d_{Q[Y]}(x)<\gS ~\forall x$
(note $d_{Q[Y]}(x)=0$ if $x\not\in Y$)---we may apply
the conclusion of Lemma~\ref{L2.3'} to $Q[Y]$ to obtain
\[\deff_G(\Pi) \geq 2r^2|Q[Y]|\geq r |Q|.\]

We may thus assume (w.l.o.g.) that at least (say) $|Q|/r$
of the edges of $Q$ meet $Z:=Z_2\sm X$, which with \eqref{dQ}
gives
\beq{Qleq4}
|Q|\leq r|Z|\gz np.
\enq
But we will show that if this is true
then we can
obtain a cut significantly larger than $\Pi$
by moving an appropriate subset of $Z$ to $A_2$.
The main point here is that vertices of $Z$ must have
many neighbors in $A_1$.
Set $\gl=\gl_r = (\cc_r-\aaa_r-2\gz)$.
We assert that
\beq{Zdeg}
d_{A_1}(x) > \gl np ~~~\forall x\in Z.
\enq
{\em Proof.}  For $x\in Z$ we have
$D_\Pi(x)\geq \cc_r n^2p^2$ (since $x\not\in X$)
and (a little crudely, using \eqref{deg1} to bound $d(x)$
and $x\in Z_2$ to bound $d_{A_2}(x)$)
\begin{eqnarray*}
D_\Pi(x) &<& (d_{A_1}(x) +d_{A_2}(x))d(x)
+\Cc{r-3}{2}\left(\tfrac{d(x)}{r-3}\right)^2\\
&<& (1+o(1))[(d_{A_1}(x)+\gz np)np +\aaa_rn^2p^2]\\
&<& \left[\tfrac{d_{A_1}(x)}{np}
+\aaa_r +2\gz\right]n^2p^2;
\end{eqnarray*}
so we have \eqref{Zdeg}.\qqqed

\medskip
Now choose $W \sub Z$ so that $\nabla_G(W,A_1\sm W)$
contains at least half the edges of $G[A_1]$ meeting $Z$.
(This is true on average for
$W$ chosen uniformly from the subsets of $Z$,
so such a choice is possible.)
Let
\[\Pi'=(A_1\sm W,A_2\cup W,A_3\dots A_{r-1}).\]
Then
\begin{eqnarray}
|\Pi'_G|-|\Pi_G| & = & |\nabla_G(W,A_1\sm W)|-|\nabla_G(W,A_2)|
\nonumber\\
&\geq & \sum_{x\in Z}(d_{A_1}(x)/4 - d_{A_2}(x))\nonumber\\
&\geq &|Z|(\gl/4-\gz)np.\label{Ztfrac}
\end{eqnarray}
Thus $\deff_G(\Pi)$ is at least the r.h.s. of \eqref{Ztfrac},
and according to \eqref{Qleq4}
(and \eqref{gzsat}) this is larger than $r|Q|$.
This completes the proof of Lemma~\ref{L2.3}\qed

\mn
{\em Remark.}
Constants aside, the value of $\gS$ in \eqref{s} cannot be increased without
invalidating the proof of Lemma~\ref{L2.3'} (see the bound on $\ov{\gD}$ in the proof of
Lemma~\ref{basic'}),
while Lemma~\ref{xKLS} (not just its proof) is false for smaller values of
the bound in \eqref{tlb}.
But the above
proof of Lemma~\ref{L2.3}
uses Lemma~\ref{xKLS} to bound degrees in $Y$
by an instance, $\gS$, of the latter bound, supporting application of Lemma~\ref{L2.3'};
so the fact that the $\gS$ needed for this application is not less than what's
affordable in \eqref{tlb} is crucial, and it would be nice to somehow
understand that this is not just a lucky accident.

\section{Proof of Lemma~\ref{xKLS}}\label{PLxKLS}

\mn
For given disjoint $ S_1\dots S_{r-1}\sub V$
with $|S_i|=s_i$, let
\[
\B(S_1\dots S_{r-1})
= \{\tau(S_1\dots S_{r-1})
< (1-\vr)s_1\cdots s_{r-1} p^{\Cc{r-1}{2}}\}.
\]
We will show that for any $ S_1\dots S_{r-1}$ with sizes as in
Lemma~\ref{xKLS},
\beq{KLS}
\Pr(\B(S_1\dots S_{r-1})) < \exp[-(3 \log r)np],
\enq
but first show that this does give the lemma.
By Proposition~\ref{vdegree} it is enough to bound the probability that
the conclusion of Lemma~\ref{xKLS} fails at some $v$ with $d(v)\leq 2np$;
this probability is at most
\beq{dv2np}
\sum_W\Pr(N_v=W)
\sum_{S_1\dots S_{r-1}}\Pr(\B(S_1\dots S_{r-1})),
\enq
with the first
sum over $W\sub V\sm \{v\}$ of size at most $2np$
and the second over disjoint $S_1\dots S_{r-1}\sub W$ obeying
the size requirements of the lemma.
But according to
\eqref{KLS} the expression in \eqref{dv2np} is less than
\[
\Cc{n}{2np}(r-1)^{2np} \exp[-(3 \log r)np] ~ <~ \exp[-(\log r)np] ~=~ o(1/n)
\]
(which is needed since we multiply by $n$ to account for
the choice of $v$).

\mn
{\em Proof of} \eqref{KLS}.
This is a straightforward application of Theorem~\ref{TJanson}, the
notation of which we now follow.
With $\gG=\nabla(S_1\dots S_{r-1})$ and
$A_1\dots A_m$ the (edge sets of) copies of $K_{r-1}$ in $\gG$,
we have
\beq{mubd}
\mu = s_1\cdots s_{r-1} p^{\Cc{r-1}{2}}
\enq
and
\begin{eqnarray}
\ov{\gD} &<& \sum_{i=2}^{r-1}\sum_{I\in \Cc{[r-1]}{i}}
\prod_{j=1}^{r-1}s_j^2\prod_{j\in I}s_j^{-1}
\cdot p^{2\Cc{r-1}{2}-\Cc{i}{2}}.\nonumber
\end{eqnarray}
For \eqref{KLS} ({\em via} Theorem~\ref{TJanson}) we need
$(\vr \mu)^2 > 3\log r\cdot np\cdot 2\ov{\gD}$, or
(equivalently)
\beq{vr2}
\vr^2 > (6\log r) np ~
\sum_{i=2}^{r-1}
\sum_{I\in \Cc{[r-1]}{i}}
\prod_{j\in I}s_j^{-1}
\cdot p^{-\Cc{i}{2}}.
\enq
Setting $s_1^*=\min\{s_1,\theta np\}$ and using $s_i>\theta np$, we find that
the r.h.s. of \eqref{vr2} is less than
\[
(6\log r) 2^r np ~
\sum_{i=2}^{r-1}
\left[ s_1^* ~(\theta np)^{i-1}
p^{\Cc{i}{2}}\right]^{-1}
 ~~~~~~~~~~~~~~~~~~~~~~~~~~~~~~~~
\]
\[
 ~~~~~~~~~~~~~~~~~~~~~~~~~~~~~~~~
 = (6\log r) 2^r  ~
\sum_{i=2}^{r-1}
 (s_1^*)^{-1}\theta ^{1-i} n^{2-i}p^{-\Cc{i+1}{2}+2},
\]
so that the inequality holds provided
\beq{s1*}
s_1^*  > \vr^{-2}
(6\log r) 2^r  ~
\sum_{i=2}^{r-1}
\theta ^{1-i} n^{2-i}p^{-\Cc{i+1}{2}+2}.
\enq
It is also easy to see (e.g. by considering ratios
of consecutive terms) that the largest summand
in \eqref{s1*} is either the first or the last; so,
again without being too careful,
we may (upper) bound the entire
r.h.s. of \eqref{s1*} by the expression in \eqref{tlb}
(which gives the lemma since, as already noted, this expression
is less than $\theta np$).\qed

\section{Proof of Lemma~\ref{xKLS'}}\label{PLxKLS'}

(A reminder:
rooted graphs and some of the other notions and notation used
here were introduced in Section~\ref{SKV}.)

Set $\I=\{(i,j):1\leq i<j\leq r-1\}$ and
write ``$\prec$" for ``reverse lexicographic" order on
$\I$; that is, $(i,j)\prec (k,l)$ if either
$j<l$ or $j=l$ and $i<k$.
For $(i,j)\in \I$, write
$\vs(i,j)$ for the index of $(i,j)$ under ``$\prec$";
for example $\vs(2,3)=3$ and $\vs(r-2,r-1)=\C{r-1}{2}$.

For $(i,j)\in \I$,
let $H_{ij}$
be the rooted graph
with vertex set $\{u_0,u_i\dots u_j\}$, edge set
\[\{u_0u_k:k\in [j]\}\cup\{u_ku_l: (k,l)\in \I,(k,l)\prec (i,j)\}\]
(so all edges except those joining $j$ to $[i,j-1]$)
and root sequence $(u_0,u_i,u_j)$.
Set
\beq{SijDef}
S_{ij} = n^{j-1}p^{\vs(i,j)+j-1}
\enq
and notice that
\beq{Sij}
S_{ij} \geq  n^{j-1}p^{\C{j+1}{2}-1}
=(np^{(j+2)/2)})^{j-1} =
\left\{\begin{array}{ll}
\glr &\mbox{if $j=r-1$,}\\
n^{\gO(1)} &\mbox{otherwise.}
\end{array}\right.
\enq   

We need one auxiliary result:
\begin{prop} \label{VuCor}
For any $\vt>0$,  
w.h.p.
\[N(H_{ij},G;v,x,y) < \vt S_{ij}/\log n\]
for all $(i,j)\in \I$ and $v,x,y\in V$.
\end{prop}
\nin

\mn
{\em Proof.}
This is an application of Proposition~\ref{KVcor1},
in which, having fixed $(i,j)\in\I$, we set
$\theta= \min\{\vt,(r-3)/(r+1)\}$,
$H=H_{ij}$ (so $s=3$), $S=S_{ij}$ and $(x_1,x_2,x_3)=(v,x,y)$.

From \eqref{Sij} we have $S> K\log n$ for any fixed $K$
(and large enough $C$), while,
since $(v_H,e_H) = (j-2, \vs(i,j)+j-3)$,
the combination of \eqref{SijDef}
and \eqref{E0}
(which said $\E_0\sim n^{v_H}p^{e_H}$)
gives
\[
\E_0 < (1+o(1))(np^2)^{-1}S < n^{-\theta}S.
\]
Thus Proposition~\ref{VuCor} will follow from
\beq{Hisbal}
\mbox{$H$ is balanced.}
\enq

\mn
(As will appear below, $H$ is {\em strictly}
balanced unless $(i,j)= (2,4)$.)

\mn
{\em Proof}.  This is a routine verification
and we aim to be brief.
It is enough to show
\beq{yourboat}
\rho(H)\geq \rho(H[k])  ~~\forall 1\leq k \leq j-2,
\enq
where we write $H[k]$ for the subgraph
of $H$ induced by $\{0\dots k\}\cup \{i,j\}$
(so $H[j-1]=H$) and exclude the case $k=i=1$ since it gives
$v(H[k])=0$.
One easily checks that
\[
v(H[k]) = \left\{\begin{array}{cl}

k & \mbox{if $1\leq k<i$,}\\
k-1 &\mbox{if $i\leq k<j$,}
\end{array}
\right.
\]
\[
e(H[k]) = \left\{\begin{array}{cl}
\C{k+1}{2}+2k & \mbox{if $1\leq k<i$,}\\
\C{k+1}{2}+i-2 &\mbox{if $i\leq k<j$,}
\end{array}
\right.
\]
and (consequently)
\[
\rho(H[k]) = \left\{\begin{array}{ll}
\frac{k+5}{2} & \mbox{if $1\leq k<i$,}\\
\frac{k^2+k+2i-4}{2(k-1)}=:f_i(k)
 &\mbox{if $i\leq k<j$.}
\end{array}
\right.
\]
It follows (with a tiny calculation for the third assertion) that
\eqref{yourboat} holds:
strictly if $k\leq i-2$;
with equality if $k=i-1 $ or $k=i=2$; and strictly otherwise
(so if $k\geq i$ and $(k,i)\neq (2,2)$).
This completes the proofs of \eqref{yourboat} and \eqref{Hisbal},
and also shows that we have
strict inequality in the former unless
$k=i=2$ and $j=4$, so, as mentioned
mentioned earlier, strict balance
in the latter unless $(i,j)=(2,4)$.\qed

\mn
{\em Proof of Lemma~\ref{xKLS'}.}
To somewhat lighten the notational load, set
$\ga_r = \frac{2^{r-3}}{(r-3)!}$,
$\gb_r= \frac{(r-3)!}{(r-3)^{r-3}}$ and
$\gc_r=\ga_r\gb_r= (\frac{2}{r-3})^{r-3}$.

\medskip
In what follows, we assume that
$v\in V$ and that $S,T$ are disjoint subsets
of $V\sm \{v\}$
satisfying the size requirements of Lemma~\ref{xKLS'}.
Of course we may also assume $T\neq \0$,
since \eqref{TRST} is vacuous if $T=\0$.
Let
\beq{TvST}
\T(v,S,T) =
\{\tau(S,T,N_v\sm (S\cup T))\geq
4\ga_r\pi |T| \gL_r(n,p)\}\enq
and
\[\R =\bigcup(\{S,T\sub N_v\}\wedge\T(v,S,T)),\]
with the union over $v,S,T$ as above.
Notice that $\ga_r\leq 2$ for all $r$,
so that the expression
$4\ga_r\pi |T| \gL_r(n,p)$ in \eqref{TvST} is at most
the bound in \eqref{TRST};
thus to prove Lemma~\ref{xKLS'} it is enough to show
\beq{L2.2toshow}
\Pr(\R) =o(1).
\enq

Set $\vt = \frac{\gc_r\eps}{80}\C{r-1}{2}^{-2}$ and let
\[
\Q=
\{N(H_{ij},G;v,x,y) < \vt S_{ij}/\log n
~\forall i,j,v,x,y\}
\wedge \{d(v)< 2np ~\forall v\}.
\]
According to Proposition~\ref{VuCor}
and Lemma~\ref{vdegree} we have
\beq{Qbar}
\Pr(\ov{\Q})=o(1).
\enq

Now
\beq{Rsub}
\R\sub \ov{\Q}\cup\bigcup_v\bigcup_W\bigcup_{S,T}
\left[\{N_v=W\}\wedge \T(v,S,T)\wedge\Q\right]
\enq
and
\beq{Rsum}
\Pr(\R)\leq \Pr(\ov{\Q})+
\sum_v\sum_W\Pr(N_v=W)\sum_{S,T}
\Pr(\T(v,S,T)\wedge\Q|N_v=W),
\enq
where (in both \eqref{Rsub} and \eqref{Rsum})
$W$ runs over subsets of $ V\sm \{v\}$ of size at most $2np$
and $(S,T)$ over pairs of disjoint subsets of $W$ with
sizes as in Lemma~\ref{xKLS'}.

Thus, in view of \eqref{Qbar}, we will have \eqref{L2.2toshow}
if we can show that the inner sums in \eqref{Rsum} are small;
for example, it is enough to show that
for fixed $v,W,S,T$ (as above), with $R:=W\sm (S\cup T)$
and
\[
\T(S,T,R):=\{\tau(S,T,R)\geq 4\ga_r\pi |T| \gL_r(n,p)\},
\]
we have
\beq{L22toshow}
\Pr(\T(S,T,R)\wedge\Q|N_v=W)<  \exp[- 4(\pi/\eps) |T| \log n].
\enq
This suffices since for each $t>0$
the number of choices for $S,T$ with $|T|=t$
is less than
$\Cc{n}{t}\Cc{n}{\pi t/\eps} < \exp[2(\pi/\eps)t\log n ],
$
which with \eqref{L22toshow}---recall we assume
$T\neq \0$---bounds the inner sums in
\eqref{Rsum} by $\sum_{t\geq 1}n^{-2t}$.

\mn
{\em Remark.}  The bound on the number of
$(S,T)$'s could be made a little smaller
since $S,T$ are chosen from $W$ rather than from all of $V$, but
there is little to be gained by this
(unlike in the proof of Lemma~\ref{xKLS} where the difference
was crucial); rather, the point here
is that, since $W,S,T$ determine $R$, we avoid
paying an unaffordable $\exp[\gO(|R|\log n)]$
to account for choices of $R$.

\medskip
To slightly streamline some of our expressions, we
now set, for an event $\A$,
$\PPP(\A) =\Pr(\A\wedge \Q|N_v=W)$.
For reasons that will appear below
(see \eqref{Qvsij} and the lines immediately following it),
we will derive \eqref{L22toshow} from the
following multipartized version.
\begin{lemma}\label{Ri's}
For any v,W,S,T as above and partition
$R_1\cup\cdots\cup R_{r-3}$ of $R=W\sm (S\cup T)$,
\beq{tSTR1}
\mbox{$\PPP(\tau(S,T,R_1\dots R_{r-3}) > 2\gc_r\pi |T|\gL_r(n,p))
< \exp[-5 (\pi/\eps)|T|\log n].$}
\enq
\end{lemma}
\mn
{\em Remark.}  Note---{\em cf.} the preceding remark---this
does {\em not} say that w.h.p. (under $\PPP$) we have
$\tau(S,T,R_1\dots R_{r-3}) \leq 2\gc_r\pi |T| \gL_r(n,p)$
for all relevant choices of $S,T$ and $R_i$'s,
since for small $T$
the number of choices for the $R_i$'s
overwhelms the bound in \eqref{tSTR1}.

\medskip
To see that Lemma~\ref{Ri's} implies
\eqref{L22toshow},
choose, independently of $G$,
a random (uniform) ordered partition
$R_1\cup\cdots \cup R_{r-3}$ of $R$.
Given any specification of $G$, say with $\tau(S,T,R)=\tau$,
we have
\[\E\tau(S,T,R_1\dots R_{r-3}) = \gb_r\tau,\]
whence, by Markov's Inequality,
\begin{eqnarray*}
\Pr(\tau(S,T,R_1\dots R_{r-3})<\vs\tau)
&=&
\Pr(\tau-\tau(S,T,R_1\dots R_{r-3})> (1-\vs)\tau)\\
&<& (1-\gb_r)/(1-\vs) = 1 - (\gb_r-\vs)/(1-\vs)
\end{eqnarray*}
for any $\vs\in (0,\gb_r)$,
where ``$\Pr$" now refers only to the choice of partition.

Applying this with $\vs =\gb_r/2$ (and recalling $\gc_r=\ga_r\gb_r$)
gives,
with the natural extension of $\PPP$
to probabilities involving both $G$ and the random
partition,
\[
\mbox{$\PPP(\tau(S,T,R_1\dots R_{r-3})$} >
\mbox{$ 2\gc_r\pi t \gL_r(n,p))$}
~~~~~~~~~~~~~~~~~~~~~~~~~~~~~~~~~~~~~
\]
\[
~~~~~~~~~~~~~~~~~~~~~~~~~~~~~~~~~~~~~
> \mbox{$(\gb_r/2)
\PPP(\tau(S,T,R) > 4\ga_r\pi t \gL_r(n,p))$},
\]
and combining this with \eqref{tSTR1} we have
\[
\mbox{$\PPP(\tau(S,T,R) > 4\ga_r\pi t \gL_r(n,p))
< (2/\gb_r) \exp[-5 (\pi/\eps)t\log n]$}
\]
(which is less than the bound in \eqref{L22toshow}).
\qed

\mn
{\em Proof of Lemma~\ref{Ri's}.}
Fix
$v,S,T,W $ and $R_1\dots R_{r-3}$ as in the lemma
and
set $|S|=s$, $|T|=t$, $|R_i|=r_i$,
$\Psi=\tau(S,T,R_1\dots R_{r-3})$ and
$
\T =\{\Psi > 2\gc_r\pi t \gL_r(n,p)\}.
$

In what follows we will usually be considering
variants of $\Q$ rather than $\Q$ itself, so
abandon the notation $\PPP$ used above; but
we do continue to condition
on $\{N_v=W\}$ and omit this conditioning from the notation.
Thus our target inequality
\eqref{tSTR1} (for our specified $v,S,T,W ,R_1\dots R_{r-3}$)
becomes
\beq{sSTR2}
\Pr(\T\wedge \Q)< \exp[-5 (\pi/\eps)t\log n].
\enq

\medskip
Set $R_{r-2}=S$, $R_{r-1}=T$
and, for $(i,j)\in \I$,
\[
J_{ij} =G\cap \nabla(R_i,R_j).
\]
Note that $\T$ is determined by
$\cup J_{ij}$.

We choose the $J_{ij}$'s in the order given by $\prec$
and set
\[\Psi_\vs=\E[\Psi|(J_{ij}:\vs(i,j)\leq \vs)];\]
in particular $\Psi_0=\E \Psi$
and $\Psi_{\C{r-1}{2}}=\Psi$.
Notice that
\begin{eqnarray}
\Psi_0 &=&
\mbox{$st\prod_{i=1}^{r-3}r_i\cdot p^{\Cc{r-1}{2}}$}
\leq
st(|W|/(r-3))^{r-3}p^{\Cc{r-1}{2}}\nonumber\\
&\leq&
\pi tnp(2np/(r-3))^{r-3}p^{\Cc{r-1}{2}}
= \gc_r \pi t\gL_r(n,p) =:\mu. \label{Psi0}
\end{eqnarray}

Given
\[
G_{ij}:=\nabla(v,W)\cup \bigcup\{J_{kl}:(k,l)\prec (i,j)\}
\]
(note $\nabla(v,W)\sub G$),
we may write
\beq{Psivs}
\Psi_{\vs(i,j)}= \sum\{\xi_{xy}w_{xy}:(x,y)\in R_i\times R_j\},
\enq
where the $\xi_{xy}$'s are independent Ber($p$) r.v.'s and
\beq{wxy}
\mbox{$w_{xy} = M_{xy}\prod_{l=j+1}^{r-1}r_l \cdot
p^{\Cc{r-1}{2}-\vs(i,j)}$,}
\enq
with $M_{xy}$ the number of copies $\vp$ of $H_{ij}$
in
$
G_{ij}$
having
$\vp(u_0)=v$, $\vp(u_i)=x$, $\vp(u_j)=y$
and $\vp(u_l)\in R_l$ for $l\in [j-1]\sm \{i\}$.
(The exponent of $p$ in \eqref{wxy} is the number of
$(k,l)\in \I$ with $(k,l)\succ (i,j)$,
so the number of $J_{kl}$'s that are chosen after $J_{ij}$.)
Of course
\beq{Mxy}
M_{xy}\leq N(H_{ij}; G_{ij},v,x,y).   
\enq

Define events $\Q_\nu$ ($0\leq \nu < \C{r-1}{2}$) by
\beq{Qvsij}
\Q_{\vs(i,j)-1}=
\{N(H_{ij},G_{ij};v,x,y)< \vt S_{ij}/\log n
~~\forall (x,y)\in R_i\times R_j\}
\enq
(with $\Q_0$ the full probability space).
Then
$\Q_\vs\supseteq \Q$ for all $\vs$
and---the point of the
fussy definitions---$J_{ij}$ is independent of
$\{N_v=W\}\wedge \Q_{\vs(i,j)-1}$.

If $\Q_{\vs(i,j)-1}$ holds, then \eqref{wxy} and \eqref{Mxy}
(and the definition
of $S_{ij}$ in \eqref{SijDef}) give
\beq{wxy'}
w_{xy}< 2\vt \glr/\log n
~~~\forall (x,y)\in R_i\times R_j
\enq
(using
$\prod_{l=j+1}^{r-1}r_l < 2(np)^{r-j-1}$).

\medskip
For $\vs \in [\C{r-1}{2}]$,
let
\[
\T_\vs = \{\Psi_\vs-\Psi_{\vs-1}> \Cc{r-1}{2}^{-1}\mu\}.
\]
In view of \eqref{Psi0}, we have $\T \sub \bigcup \T_\vs$
(with the union over $\vs\in [\Cc{r-1}{2}]$), whence,
using $\Q\sub Q_\vs$,
\begin{eqnarray*}
\T\wedge \Q &\sub&
\bigcup_\vs(\T_\vs
\wedge \ov{\T_1}\wedge \cdots \wedge \ov{\T}_{\vs-1})\wedge \Q\\
&\sub &\bigcup_\vs(\T_\vs\wedge \Q_{\vs-1}
\wedge \ov{\T_1}\wedge \cdots \wedge \ov{\T}_{\vs-1});
\end{eqnarray*}
so we
will have \eqref{tSTR1} if we show, for $\vs\in [\Cc{r-1}{2}]$,
\beq{L22toshow'}
\Pr(\T_\vs|\Q_{\vs-1}
\wedge\ov{\T_1}\wedge\cdots \wedge\ov{\T}_{\vs-1})
< \Cc{r-1}{2}^{-1} \exp[- 5(\pi/\eps) t \log n].
\enq
(Of course the first factor on the r.h.s. is unimportant.)

The main effect of
the conditioning in \eqref{L22toshow'} is the inequality
\eqref{wxy'}
implied by $\Q_{\vs-1}$.
A second effect
is that nonoccurence of earlier $\T_\varrho$'s bounds
$\Psi_{\vs-1}$ above by
$2\mu$.

Now
let $\psi=2\mu$,
$\gl =\Cc{r-1}{2}^{-1}\mu$,
$\eta =(2\Cc{r-1}{2})^{-1}$
and
$z= 2\vt \glr/\log n$
(the bound in \eqref{wxy'}).
Then
$\T_\vs$ is the event that $\Psi_\vs$, which
(again, given $G_{ij}$ where $\vs=\vs(i,j)$) is just the
sum in \eqref{Psivs}, is greater than $\Psi_{\vs-1}+\gl$,
where we have $\E\Psi_\vs =\Psi_{\vs-1}\leq \psi$.
We may thus apply Lemma~\ref{Chernwts}
to bound the l.h.s. of
\eqref{L22toshow'} by
\[\exp[-\eta\gl/(4z)]
=\exp[-5(\pi/\eps) t\log n].\]

\section{Proof of Lemma~\ref{L2.2B}}\label{PL2.2B}

As noted earlier, proving Lemma~\ref{L2.2B}
turned out to be quite a bit
trickier than seemed likely on first inspection.
Most interesting here are the roundabout approach
{\em via} Lemma~\ref{Rlemma} (discussed a bit below)
and the use of Lemma~\ref{TRW} in the proof of Lemma~\ref{Rlemma}.
%
(While
it had seemed to us since \cite{DKMantel}
that a proof of Theorem~\ref{MT} for $r=4$
would extend fairly automatically to general $r$,
this turned out to be not quite true, the one point
requiring significant additional ideas being the proof of
Lemma~\ref{L2.2B}.)

One curious point here is that,
while one expects things to get easier as $p$ grows,
our main line of argument runs into difficulties
when $p$ is too far above the lower bound in \eqref{p}.
On the other hand---now more in line with intuition---the
statement for larger $p$ follows quite
easily once we have the following
quantified
version for small $p$.

\begin{lemma}\label{L2.2B'}
For each $\gl>0$ there is a $\vr>0$ such that
for each L, if
\beq{smallp}
p=Cn^{-\tfrac{2}{r+1}}\log^{\tfrac{2}{(r+1)(r-2)}}n,
\enq
with a sufficiently large (fixed) C, then
with probability at least $1-n^{-L}$,
\beq{kappaST'}
\kappa(S,T) < \gl|S|\glr
\enq
whenever
$S\sub \Cc{V}{2} $ and $T\sub G$ satisfy
$\gD_S\leq 2np,$ $V(S)\cap V(T) = \0$
and
$|T|\leq |S|<\vr n^2p$.
\end{lemma}

\mn
{\em Proof of Lemma~\ref{L2.2B} given Lemma~\ref{L2.2B'}.}
This is similar to the derivation of \eqref{L22toshow}
from Lemma~\ref{Ri's}.
Fix $\pi$ as in Lemma~\ref{L2.2B} and
let $\vr$ be the value corresponding to
$\gl =2^{-(\C{r}{2}+2)}\pi$ in Lemma~\ref{L2.2B'}.
We show that Lemma~\ref{L2.2B} holds with $\eps = \vr/2$.

For $p$ as in \eqref{smallp}
and $q>p$,  let $\gz = p/q$,
$G=G_{n,q}$ and $G'=G_\gz$ ($\sim G_{n,p}$),
and write $\kappa(\cdot)$ and $\kappa'(\cdot)$ for
$\kappa_G(\cdot)$ and $\kappa_{G'}(\cdot)$ respectively.
We first observe (this is just for convenience) that we may confine our attention to
$T$'s that are not too small.
According to Corollary \ref{KrCor},
there is a fixed $K$ so that w.h.p.
$G$ satisfies
\beq{kST}
\kappa(S,T)< |S||T|\max\{2n^{r-4}q^{\C{r}{2}-6},K\log n,\}
~~~ \forall S,T\sub \Cc{V}{2}.
\enq
%
But if \eqref{kST} holds then
\eqref{kappaST} is automatic whenever
\beq{Tsmall}
|T| \leq \frac{\pi\glq}{\max\{2n^{r-4}q^{\C{r}{2}-6},K\log n\}}~.
\enq
Note also that the bound here is fairly large
compared to $\gz^{-1}= q/p$, e.g. since each of
$\glq\cdot p/q>\glr$ and $n^2q^5\cdot p/q> n^2p^5$ is at least
a large multiple of $\log n$.

Thus, in view of Lemma~\ref{L2.2B'},
it is enough to show that if $G$ violates the conclusion of
Lemma~\ref{L2.2B} (with $q$ in place of $p$)
at some $S,T$ of size at least the r.h.s. of \eqref{Tsmall}---so in particular
with $\gz |T|$ slightly large---then $G'$
violates the conclusion of Lemma~\ref{L2.2B'}
with probability at least (say)
$n^{-r}$.

To see this, suppose a violation for $G$ occurs at $(S,T)$
with sizes as above.
We then observe that we may choose
$S'\sub S$ with (say)
\beq{S'}
\mbox{$\gD_{S'}\leq 2np$, $~\gz |S|/2< |S'|<2\gz |S|~$ and
$~\kappa(S',T)\geq \gz \kappa(S,T)/2$}
\enq
(so $|S'|< 2\gz \eps n^2q = \vr n^2p$).
Existence of $S'$ is given by Proposition~\ref{PCol},
as follows.
Set $m = \gD_S+1$ and
let $S_1\dots S_m$ be (the color classes of) an
equitable $m$-coloring of $S$, with
$\kappa(S_1,T)\geq \cdots \geq\kappa(S_m,T)$.
Then
$S' = S_1\cup\cdots \cup S_{\lfloor \gz m-1\rfloor}$ satisfies \eqref{S'}.

Now set
$u= \min\{\lceil \gz |T|\rceil, |S'|\}$ and
$v=\C{r}{2}-1$,
and let $T'=G'\cap T$.
We claim that with probability at least (say) $n^{-r}$,
\beq{T'gz}
\mbox{$|T'|=u~~$ and
$~~\kappa'(S',T')\geq (\gz/2)^v\kappa (S',T)/2,
$}
\enq
in which case $S',T'$ (which clearly satisfy the conditions
following \eqref{kappaST'}) violate \eqref{kappaST'},
since
\[
(\gz/2)^v\kappa (S',T)/2
~\geq ~(\gz/2)^v\gz\kappa (S,T)/4
~\geq ~\gl|S'|\glr
\]
(where we used
$
\gk(S,T)~\geq ~\pi |S|\glq
~\geq ~ \pi (2\gz)^{-1}|S'|\glr \gz^{-v}
$).

For the claim,
set
$\kappa' =\kappa'(S',T')$,
$\mu =\gz^v\kappa (S',T)$ ($= \E\kappa' $),
$\mu'=2^{-v}\mu$
and $\Q =\{|T'|=u\}$.
The probability in question is
\[
\Pr(\Q )
\Pr(\kappa'\geq \mu'/2|
\Q ) ~>~ n^{-1}\Pr(\kappa'\geq \mu'/2|
\Q )
\]
(since $\Pr(\Q) =\gO((n^2p)^{-1/2})> 1/n$).
We also have
\[
\E[\kappa'|\Q ]~\geq ~\tfrac{(u)_v}{(|T|)_v}\kappa(S',T)~ > ~(1-\vs)\mu'
\]
for some small constant $\vs$.
(Here we use the assumption that $u$ is fairly large.
The expectation is typically more like $\gz^v\kappa(S',T)$, since most
relevant edges are in $G\sm T$,
whereas the bound allows them to all be drawn from $T$.)
Markov's Inequality thus gives
\begin{eqnarray*}
\Pr(\kappa'< \mu'/2|\Q ) &=&
\Pr(\kappa(S',T)-\kappa'> \kappa(S',T)-\mu'/2|\Q )\\
&<& (\kappa(S',T)-\mu'/2)^{-1}(\kappa(S',T)-\E[\kappa'|\Q ])\\
&<& 1-\mu'/(3\kappa(S',T)) ,
\end{eqnarray*}
so \eqref{T'gz} holds with probability
at least $n^{-1}(\gz/2)^{\C{r}{2}}>
n^{-1}p^{\C{r}{2}}> n^{-r}$.\qed

\medskip
We now turn to the proof of Lemma~\ref{L2.2B'}.
Though the statement here, like that of Lemma~\ref{L2.2B},
is natural,
it seems resistant to frontal assault
(the difficulties are reminiscent of those associated
with upper tail bounds---see e.g. \cite{Chat,DK1} and the
history reviewed in
\cite{IUT} or
\cite[Sec. 2.6]{JLR} for a sort of case study of one such question---though
nothing we know from that arena seems to help here);
so the argument will be
somewhat
indirect.

Recall that $\kappa(S) $ counts choices of $xy\in S$ and
$Z\in  \C{V\sm \{x,y\}}{2}$ with $\C{Z\cup\{x,y\}}{2}\sm \{\{x,y\}\}\sub G$.
The value of $\kappa(S)$ is (essentially) known ({\em cf.} \eqref{kxy}), so we
also know the size of the multiset, say $M$, of edges
appearing in the various $Z$'s (namely $|M|=\kappa(S)\C{r-2}{2}$).
The sense of
Lemmas~\ref{L2.2B} and \ref{L2.2B'} is that if $T$ is only a small part of $G$
then few of these edges should come from $T$.
(Note, though, that if $|S\ll n^2p/\glr$ then {\em all} of the edges in question
lie in some $T$ of size
$\Theta(|S|\glr)\ll n^2p$;
so, as noted in Section~\ref{S2.2}, bounding $|T|$ without reference to $|S|$ will not do here.)

It is thus natural to try to prove Lemma~\ref{L2.2B'} by
controlling
pairs that appear
too often in $M$.
When $r=4 $
there is in fact a (rather long,
martingale-based) argument along these lines that does work;
but we were unable to get that argument, or
any such ``natural" approach to work
for larger $r$,
and instead approach the problem from the other direction,
showing (roughly) that most of $M$ is spread among edges that do
{\em not} have very high multiplicities.
This works best when $S$ (called $R$ below) is small enough that few
of the edges of $G$ are counted even once in $M$.
In this case we are able to show, using Theorem~\ref{TRW},
that the number counted {\em at least} once is close to $|M|$,
leaving little room for higher multiplicities.
(This also requires a stricter bound on $\gD_S$.)
The result for larger $S$ then follows easily, though the
partitioning step that accomplishes this seems rather wasteful.

\bigskip
For $R\sub \C{V}{2}$ and $\{u,v\}\in \C{V}{2}$,
set $\gs_R(u,v)= {\bf 1}_{\{\kappa(R,u,v)>0\}}$
(thus $\gs_R(u,v)=1$ iff there are
$xy\in R$ with $\{x,y\}\cap \{u,v\}=\0$ and
$Z\in \C{V\sm \{x,y,u,v\}}{r-4}$ such that all pairs
from $\{u,v,x,y\}\cup Z$ other than $xy$ are edges of $G$;
in particular $\gs_R(u,v)=0$ if $uv\not\in G$)
and
\[
\gs(R) =\sum\{\gs_R(u,v):\{u,v\}\in \Cc{V}{2}\}.
\]

\mn
The next assertion easily implies
Lemma~\ref{L2.2B'}.

\begin{lemma}\label{Rlemma}
For each $\gb> 0$ there is an $\eta>0$ such that
for each L, if $p$ is as in \eqref{smallp} with a
sufficiently large (fixed) C, then with probability at least
$1-n^{-L}$
\beq{Rreq}
\gs(R) > \tfrac{1- \gb}{2(r-4)!} |R|\gL_r(n,p)
\enq
whenever $R\sub \C{V}{2}$,
\beq{Rreq'}
|R|< \eta n^2p/\gL_r(n,p)=
\eta n^{-(r-4)}p^{-(r^2-r-4)/2},
\enq
and
\beq{gDR}
\gD_R \leq \gb^{-1}np/\gL_r(n,p)=
\gb^{-1} n^{-(r-3)}p^{-(r^2-r-4)/2}.
\enq
\end{lemma}

\mn
{\em Proof of Lemma~\ref{L2.2B'} given Lemma~\ref{Rlemma}.}
Let $\gb=\gl/5$ and
let $\eta$ be the value corresponding to this $\gb$
in Lemma~\ref{Rlemma}.
We show that
$\vr =\gb\eta$ meets the requirements of
Lemma~\ref{L2.2B'}.

By Corollary \ref{KrCor} we know that for any $M$ we have
\beq{kxy}
\kappa(xy) <\tfrac{1+ \gb}{(r-2)!} \gL_r(n,p)  ~~\forall x,y
\enq
with probability at least $1-n^{-M}$, provided $C$
(in \eqref{smallp}) is sufficiently large
(namely, large enough so that the bound in \eqref{kxy} is at least
$K\log n$,
with $K$ as in the corollary).
So in view of Lemma~\ref{Rlemma}, it is enough to show
that
\eqref{kappaST'}
holds (for $S,T$ as in Lemma~\ref{L2.2B'}) 
whenever \eqref{kxy} is true and
\eqref{Rreq} holds for $R$ as in
Lemma~\ref{Rlemma}; we therefore assume these and proceed.

Set $m=4\gb \gL_r(n,p)$ and
let $R_1\cup \cdots \cup R_m$ be a partition of $S$
with, for each $i$ (and with more precision than is really necessary),
\beq{Ribd}
|R_i|<\lceil 2|S|/m \rceil< \eta n^2p/\glr
\enq
and
\beq{dDRibd}
\gD_{R_i}\leq \lceil(\gD_S+1)/m\rceil
<\gb^{-1}np/ \gL_r(n,p).
\enq
(Again, existence of such a partition is given by Proposition~\ref{PCol}, briefly
as follows.
If $m\geq \gD_S+1$, we can take the $R_i$'s themselves to be matchings;
otherwise
we equipartition $S$ into $\gD_S+1$ matchings and take each $R_i$ to be a union of
$\lfloor (\gD_S+1)/m\rfloor$ or $\lceil (\gD_S+1)/m\rceil$ of the matchings
(whence \eqref{dDRibd}),
with $\max|R_i|-\min|R_i|\leq \lceil |S|/(\gD_S+1)\rceil$,
which is easily seen to imply \eqref{Ribd}.)

Then \eqref{Rreq} gives
\[
\gs(R_i) > \tfrac{1- \gb}{2(r-4)!} |R_i|\gL_r(n,p)
~~~\forall i,
\]
while from \eqref{kxy} we have
\[ 
\kappa(R_i) \Cc{r-2}{2} = \Cc{r-2}{2}\sum_{xy\in R_i}\kappa(xy)
<
\tfrac{1+ \gb}{2(r-4)!} |R_i|
\gL_r(n,p)  ~~~\forall i.
\] 
Thus (again, for each $i$)
\begin{eqnarray*}
\kappa(R_i,T)
&\leq &
|T|+ \sum\{(\kappa(R_i,u,v)-1)^+:\{u,v\}\in\Cc{V}{2}\}
\\
&=&
|T| + \kappa(R_i) \Cc{r-2}{2} - \gs(R_i)
~< ~|T| + \gb|R_i|\gL_r(n,p)
\end{eqnarray*}
%
(instead of $\gb$ we could write $\gb/(r-4)!$) and, finally,
\[
\kappa(S,T) = \sum\kappa(R_i,T)
\leq \gb \gL_r(n,p)(4|T| + \sum |R_i|) \leq 5\gb \gL_r(n,p)|S|.
\]\qed

\mn
{\em Proof of Lemma~\ref{Rlemma}.}
It is enough to show that for any $R\sub \C{V}{2}$
satisfying \eqref{Rreq'} and \eqref{gDR},
and any fixed $K$,
\beq{gsRtoshow}
\Pr(\gs(R)< \tfrac{1- \gb}{2(r-4)!} |R|\gL_r(n,p))
< \exp [- K |R|\log n],
\enq
provided $C$ is large enough.

This is the promised
application of Theorem~\ref{TRW}.
Following the notation used there,
we let $i$ run over $\C{V}{2}$ and,
for a given $i=uv$, let $j$ run over pairs $(xy,W)$
with
\beq{xyW}
\mbox{$xy\in R$, $\{x,y\}\cap \{u,v\}=\0$
and $W\in \C{V\sm \{x,y,u,v\}}{r-4}$.}
\enq
For such
$i=uv$ and $j= (xy,W)$, we set
\[
B_{ij}=
\{K(xy;W\cup \{uv\})\sub G\}
\]
(where, as in Lemma~\ref{newHlemma},
$K(xy,Z) := \C{\{x,y,\}\cup Z}{2}\sm \{xy\}$)
and observe that $I_i=\gs_R(u,v)$.
We then need estimates for $\mu$, $\ov{\gD}$ and $\gc$,
the first of which is easy:
\[\mbox{$\mu
 ~~(=\sum_{i,j}\E I_{ij}) ~
= |R|\Cc{n-2}{2}\Cc{n-4}{r-4}p^{\Cc{r}{2}-1}
\sim \tfrac{1}{2(r-4)!} |R|
\gL_r(n,p)$.}\]

The second is not quite so easy, but has essentially
already been worked out in Lemma~\ref{newHlemma},
the only difference being that
each pair $(K,L)$ appearing in
the $\ov{\gD}$ of
\eqref{ovgD} is counted $\C{r-2}{2}^2$ times in the present
\[
\ov{\gD} := \sum\sum\{\E I_{ij}I_{kl}:
(i,j)\sim (k,l)\}.
\]

\mn
(Each $K=K(xy,Z)$ in \eqref{ovgD}
corresponds to $\C{r-2}{2}$ pairs
$(i,j)$ with $i=uv$ for some $\{u,v\}\in \C{Z}{2}$
and $j= (xy,Z\sm \{u,v\})$.)
%
Lemma~\ref{newHlemma} thus implies that for (the present) $\ov{\gD}$
we have
\beq{gDbd}
\ov{\gD} < \gz |R|n^{2r-4}p^{r^2-r-2}/\log n
\enq
for any desired $\gz>0$ provided $C$ is sufficiently large.
(In more detail: given $\gz$, take
$\xi = \C{r-2}{2}^{-2}\gz$, let $\vt$ be the value associated
with this $\xi$ in Lemma~\ref{newHlemma}, and choose $C$
so that
$C^{\C{r}{2}-1}> (\gb\vt)^{-1}$
(see \eqref{dDRibd})
and $C$ is large enough to support
Lemma~\ref{newHlemma} (with this $\xi$ and $\vt$).)

\medskip
For consideration of
$\gc$ ($ =\sum_i\sum_{\{j,k\}}\E I_{ij}I_{ik}$),
we temporarily fix $i=uv$.
A relevant $\{j,k\}$ is then an (unordered) pair of
{\em distinct} pairs of the form $(xy,W)$ as in \eqref{xyW}.
(The two pairs may, for example, use the same $r-2$ vertices, but
must then have different $xy$'s.)

The argument here is similar to that for Lemma~\ref{newHlemma}.
We may
classify a pair $\{j,k\}$ with
$j= (xy,W)$, $k= (x'y',W')$, according to
\[a(j,k) =|\{x,y\}\cap\{x',y'\}|\in \{0,1,2\}\]
and
\[b(j,k) =|(W\cup\{x,y\})\cap(W'\cup\{x',y'\})|
\in \{0,1\dots r-2\},\]
noting that if $a(j,k)=2$ then $b(j,k)\leq r-3$
(and, of course, $b(j,k)\geq a(j,k)$).
It will also be helpful to set
$K_j=K(xy;W\cup \{uv\})$
($=\C{W\cup \{u,v,x,y\}}{2}\sm \{xy\}$).

Write $\g_i$ for the set of relevant $\{j,k\}$'s
and let
\[N(a,b)= |\{\{j,k\}\in \g_i:
a(j,k) =a, b(j,k) =b\}|.\]
Then
\beq{N0b}
N(0,b) < B|R|^2n^{2r-8-b}
\enq
and
\beq{N1b}
\mbox{$N(1,b) < B\sum_xd^2_R(x)n^{2r-7-b}
\leq B|R|\gD_{_R} n^{2r-7-b}$}
\enq
for a suitable $B=B_r$, and
\[N(2,b) < |R|n^{2r-6-b}.\]
(These are the same as the bounds we used for
the $N(a,b)$'s in the proof of Lemma~\ref{newHlemma},
except that what was $r$ is now $r-2$, since we
have set aside the common pair $i=uv$.
Note also that we have slightly different constraints
on the possibilities for $(a,b)$:
in the earlier discussion $(a,b)=(2,2)$ was excluded
because we wanted only {\em overlapping}
pairs $(K,L)$ (that is, pairs sharing an edge not in $R$);
in the present situation $(2,2)$ is
allowed,
but we exclude the possibility $j=k$, a.k.a.
$(a,b)=(2,r-2)$.)

\medskip
On the other hand,
\[
\E I_{ij}I_{ik} = p^{r^2-r-2-|K_j\cap K_k|}
\]
and, for $b(j,k)=b$,
\[
|K_j\cap K_k| ~\left\{\begin{array}{ll}
=\Cc{b+2}{2}-1&\mbox{if $a(j,k)=2$,}\\
\leq \Cc{b+2}{2}&\mbox{otherwise.}
\end{array}\right.
\]
(E.g. if $b=r-2$ the truth is $|K_j\cap K_k|= \C{r}{2}-2$,
but we don't need this.)

\medskip
Let $\mu'=\C{n}{2}^{-1}\mu$
(the average over $i$ of $\sum_j \E I_{ij}$)
and $D=2B(r-4)!$
(with $B$ as in \eqref{N0b}, \eqref{N1b}), chosen so that
\[\mu' D> B |R| n^{r-4} p^{\C{r}{2}-1}.\]

Setting (again with $i$ fixed)
\[S(a,b) =\sum\{\E I_{ij}I_{ik}:a(j,k)=a,b(j,k)=b\}\]
and
combining the above observations
(and little calculations) yields
\begin{eqnarray*}
S(0,b)
&<&  B|R|^2n^{2r-8-b}p^{r^2-r-2-\Cc{b+2}{2}}\\
& < & \mu' \cdot D|R|n^{r-b-4}p^{\Cc{r}{2}-\Cc{b+2}{2}-1}
\\
& < & \mu' \cdot \eta D(np^{(b+3)/2})^{-b}
\end{eqnarray*}
(using \eqref{Rreq'}),
\begin{eqnarray*}
S(1,b)
&<&  B|R|\gD_R n^{2r-7-b}p^{r^2-r-2-\Cc{b+2}{2}}\\
& < & \mu' \cdot D\gD_R n^{r-b-3}p^{\Cc{r}{2}-\Cc{b+2}{2}-1}
\\
& < & \mu' \cdot  D \gb^{-1} (np^{(b+3)/2})^{-b}
\end{eqnarray*}
(using \eqref{gDR}), and
\begin{eqnarray*}
S(2,b)
&<&  |R| n^{2r-6-b}p^{r^2-r-2-\{\Cc{b+2}{2}-1\}}\\
& < & \mu' \cdot Dn^{r-b-2}p^{(r-b-2)(r+b+1)/2}
\\
& = & \mu' \cdot D(np^{(r+b+1)/2})^{r-b-2}.
\end{eqnarray*}
In particular, recalling \eqref{smallp} and the exclusion
of $(a,b)=(2,r-2)$ (and the trivial $b\geq a$),
we have
\[
S(a,b) <\left\{\begin{array}{ll}
D\eta\mu'&\mbox{if $(a,b)= (0,0)$,}\\
o(\mu')&\mbox{otherwise,}
\end{array}\right.
\]
and, now letting $i$ vary and summing over $i$ and $(a,b)$,
\beq{gcbd}
\gc \leq O(\eta \mu)
\enq
(where the implied constant doesn't depend on $\eta$).

\medskip
Finally, taking $t=\gb \mu/2$, and applying
Theorem~\ref{TRW} (using the bounds \eqref{gDbd} and
\eqref{gcbd} and
noting that $X:=\sum I_i= \gs(R)$),
we find that the l.h.s. of \eqref{gsRtoshow} is less than
\begin{eqnarray*}
\Pr(\gs(R)< \mu-t) &<&
\exp[-(t-\gc)^2/(2\ov{\gD})]\\
&<&
\exp[-\tfrac{(\gb/2 - O(\eta))^2}{4\gz((r-4)!)^2}|R|\log n],
\end{eqnarray*}

\mn
which is less than the r.h.s. of \eqref{gsRtoshow} for
suitable $\eta$ and $\gz$.
(Recall---see \eqref{gDbd}---$\gz$ can be made as small as
we like {\em via} a suitable choice of $C$.)\qed

\section{Rigidity and correlation}\label{RandC}

One reason for the difficulty of the
problem treated in this paper is surely the difficulty
of understanding maximum cuts themselves,
an issue whose centrality is perhaps clearer in \cite{BPS};
our ways of dealing with (or avoiding) it
in Sections~\ref{PXdefect} and \ref{PL2.3'}
are based on the notions and soft observations developed here.

\medskip
We again
use $H$ for a general graph on $V$ and
$G$ for $G_{n,p}$.
Let $\cee$ be a collection of {\em balanced} cuts.
The discussion in this section makes sense more
generally,
but all $\cee$'s used in what follows will be
of the type
\beq{cee}
\cee(X) :=\{\Pi=(A_1\dots A_{r-1}):
\mbox{$\Pi$ is balanced, $X\sub A_1$}\}
\enq
for some $X\sub V$.
%
We will often use this with $X=V(Q)$ for some
$Q\sub \C{V}{2}$, in which case we also write
$\cee(Q) $ for $\cee(X).$

For a graph $H$, we use
\[b(\cee, H)= \max\{|\Pi_H|:\Pi\in \cee\}\]
and
\[\max(\cee,H)=\{\Pi\in \cee:
|\Pi_H|=b(\cee, H)\}\]
---we will speak of ``max cuts"---and, for $\Pi\in \cee$, define the {\em defect of}
$\Pi$ {\em relative to $(\cee,H)$}
to be
\[\deff_{\cee,H}(\Pi)= b(\cee, H)- |\Pi_H|.\]

\medskip
Given $\cee$
and $H$, we
define an equivalence relation
``$\equiv$" (or ``$\equiv_{(\cee,H)}$") by:
\[x\equiv y ~~~\mbox{iff}~~~
\Pi(x)=\Pi(y) ~\forall \Pi\in \max(\cee,H),\]
where $\Pi(x)$ is the block of $\Pi$
containing $x$.
Equivalence classes are
$(\cee,H)$-{\em components}, or simply
{\em components} if the identities of $\cee$ and $H$ are
clear.
(Of course if $\cee=\cee(X)$ then $X$ is automatically contained in
some component, whatever the value of $H$.)

Given $\cee$,
say
$H$ is {\em rigid}
if
the number of equivalent pairs under $\equiv_{(\cee,H)}$
is at least $(1-\ga)n^2/(2(r-1))$.
(Recall $\ga$ was one of the basic constants previewed at the end of Section~\ref{Usage}.)
\begin{prop}\label{Pcore}
If H is rigid then there are distinct $(\cee,H)$-components
$S_1\dots S_{r-1}$ of size greater than $n/r$.
\end{prop}

For a rigid $H$ we will call the (necessarily unique)
collection $\{S_1\dots S_{r-1}\}$ in Proposition~\ref{Pcore}
the {\em core} of $H$.
(Note that, in contrast to our usage for cuts, we
think of the core as {\em un}ordered.)
Of course a nonrigid $H$ may also admit
$S_1\dots S_{r-1}$ as in the proposition; but it
will be convenient in what follows to regard only rigid graphs
as having cores, so if we speak of the core of $H$,
then $H$ is rigid by definition.

\mn
{\em Proof of Proposition} \ref{Pcore}.
This is given by the following assertion,
applied when $H$ is rigid with $(\cee,H)$-components
$S_1\dots S_m$.

\mn
{\em Claim.}
If $S_1\cup\cdots\cup S_m$ is a partition of $V$
with $s_i:=|S_i|\leq (1+\gd)n/(r-1)$ $\forall i$ and
$\sum\C{s_i}{2} > (1-\ga)n^2/(2(r-1))$, then
some $r-1$ of the $s_i$'s are greater than
$(1-r\ga)n/(r-1)$.

\mn
(So we actually get
the proposition with $(1-r\ga)n/(r-1)$ in place of $n/r$;
but $n/r$ is convenient and
sufficient for our purposes.)

\medskip
For the proof of the claim, set $\gl=r\ga$.
Among $(S_1\dots S_m)$'s for which
the conclusion fails,
$\sum\C{s_i}{2}$ is maximum when $m=r$,
$s_1=\cdots =s_{r-2}=(1+\gd)n/(r-1)$ and
$s_{r-1}=(1-\gl)n/(r-1)$ (so $s_m= n-(s_1+\cdots +s_{r-1})$).
This gives
\begin{eqnarray*}
(1-\ga)\tfrac{n^2}{2(r-1)} &<&\sum\Cc{s_i}{2}\\
&<& \left[(r-2)(1+\gd)^2 + (1-\gl)^2 + (\gl-(r-2)\gd)^2\right]
\tfrac{n^2}{2(r-1)^2}\\
&<& (1-\ga)
\tfrac{n^2}{2(r-1)},
\end{eqnarray*}
a contradiction
(where we used $\ga \gg \gd$
for the final inequality).\qed

For (rigid) $H$ with core $\{S_1\dots S_{r-1}\}$,
we say $Q\sub \C{V}{2}$ is {\em in the core}
if $V(Q)$ is contained in
one of $S_1,S_2,\ldots, S_{r-1}$.
(We will only use this with $\cee=\cee(Q)$, but note---a
point that will cause some trouble below---this
does not guarantee that
$Q$ is in the core.)

For any $H$, set
\beq{crit}
\crit(H)~~ (=\crit_\cee(H))~ =
H\cap \bigcap \{\ext(\Pi):\Pi\in \max(\cee,H)\};
\enq
thus $e\in H$ is in $\crit(H)$ iff $b(\cee, H-e)< b(\cee,H)$.
Notice in particular that if $\{S_1,S_2,\ldots, S_{r-1}\}$
is the core of $H$, then
$ 
\nabla(S_1,S_2,\ldots, S_{r-1})\sub \crit(H).
$ 

\medskip
The next two lemmas are the promised applications
of Harris' Inequality (see Section \ref{SecHarris}).
As mentioned earlier, these were  suggested by the
way Harris is used in \cite{BPS};
the crucial new idea here appears in \eqref{uniquecore},
where uniqueness of the core bounds the sum of probabilities by 1
(and again in the proof of Lemma~\ref{transfer2}, where
uniqueness is arranged in a simpler way).

We again write $G$ for $G_{n,p}$.

\begin{lemma}\label{transfer}
Fix $X\sub V$.
Suppose that for each collection
$\{T_1\dots T_{r-2}\}$ of disjoint subsets of $V\sm X$,
the event $F(T_1\dots T_{r-2})$ ($=F(\{T_1\dots T_{r-2}\})$) is
decreasing in and determined by
$\nabla(X,T_1\dots T_{r-2})$,
and that
\beq{FST}
\Pr(F(T_1\dots T_{r-2}))< \xi ~~\mbox{whenever} ~
|T_1|\dots |T_{r-2}| > n/r.
\enq
Given $\cee$, let $\R$ be the event
that $G$ is rigid, say with core $\{S_1\dots S_{r-1}\}$,
$X\sub S_1$, and $F(S_2\dots S_{r-1})$ holds.
Then $\Pr(\R)< \xi$.
\end{lemma}
\nin
{\em Proof.}
For disjoint $S_1\dots S_{r-1}\sub V$ with
$X\sub S_1$, set
\[
E(S_1\dots S_{r-1})=
\{\mbox{$G$ has core $\{S_1\dots S_{r-1}\}$}\}.
\]
The main point (justified below) is that,
for any such $S_1\dots S_{r-1}$,
\beq{EST}
\mbox{$E(S_1\dots S_{r-1})$
is increasing in
$\nabla(X,S_2\dots S_{r-1})$},
\enq
whence,
by Theorem~\ref{Harris} (applied to the indicators of $E$
and $F$),
\[
\Pr(E(S_1\dots S_{r-1})\wedge F(S_2\dots S_{r-1}))
~~~~~~~~~~~~~~~~~~~~~~~
\]
\beq{PrEs}
~~~~~~~~~~~~~~~~~
\leq \Pr(E(S_1\dots S_{r-1}))
\Pr(F(S_2\dots S_{r-1})).
\enq
This gives the lemma, since
\begin{eqnarray}
\Pr(\R)&=&
\mbox{$\sum\Pr(E(S_1\dots S_{r-1})
\wedge F(S_2\dots S_{r-1}))$}\nonumber\\
&<& \mbox{$\xi\sum\Pr(E(S_1\dots S_{r-1}))
~\leq ~\xi$,}\label{uniquecore}
\end{eqnarray}
where the sums are over $(S_1\dots S_{r-1})$ as above
(that is, the $S_i$'s are disjoint with $X\sub S_1$)
and the first inequality uses \eqref{PrEs} and \eqref{FST}
(the latter applicable because $E(S_1\dots S_{r-1})$ implies $|S_i|>n/r$ $\forall i$).

The reason for \eqref{EST} is simply that if
$E(S_1\dots S_{r-1})$ holds, then
adding a pair from
$\nabla(X,S_2\dots S_{r-1})$
(or, for that matter, $\nabla(S_1\dots S_{r-1})$)
to $G$ doesn't affect the set of max cuts:
any such pair is in $\ext(\Pi)$ for every $\Pi\in \max(\cee,G)$,
so each such $\Pi$
remains a max cut, and, moreover (since $b$ increases),
no new cuts are added to $\max(\cee,G)$.\qed

\begin{lemma}\label{transfer2}
Fix $X\sub V$ and an order ``$\prec$" on $\cee=\cee(X)$.
Suppose that for each $(r-2)$-tuple $(B_1\dots B_{r-2})$ of
disjoint subsets of $V\sm X$,
$F(B_1\dots B_{r-2})$
is an event decreasing in and determined by
$\nabla(X,B_1\dots B_{r-2})$,
and that
\beq{FST'}
\Pr(F(B_1\dots B_{r-2}))< \xi ~~\mbox{whenever}
~|B_1|\dots |B_{r-2}| > (1-\gd)n/(r-1).
\enq
Let
$\R$ be the event that
for the first member (under $\prec$),
say $(A_1\dots A_{r-1})$, of $\max(\cee,G)$,
$F(A_2\dots A_{r-1})$ holds.
Then $\Pr(\R)< \xi$.
\end{lemma}

\nin
{\em Proof.}
This is similar to the proof of
Lemma~\ref{transfer}.
For $A_1\dots A_{r-1}$ partitioning $ V$,
the event
\[
E(A_1\dots A_{r-1})=
\{\mbox{$(A_1\dots A_{r-1})$
is the first member of $\max(\cee,G)$}\}
\]
(which in particular implies $X\sub A_1$)
is increasing in $\nabla_G(X,A_2\dots A_{r-1})$.
(If $E(A_1\dots A_{r-1})$ holds, then
adding a pair from
$\nabla_G(X,A_2\dots A_{r-1})$
doesn't remove $(A_1\dots A_{r-1})$ from,
or add any new members to, the
set of max cuts
(though here the set of max cuts may shrink).

The rest of the earlier argument applies without modification.
(Note that, since members of $\cee$ are required to be balanced,
$E(A_1\dots A_{r-1})=\0$ unless $|A_i|>(1-\gd)n/(r-1)$
$\forall i$.)\qed

\section{Proof of Lemma~\ref{Xdefect}}\label{PXdefect}

{\em Note.}  Here and in Section~\ref{PL2.3'}
it will sometimes be better to speak of a {\em set of graphs} rather than
an event; in particular this will be helpful when the discussion involves
more than one random graph.
The default remains our usual $G=G_{n,p}$; that is, when we say without
qualification that some event holds, we mean it holds for $G$.

\medskip
We first observe that it is enough to prove Lemma~\ref{Xdefect} for
\beq{t}
t<Kp^{-1}
\enq
for a suitable fixed $K$:

\begin{prop}\label{propt}
There is a K such that w.h.p. no balanced cut admits more
than $Kp^{-1}$ bad vertices.
\end{prop}
\mn
{\em Proof.}
By Proposition~\ref{vdegree} it is enough to bound the
probability that some balanced $\Pi$
admits
$Kp^{-1}$ bad vertices $x$ with
\beq{dxagain}
d(x) = (1\pm o(1))np.
\enq
Here we use $\cc_r< \bbb_r$.
(Recall these were defined in \eqref{abphi}.)
If
$x$ satisfying \eqref{dxagain} is bad for $\Pi=(A_1\dots A_{r-1})$
(so $x\in A_1$),
then we assert that, writing $d_i $ for $d_{A_i}(x)$,
at least one of $d_2\dots d_{r-1}$
is less than $(1-\vs)np/(r-1)$,
for some (fixed) $\vs=\vs_r>0$.  To see this without too much calculation,
consider the ``ideal" version in which $d(x)=np$
and each of $d_2\dots d_{r-1}$ is at least
$np/(r-1)$.  In this case $D_\Pi(x)$
(defined following \eqref{D})
is minimum
when $d_2=\cdots =d_{r-2}=np/(r-1)$ and $d_{r-1}= 2np/(r-1)$,
in which case
({\em cf.} the remark preceding Lemma~\ref{Xdefect})
\[
D_\Pi(x) =
\Cc{r-3}{2}(\tfrac{np}{r-1})^2 + (r-3)\cdot 2(\tfrac{np}{r-1})^2
= \tfrac{r(r-3)}{2(r-1)^2}n^2p^2 = \bbb_rn^2p^2.
\]
It is then clear that for a small enough
$\vs$ ($=\vs_r>0$),
replacing these ideal assumptions by
\eqref{dxagain} and $d_i> (1-\vs)np/(r-1)$
($i\in [2,r-1]$)
still forces
$D_\Pi(x) > \cc_r n^2p^2$, contradicting the assumption
that $x$ is bad for $\Pi$.

So if $\Pi$ admits at least $t$ bad vertices,
then there are an $i\in [2,r-1]$ and some
$T\sub A_1$ of size at least $t/(r-2)$, each of whose vertices
$x$ has $d_{A_i}(x)<(1-\vs)np/(r-1)$.  But if $\Pi$ is
balanced (so $|A_i|>(1-\gd)n/(r-1)$),
then this implies (say)
\[
|\nabla(T,A_i)| < (1-\vs+2\gd) |T||A_i|p,
\]
which, for the $K$ corresponding to $\eps = \vs-2\gd$
and $c = (1-\gd)/(r-1)$ in Proposition~\ref{nablaST},
violates the conclusion of that proposition
if $t > Kp^{-1}$.
\qed

\medskip
{\em We assume for the rest of
this section that $t$ ranges over values
satisfying \eqref{t}
and
$X$ over $t$-subsets of $V$}.
Set
\[\vt = (\bbb_r-\cc_r)/3.\]
We will prove Lemma~\ref{Xdefect} with $\nu =\vt \eps$, with
$\eps$ as in the paragraph following the proof of Lemma~\ref{lemmaG}.

Let $\Q$ be the event that the conclusions of
Propositions \ref{vdegree} and \ref{propG[X]}(b)
hold and
\[\R_X =\{\exists
\Pi\in \cee(X):
\deff_G(\Pi)<\vt tn^{3/2}p^2 ~\mbox{and each $x\in X$
is bad for $\Pi$}\};\]
so we should show
$\Pr(\cup \R_X) =o(1) $.
We have
\beq{PrU}
\Pr(\cup \R_X) ~\leq ~\Pr(\ov{\Q}) + \sum\Pr(\R_X\wedge\Q)
~<~ o(1)+ \sum\Pr(\R_X\wedge\Q),
\enq
so
will be done if we show, for each $t$ and $X$,
\beq{QXD}
\Pr(\R_X\wedge\Q) < \exp[-\gO(tnp)].
\enq

\medskip
From this point we fix ($t$
and) $X$
and set $W=V\sm X$, $\R=\R_X$.
The strategy will involve choosing $G$ in two stages.
We hope to arrange that the output, $G'$, of the first stage
admits some ``reference" cut, say $\Pi^*$, that is both maximum in $G'$ and
poised to gain many edges in the second stage,
whereas it is likely that each ``bad" cut (i.e. one for which $X$ is bad)
sees significantly fewer additions.
If this does happen then $\Pi^*$
will (typically) be significantly
larger than any bad cut once we add the contributions of the second stage.
(At a nontechnical level this echoes
what was perhaps the main idea
of \cite{DKMantel}; see the proof of (32) there
and item {\bf D} in Section~\ref{Remarks} below.)

The second stage will be confined to
the (random) set of pairs $uv$ having at least---usually exactly---one
common neighbor in $X$, so that
the likely number of additions to a particular $\Pi$ grows
with the number of such pairs in $\ext(\Pi)$,
roughly
the sum over $x\in X$ of the $D_\Pi(x)$'s.
Thus a $\Pi$ for which $X$ is bad will tend to suffer in the second stage;
the more interesting question
is, how do we know that there is some max cut
for which the $D_\Pi(x)$'s are large?
The answer is provided by
Lemma~\ref{transfer2}, but circuitously.

The lemma easily gives the desired cut for our usual
$G=G_{n,p}$, or, more generally, for a $G$ with edges chosen
independently with large enough probabilities.
But $G'$ will not be of this type and the lemma
seems not to apply directly;
instead we apply it to an (auxiliary) copy, $H$, of $G_{n,p}$,
and then
couple $H$ with $G'$.
This looks unpromising
since the distributions are not at all close,
but succeeds roughly because
even the tiny probability that the two graphs coincide
is much larger than the probability that the event of interest fails for $H$.

\medskip
Fix some order ``$\prec$" on $\cee:=\cee(X)$
and let $\T$ be the set of graphs $H$ (on $V$) for
which the first member, say
$(A_1\dots A_{r-1})$, of $\max(\cee,H)$
satisfies
\beq{xinX}
|\{x\in X: \min\{|N_H(x)\cap A_i|:i\in [2,r-1]\}
< (1-\vt)np/(r-1)\}| > \vt t.
\enq

\nin
{\em Remark.}
As suggested earlier, the argument
will now
involve some interplay of different
random graphs, and we need to be clear as
to which graph is meant when we speak of
membership in $\T$.
In particular the next lemma refers
to a {\em generic} copy of $G_{n,p}$; it
will be used to prove that a similar statement holds for
a slightly mongrelized version of the graph we're really
interested in.

\begin{lemma}\label{lemmaG}
$\Pr(G_{n,p}\in \T)<\exp[-\gO(tnp)]$
\end{lemma}
\nin
(where the implied constant depends on $\vt$).

\mn
{\em Proof.}
For
disjoint $B_1\dots B_{r-2}\sub W$, let
$F(B_1\dots B_{r-2}) $
be the event
\[\{|\{x\in X: \min\{d_{B_i}(x): i\in [r-2]\}
< (1-\vt)np/(r-1)\}| > \vt t\}.\]
According to Lemma~\ref{transfer2} it is enough to show that
\eqref{FST'} holds for some $\xi = \exp[-\gO(tnp)]$.

To see this, fix
$B_1\dots B_{r-2}$ as above with $|B_i|> (1-\gd)n/(r-1)$ $\forall i$.
If $F(B_1\dots B_{r-2}) $ holds then there are $i\in [r-2]$
and $Y\sub X$ with $|Y|> \vt t/(r-2)$ such that
$d_{B_i}(x)< (1-\vt)np/(r-1)$ $\forall x\in Y$.
By Theorem~\ref{Chern},
the probability that this occurs for a given $i$ and $Y$ is
less than $\exp [-\gO(tnp)]$
(where, again, the implied constant---roughly $\vt^2/(2r)$---depends on $\vt$);
so, accounting for the
number of possibilities for $i$ and $Y$,
we have
\[
\Pr(F(B_1\dots B_{r-2})) < r2^t\exp[-\gO(tnp)]=\exp[-\gO(tnp)].
\]
\qqqed

\medskip
We now generate $G$
(the version of $G_{n,p}$ in which we're
really interested)
in stages.
For $L\sub \nabla(X,W)$, set
$P_L =\cup_{x\in X}\C{N_L(x)}{2}$
and $Q_L = \C{X}{2}\cup (\C{W}{2}\sm P_L)$.
Fix $\eps>0$ with $\eps^2$ small compared
to the implied constant in Lemma~\ref{lemmaG},
and set $q=\eps n^{-1/2}$.
We choose edges of
$G$ in
the following order.

\mn
(i)  Choose $L=\nabla_G(X,W)$
and set $P_L=P$ and $Q_L=Q$.

\mn
(ii)  Choose $G\cap Q$.

\mn
(iii)  Choose edges in $P$
{\em with probability}
$p'$, where $1-p'=(1-p)/(1-q)$ (so $p'\approx p-q$),
these choices made independently.

\mn
(iv)  Choose additional edges in $P$ (again independently)
with probability $q$.

\mn
(Note that the resulting $G$ is indeed a copy of $G_{n,p}$.)
Let $G'$ be the output of (i)-(iii),
and
\beq{sss}
\sss =\{|G'\cap P|\leq 2|P|p\}.
\enq

Let $\Qq\supseteq \Q$
be the event that $G$
satisfies the
conditions:
\beq{Q*1}
d(x)=(1\pm o(1))np ~~\forall x\in X;
\enq
\beq{Q*2}
d(x,y)=(1\pm o(1))np^2 ~~\forall x,y\in X ~(x\neq y);
\enq
\beq{Q*3}
|G[X]| < t\log n.
\enq
%
%

\mn
Note that membership of $G$ in $\Q^*$ depends only on the edges chosen
in (i) and (ii), so $G\in \Q^*$ is the same as $G'\in \Q^*$;
this allows us to continue to use notation like $d(x)$, $d_A(x,y)$, $N(x)$ for $x,y\in X$
without ambiguity.

\medskip
We need two easy consequences of $\Qq$
(actually of \eqref{Q*1} and \eqref{Q*2}): first,

\beq{|P|}
|P| = (1\pm o(1))tn^2p^2/2,
\enq
and second,
for any disjoint $S,T\sub W$,
\beq{YST}
|P\cap \nabla (S,T)| >
\sum_{x\in X}d_S(x)d_T(x)- o(tn^2p^2).
\enq

\mn
{\em Proof of} \eqref{|P|}.
We have
\beq{Pub}
|P|\leq \sum_{x\in X} \Cc{d_W(x)}{2} \leq \sum_{x\in X}\Cc{d(x)}{2}
<(1+o(1))t n^2p^2/2
\enq
(with the last inequality given by \eqref{Q*1}).
For a lower bound we may use
\begin{eqnarray*}
|P| &\geq&
|\cup_{x\in X}\Cc{N(x)}{2}|- (|\Cc{X}{2}|+|\nabla(X,W)|)
\\
&\geq &\sum_{x\in X} \Cc{d(x)}{2}
-\sum_{\{x,y\}\in \Cc{X}{2}}
\Cc{d(x,y)}{2}- tn.
\end{eqnarray*}
By \eqref{Q*1} and \eqref{Q*2} the first sum is $(1\pm o(1))t n^2p^2/2$ and the
second is at most $(1+o(1))t^2n^2p^4/4$.
Combining these observations (and recalling \eqref{t} and
\eqref{po1})
gives \eqref{|P|}.\qqqed

\mn
{\em Proof of} \eqref{YST}.
This is similar.  We have
\begin{eqnarray*}
|P\cap \nabla (S,T)| &=&
|\bigcup_{x\in X}\nabla (N_S(x),N_T(x))| \\
&>& \sum_{x\in X}d_S(x)d_T(x)-
\sum\{d_S(x,y)d_T(x,y):\{x,y\}\in \Cc{X}{2}\},
\end{eqnarray*}
and (again using \eqref{Q*2}, \eqref{t} and
\eqref{po1}) the subtracted term is less than
$|X|^2n^2p^4 = o(tn^2p^2)$.
\qqqed

\medskip
Returning to \eqref{QXD}, we have
\begin{eqnarray*}
\Pr(\R\wedge\Q)&\leq  &\Pr(\R\wedge\Qq) \\
&\leq &\Pr(\ov{\sss} \wedge\Qq) + \Pr(\Qq\wedge (G'\in \T ))
+ \Pr(\R|\Qq\wedge \sss\wedge (G'\not\in \T)),
\end{eqnarray*}
(recall $\sss$ was defined in \eqref{sss})
and, from \eqref{|P|} and Theorem~\ref{Chern},
\[
\Pr(\ov{\sss} \wedge\Qq)\leq \Pr(\ov{\sss}|\Qq) < \exp[-\gO(tn^2p^3)].
\]
Thus \eqref{QXD} (and
Lemma~\ref{Xdefect}) will follow from the next two assertions, the
more interesting of which is the first.

\mn
{\em Claim} 1.
$\Pr(\Qq\wedge (G'\in \T )) < \exp[-\gO(\eps^2 tnp)]$.

\mn
{\em Claim} 2.
$\Pr(\R|\Qq\wedge \sss\wedge (G'\not\in \T)) < \exp[-\gO(tn^{3/2}p^2)]$.

\mn
(Note $n^{3/2}p^2\gg np$.
The implied constant in Claim 1 doesn't depend on $\eps$; it could of course
absorb the $\eps^2$, but we prefer the current form as it
better reflects the source of the bound.)

\mn
{\em Proof of Claim} 1.
This is achieved by a comparison (coupling) of $G'$ and
$G_{n,p}$.
Let $H$ consist of the edges chosen in (i) and (ii) together
with edges in $P$ chosen independently (of $G'$ and each other),
each with probability $p$.  Then $H\sim G_{n,p}$, so by
Lemma~\ref{lemmaG} we have
\beq{PrHT}
\Pr(H\in\T)<\exp[-\gO(tnp)].
\enq

\medskip
Let $\g = \{K\in \Qq:|K\cap P|< |P|(p-2q)\}$.
(We again note that membership of $K$ in $\Q^*$ depends
only on the edges of $K$ incident with $X$.)
By Theorem~\ref{Chern}
(recalling that $\Qq$ implies
\eqref{|P|})
we have
\beq{G'ing}
\Pr(G'\in \g) < \exp[-\gO(\eps^2  tnp)].
\enq

On the other hand, we assert,
\beq{Kgood}
K\in \Qq\sm \g ~\Ra ~ \Pr(G'=K)< \exp[O(\eps^2ntp)]\Pr(H=K).
\enq
{\em Proof.}
Fix $K\in \Qq\sm\g$, say with $\nabla_K(X,W)=L$, and
let $P=P_L$,
$m=|P|$ ($\sim tn^2p^2/2$ since $K\in \Qq$)
and $k = |K\cap P|$ ($> m(p-2q)$).
Then
\[
\Pr(G'=K) = \Pr(G'\sm P=K\sm P)\Pr(G'=K|G'\sm P=K\sm P)
\]
and
\[\Pr(G'=K|G'\sm P=K\sm P) =\Pr(B(m,p')=k)\Cc{m}{k}^{-1}.\]
Repeating this with $H$ in place of $G'$ and using
$\Pr(G'\sm P=K\sm P)=\Pr(H\sm P=K\sm P))$ gives
\[
\frac{\Pr(G'=K)}{\Pr(H=K)} = \frac{\Pr(B(m,p')=k)}{\Pr(B(m,p)=k)}.
\]
The r.h.s. is less than 1 if $k> mp$, and otherwise is
less than
\beq{binom.l.b}
[\Pr(B(m,p)=k)]^{-1} = \exp[O(\eps^2ntp)]
\enq
(routine calculation omitted), so we have \eqref{Kgood}.
\qed

\medskip
Finally, \eqref{G'ing}, \eqref{Kgood} and \eqref{PrHT} give
\begin{eqnarray*}
\Pr(\Qq\wedge (G'\in \T)) &\leq &
\Pr(G'\in \g) + \sum\{\Pr(G'=K):K\in \T\cap (\Qq\sm \g)\}\\
&< &\Pr(G'\in \g)
+\exp[O(\eps^2ntp)] \Pr(H\in \T\cap (\Qq\sm \g))\\
&<& \exp[-\gO(\eps^2  ntp)]
\end{eqnarray*}
(where the last inequality uses our assumption on $\eps$).
\qqqed

\mn
{\em Proof of Claim} 2.
Fix $G'\in \Qq\wedge \ov{\T}$ satisfying $\sss$
and let $\Pi =(A_1\dots A_{r-1}) $ be the first member of
$\max(\cee,G')$; so we are assuming \eqref{xinX} fails
(with $G'$ in place of $H$).
Let $G''=G\sm G'$.
We have
\[|\Pi_G| = |\Pi_{G'}|+ |G''\cap \nabla(A_1\dots A_{r-1})|\]
and, for any $\Pi'=(S_1\dots S_{r-1})\in \cee$,
\begin{eqnarray*}
|\Pi'_G| &=& |\Pi'_{G'}|+ |G''\cap \nabla(S_1\dots S_{r-1})|\\
&\leq & |\Pi_{G'}|+ |G''\cap \nabla(S_1\dots S_{r-1})|,
\end{eqnarray*}
whence
\begin{eqnarray}
\deff_G(\Pi')& \geq &|\Pi_G|-|\Pi'_G| \nonumber\\
&\geq & |G''\cap \nabla(A_1\dots A_{r-1})|
-|G''\cap \nabla(S_1\dots S_{r-1})|. ~~
\label{G''cap}
\end{eqnarray}
So---as we will explain in a moment---it is enough to show

\begin{lemma}\label{LCl2}
With probability $1-\exp[-\gO(tn^{3/2}p^2)]$
(where the implied constant depends on $\vt$ and
$\eps$),
\beq{G''1}
|G''\cap \nabla(A_1\dots A_{r-1})| > (1-3\vt)\bbb_rtn^2p^2q
\enq
and
\begin{eqnarray}
|G''\cap \nabla(S_1\dots S_{r-1})|
&< &|P\cap \nabla(S_1\dots S_{r-1})|q+\vt t n^2p^2q
~~~~~~~~~~~\nonumber\\
&& ~~~~~~~~~~~~~~~~~~~~~~
\forall (S_1\dots S_{r-1})\in \cee.\label{G''2}
\end{eqnarray}
\end{lemma}
\nin

\medskip
To see that Lemma~\ref{LCl2}
implies Claim 2, notice that for any $\Pi'=(S_1\dots S_{r-1})\in \cee$
for which all vertices of $X$ are bad, we have
\begin{eqnarray*}
|P\cap \nabla(S_1\dots S_{r-1})|
&\leq& \sum_{x\in X}|\Cc{N(x)}{2}\cap \nabla(S_1\dots S_{r-1})|\\
& =& \sum_{x\in X}D_{\Pi'}(x) ~<~ t\cc_r n^2p^2.
\end{eqnarray*}
\nin
So if
\eqref{G''1} and \eqref{G''2} hold then
in view of \eqref{G''cap} we have, for every such $\Pi'$,
\beq{deffG}
\deff_G(\Pi')> [(1-3\vt)\bbb_r - \cc_r-\vt]tn^2p^2q
> \vt tn^2p^2q
= \nu tn^{3/2}p^2.
\enq
Thus (more or less repeating), failure of $\R$
implies that either \eqref{G''1} or \eqref{G''2} is violated,
which by Lemma~\ref{LCl2} occurs with probability at most
$\exp[-\gO(tn^{3/2}p^2)]$, as required for Claim 2.\qed

\mn
{\em Proof of Lemma~\ref{LCl2}.}
Notice that for any $\Pi'=(S_1\dots S_{r-1})\in \cee$,
\beq{G''}
G''\cap \nabla(S_1\dots S_{r-1})
= (P\cap \nabla(S_1\dots S_{r-1}))\sm G')_q.
\enq
(Recall the r.h.s. was defined in \eqref{Xp}.
It may be helpful to observe
that we could replace
$S_1$ by $S_1\sm X$ on the r.h.s. of \eqref{G''},
since $P$ does not contain pairs
meeting $X$.)

We first consider \eqref{G''1}.
Set $D(x) = D(x;A_1\sm X,A_2\dots A_{r-1})$
(recalling that this notation was introduced in \eqref{D}).
From \eqref{YST}
we have
$|P\cap \nabla(A_1\dots A_{r-1})| > \sum_{x\in X} D(x)-o(tn^2p^2)$,
so also (since $|G'\cap P|=O(tn^2p^3)=o(tn^2p^2)$,
as follows from \eqref{|P|}, $\ov{\sss}$ and \eqref{po1}),
\beq{Pcapnabla}
|(P\cap \nabla(A_1\dots A_{r-1}))\sm G'|
> \sum_{x\in X} D(x)-o(tn^2p^2).
\enq

Let $m=(1-\vt)np/(r-1)$ and
\[Y=
\{x\in X: \min\{|N(x)\cap A_i|:i\in [2,r-1]\}
> m\}
\]
(the complement in $X$ of the set in \eqref{xinX}
when $H=G'$).
We assert that for $x\in Y$ we have
\beq{D(x)'}
D(x) > (1-\vt)\bbb_rn^2p^2.
\enq
To see this, notice that
\[  
D(x) \geq
\Cc{d_W(x)}{2} - (r-3)\Cc{m}{2}-\Cc{d_W(x)-(r-3)m}{2},\\
\]  
since we minimize $D(x)$
(subject to $x\in Y$) by taking $r-3$ of the sets
$N(x)\cap A_i$ ($i\in [2,r-1]$) to be of size $m$
and one to be of size $d_W(x)-(r-3)m$ (and $N(x)\cap (A_1\sm X)$ to be empty).
A little straightforward calculation,
using
\[d_W(x) > (1-o(1))np -|X| =(1-o(1))np \]
(which follows from $G'\in \Q^*$ and \eqref{t}),
then gives \eqref{D(x)'}
(with the ``$\vt$" actually about $2\vt/r$.)\qqqed

\medskip
Thus,
since $G'\not\in \T$ (that is,
$|Y|> (1-\vt)t$), we have
$
\sum_{x\in X} D(x)>(1-\vt)^2\bbb_rtn^2p^2,
$
which with \eqref{Pcapnabla} gives (say)
\[
|(P\cap \nabla(A_1\dots A_{r-1}))\sm G'|
> (1-2\vt)\bbb_rtn^2p^2.
\]
Finally, \eqref{G''} and Theorem~\ref{Chern}
give
\begin{eqnarray*}
\Pr(|G''\cap \nabla(A_1\dots A_{r-1})| < (1-3\vt)\bbb_rtn^2p^2q)
&<& \exp[- \vt^2 \bbb_r tn^2p^2q/2]\\
&=& \exp[-\tfrac{\eps\vt^2\bbb_r}{2}  tn^{3/2}p^2].
\end{eqnarray*}
(Here we are actually
using an easy consequence/extension
of Theorem~\ref{Chern}:
for $\xi =B(n,p)$, $s\leq np$ and any $\gl\geq 0$,
\[
\Pr(\xi < s-\gl) < \exp [- \gl^2/(2s)].
\]
To get this from Theorem~\ref{Chern}, we may, for
example, take
$\xi' =B(n,q)$ with $q=s/n\leq p$ and use
$\Pr(\xi < s-\gl)
\leq  \Pr(\xi' < s-\gl) <\exp [- \gl^2/(2s)].$)
\qqqed

\medskip
We now turn to \eqref{G''2}.
This is just a union bound but we have to be a little
careful since the most naive bound, $\exp[O(n)]$, on the number
of choices for the $S_i$'s is unaffordable for small $t$.
But this is an overcount:
since $P\sub \C{N(X)\cap W}{2}$ (where $N(X)=\cup_{x\in X}N(x)$),
we have \eqref{G''2} provided
\beq{PcapnablaT}
|G''\cap \nabla(T_1\dots T_{r-1})|
<|P\cap \nabla(T_1\dots T_{r-1})|q+\vt t n^2p^2q
\enq
whenever $(T_1\dots T_{r-1})$
is a partition
of $N(X)\cap W$ (rather than a member of $\cee$);
and, since $|N(X)\cap W|< (1+o(1))tnp$ (using $G'\in\Q^*$),
the number
of possibilities for such a $(T_1\dots T_{r-1})$ is
less than $\exp_{r-1}[(1+o(1))tnp]$.

On the other hand, for
a particular $(T_1\dots T_{r-1})$
we have
\[
|P\cap \nabla(T_1\dots T_{r-1})|
\leq \sum_{x\in X}\Cc{d_W(x)}{2} < tn^2p^2
\]
(again using $G'\in \Q^*$ to bound the $d_W(x)$'s),
which with Theorem~\ref{Chern} implies that the
probability of violating \eqref{PcapnablaT} for a given
$(T_1\dots T_{r-1})$ is less than (say)
$\exp[-\vt^2 tn^2p^2q/3]=\exp[-(\vt^2\eps/3)tn^{3/2}p^2]$.
The union bound
(and the fact that $np = o(n^{3/2}p^2)$)
now completes the argument.\qqqed

\section{Proof of Lemma~\ref{L2.3'}}\label{PL2.3'}

Write $\Q$ for the collection of nonempty $Q\sub \C{V}{2}$
satisfying \eqref{dQx}.
Lemma~\ref{L2.3'} says that w.h.p.
if $Q\in \Q$ and all pairs in $Q$
are bad for the balanced cut $\Pi=(A_1\dots A_{r-1})$
(so by definition $Q\sub \C{A_1}{2}$), then
$\deff_G(\Pi)\geq 2r^2 |Q|$.

We will show something a little stronger.
For $Q\in \Q$, let $\B_Q$ be the set of graphs $H$
for which there is some
$\Pi=(A_1\dots A_{r-1})\in \cee(Q)$
(defined following \eqref{cee}) with
\beq{QQPi}
\mbox{
$Q\sub Q_H(\Pi)$ and
$\deff_{\cee(Q),H}(\Pi)< 2r^2|Q|$.}
\enq
We show
\beq{UQBQ}
\Pr(\cup_Q\B_Q) =o(1)
\enq
(with the union over $Q\in \Q$).
This is stronger than Lemma~\ref{L2.3'} because
$\deff_{\cee(Q),G}(\Pi)$ may be smaller (and is not larger)
than $\deff_G(\Pi)$.

\medskip
Again we can only afford a union bound after restricting
the range of discourse.
Let
$\A$ be the set of graphs $H$
satisfying
\beq{G01}
d_H(x) =(1\pm \gd)np ~~\forall x\in V
\enq
(note this implies
\beq{G01'}
|H| = (1\pm \gd)\Cc{n}{2}p),
\enq
\beq{G04}
|H[S]-|S|^2p/2| < \gd n^2p ~~\forall S\sub V
\enq
and
\beq{G05}
\mbox{$\forall$ $s\in [3,r]$ and
$\{x_1\dots x_s\}\in \Cc{V}{s}$,
$\kappa(x_1\ldots x_s)=o(\gL_r(n,p))$.}
\enq

\mn
Then $\Pr(\ov{\A})=o(1)$ by Propositions~\ref{vdegree}
and \ref{propG[X]}(a)
and Corollary \ref{KrCor}.
(When $s=r$, \eqref{G05} just says that $\glr=\go(1)$.
For smaller $s$ we use Corollary \ref{KrCor}, in which,
with $\gb =(r-s)(s-2)/[2(r+1)]$,
we have
\[
Z~=~ [np^{(s+1)/2}]^{-(s-2)}\glr~<~ n^{-2\gb}\glr
\]
(using \eqref{p}).
Then either $\glr< n^\gb$, implying $Z<n^{-\gb}$,
and the bound $K$ ($=o(\glr)$) in
\eqref{NHG'} applies, or $\glr\geq n^\gb$, in which case the second
bound in \eqref{NHG'} applies and is $o(\glr)$.)

We thus have
\[
\Pr(\cup_Q\B_Q) <
o(1)+  \sum_Q\Pr(\A\cap\B_Q),
\]
and for \eqref{UQBQ} it's enough to show that for each $Q\in \Q$,
\beq{Gnptoshow}
\Pr(\A\cap\B_Q) < \exp [-3|Q|\log n].
\enq
(As elsewhere we just give the bound we need,
but the 3 could be
replaced by any constant if $C$ (in \eqref{p}) is large enough.)

For the rest of this discussion we fix $Q\in \Q$ and set
$\cee(Q)=\cee$; in particular ``rigid," ``core,"
$b(H):=b(\cee,H)$ and $\deff_H:=\deff_{\cee,H}$ now
refer to
this $\cee$.
The main line of argument in this section will
work with $G_{n,M}$ rather than
$G_{n,p}$, but for the moment we stick with the
latter.

Set $\gc'=2\gc$.
(The difference between $\gc$ and $\gc'$ should be ignored;
the extra 2, which could really be $1+o(1)$, is needed to
cover a minor detail at \eqref{k-G}.)

With $Q'$ a particular subset of $Q$ to be specified below, set,
for disjoint $T_1\dots T_{r-2}\sub V\sm V(Q)$,
\[
F(T_1\dots T_{r-2})
=\{\kappa(Q',T_1\dots T_{r-2})< \gc' |Q'|\gL_r(n,p)\}.
\]
Set
\beq{K}
K= 50\ga^{-1}r^3.
\enq

\mn
In a sense our argument attempts---not always
sucessfully---to reduce
\eqref{Gnptoshow} to a situation where
the following statement applies.
\begin{lemma}\label{basic}
Let $\R$ be the set of graphs $H$ satisfying:
$H$ is rigid, say with core
$\{S_1\dots S_{r-1}\}$, $V(Q)\sub S_1$, and $F(S_2\dots S_{r-1})$
holds in $H$.
Then for any $q>(1-2\gd) p$,
$\Pr(G_{n,q}\in \R) <\exp[-(10K\log r+1)|Q'|\log n].$
\end{lemma}

\nin
{\em Remarks.}
We will make sure that $Q'$ is a reasonably large subset
of $Q$---in some cases it will be $Q$ itself---so that
the probability here will be
smaller than the $\exp [-3|Q|\log n]$ of \eqref{GnMtoshow'}.
The reason for the $q$ is that we will see some
graphs $G_{n,M}$ with
$M$ slightly smaller than $\C{n}{2}p$.
The reason for the silly ``$+1$" will appear in
Lemma~\ref{basic''}.

\mn
{\em Proof.}
This follows immediately, {\em via}
Lemma~\ref{transfer}, from the next assertion, which is an
easy consequence of
Theorem~\ref{TJanson} and Lemma~\ref{newHlemma}.
\begin{lemma}\label{basic'}
If $T_1\dots T_{r-2}\sub V\sm V(Q)$ are disjoint
with $|T_1|\dots |T_{r-2}|>n/r$,
then for any $q>(1-2\gd) p$,
\[\Pr(G_{n,q}\models
F(T_1\dots T_{r-2})) <\exp[-(10K\log r+1)|Q'|\log n].\]
\end{lemma}
\nin

\mn
{\em Proof.}
It is of course enough to show this when $q=(1-2\gd )p$.
Notice first that for any fixed $\vt$ we have $\gS< \vt nq/\log n$
for large enough $C$ ($\gS$ as in \eqref{s},
$C$ as in \eqref{p}).
In particular we may assume that $\gD_{Q'} < \vt nq/\log n$,
where $\vt$ is chosen so that the conclusion of
Lemma~\ref{newHlemma} holds with
$\xi = (1/3)(10K\log r+1)^{-1}r^{-2(r-2)}$
and $q$ in place of $p$.

\medskip
Let $\h$ consist of all
sets of the form $K(xy,Z) = \C{\{x,y\}\cup Z}{2}\sm \{xy\}$
with $xy\in Q'$ and $Z\in \C{V}{r-2}$ meeting each of
$T_1\dots T_{r-2}$.  For $K\in \h$,
let $I_K$ be the indicator of $\{K\sub G\}$.
Then
\[
\mbox{$\mu:=\sum \E I_K
= |Q'|\prod_{i=1}^{r-2}|T_i|q^{\Cc{r}{2}-1}
>|Q'|r^{-(r-2)}\glq$}
\]
and, by our choice of $\vt$,
\begin{eqnarray}
\ov{\gD}&:=&\sum\sum \{\E I_KI_L:K,L\in \h, K\cap L\neq \0\}
\nonumber\\
&<& \xi|Q'|\glq^2/\log n.\label{ovgD'}
\end{eqnarray}
Thus, since
$F(T_1\dots T_{r-2})=\{\sum I_K< \gc' |Q'|\gL_r(n,p)\}$
(and $\gc' |Q'|\gL_r(n,p)$ is much smaller than $\mu$;
see \eqref{gcsat}),
Theorem~\ref{TJanson} gives (e.g.)
\begin{eqnarray*}
\Pr(G_{n,q}\models
F(T_1\dots T_{r-2}))& < &\exp[-\mu^2/(3\ov{\gD})]\\
&<& \exp[-\xi^{-1}r^{-2(r-2)}|Q'|\log n].
\end{eqnarray*}\qed

\medskip
We now define $Q'$ and associated sets $W_Q,Z_Q\sub V(Q)$;
these will be used to deal with a minor technical point
involving steps of ``type B" below.
(See the second paragraph of the proof of Lemma~\ref{TLlemma}.)
The choice here depends on the size of $Q$.
If $Q$ is very small, say
$|V(Q)|\leq 13$, then we take $Q'=Q$ and $W_Q=Z_Q=V(Q)$.
Otherwise, we choose $W_Q\sub V(Q)$
with $|W_Q|\leq |V(Q)|/2$ and
$|Q[W_Q]|\geq |Q|/5$
(which is possible because, as is easily verified,
if $W$ is chosen uniformly
from the $\lfloor |V(Q)|/2\rfloor$-subsets of $V(Q)$ then
$\E |Q[W]|\geq |Q|/5$), and
take $Q'=Q[W_Q]$
and $Z_Q=V(Q)\sm W_Q$.

\medskip
From this point we switch to $G_{n,M}$, noting, to begin, that
Lemma~\ref{lemma:twomodels} allows transfer of
Lemma~\ref{basic} to this setting:
\begin{lemma}\label{basic''}
For $\R$ as in Lemma~\ref{basic} and
any $M\geq (1-2\gd)\C{n}{2}p$,
\[\Pr(G_{n,M}\in \R) <\exp[-10K\log r|Q'|\log n].\]
\end{lemma}

\medskip
We assume for the rest of this section that
\beq{Mbig}
M=(1\pm\gd)\Cc{n}{2}p.
\enq
We will show
\beq{UQBQM}
\Pr(G_{n,M}\in \A\cap\B_Q) < \exp [-3|Q|\log n].
\enq
Of course this gives \eqref{Gnptoshow}, since
\beq{transfer2'}
\Pr(G_{n,p}\in \A\cap\B_Q)
\leq \max\{\Pr(G_{n,M}\in \A\cap\B_Q):M= (1\pm \gd)\Cc{n}{2}p\}.
\enq

\medskip
For the proof of \eqref{UQBQM}
we will prefer counting.  Having specified $M$,
let $\g=\g_Q$ be the set of $M$-edge graphs in $\A\cap \B_Q$.
We may rewrite \eqref{UQBQM}
as
\beq{GnMtoshow'}
|\g|< \exp [-3|Q|\log n]\Cc{\C{n}{2}}{M}.
\enq

\medskip
Set
\beq{Ld}
L=K |Q| \log n,  ~~~
\qq= 2r^2|Q|
\enq
(so $\qq$ is the defect bound in \eqref{QQPi};
recall $K$ was defined in \eqref{K}), and
\beq{gb}
\gb =[r\gd n^2p +L](r-1)/M.
\enq
We will need some weak constraints on $\gb$, e.g.
\beq{gbsmall}
\gd < \gb < 2r^2 \gd
\enq
(the upper bound since $L< Kn\gS \log n $ is much smaller than $\gd n^2p $;
see \eqref{s}).

\medskip
Fix some order
``$\prec$" on $\cee$ ($=\cee(Q)$).
We will be interested in sequences
$G_0\dots G_T$ with $G_0\in \g$,
$T\leq L$, and, for
$\Pi=(A_1\dots A_{r-1})$
the first cut as in \eqref{QQPi} (with $H=G_0$) and
$1\leq i\leq L$,

\mn
(a)  if $\crit(G_{i-1})\cap \inte(\Pi)\neq\0$, then $G_i=G_{i-1}-e$
for some $e\in \crit(G_{i-1})\cap \inte(\Pi)$
(recall ``$\crit$" was defined in \eqref{crit}); otherwise:

\mn
(b)  if $G_{i-1}$ is rigid with core
$\{S_1\dots S_{r-1}\}$ and $Q$
is not in the core, then
$G_i=G_{i-1}+e$
for some $e\in (\nabla(Z_Q,U)\cap\ext(\Pi))\sm G_{i-1}$
with $Q\sim U\in \{S_1\dots S_{r-1}\}$,
where
$Q\sim U$ means $V(Q)$ and $U$ are in the same
block of some max cut;

\mn
(c)  if $G_{i-1}$ is not rigid then
$G_i=G_{i-1}-e$
for some $e\in G_{i-1}\cap\inte(\Pi)$;

\mn
(d)  If
$G_{i-1}$ is rigid
with $Q$ in the core
(and $\crit(G_{i-1})\cap \inte(\Pi)=\0$), then $T=i-1$
(and the rest of the sequence is vacuous).

\medskip
We call sequences as above {\em legal.}
The transition from $G_{i-1}$ to $G_i$ is
the {\em $i$th step} of the sequence.
A deletion as in (a) is a step of {\em type A},
an addition as in (b) is a step of {\em type B},
and a deletion as in (c) is a step of {\em type C}.

In what follows we will show, roughly, that each
$G_0\in \g$ is the starting point of ``many" legal sequences
of some length, whereas the total number of
legal sequences of each length is ``small" (so $\g$ is small).

\medskip
For the first of these objectives, we show that there are
many choices for $G_i$ whenever $G_{i-1}$ is as in (c),
and at least one such choice if $G_{i-1}$ is as in (b).
(In reality there are also many choices in the second case---though
not necessarily as many as are guaranteed in the first---but all
we need here is that the process doesn't get stuck in situations
that demand steps of type B.
Of course it cannot get stuck at a step of type A.)


\begin{lemma}\label{AB}
In any legal sequence, fewer than $r\qq$ steps are
of types A and B, and all but at most $r\qq$ indices $i$
satisfy
\beq{ordinarystep}
\mbox{$G_i$ is not rigid and step $i$ is of type $C$.}
\enq
\end{lemma}

\mn
{\em Proof.}
This will follow from the next two assertions.

\mn
{\em Claim} 1.
Each step of type A reduces $\deff(\Pi)$
(that is, if step $i$ is of type A then
$\deff_{G_i}(\Pi) < \deff_{G_{i-1}}(\Pi)$)
and
no step increases $\deff(\Pi)$.

\mn
{\em Claim} 2.
If the $i$th step is of type B
then either

\mn
(i)  for some $j\in [r-2]$,
steps $i+1\dots i+j-1$ are of type B
and step $i+j$ is of type A, or

\mn
(ii)  $T\in \{i\dots i+r-3\}$.

\mn
{\em Proof of Claim} 1.
Deletion of an edge in $\inte(\Pi)$ (as happens in
all steps not of type B)
doesn't affect $|\Pi|$ (that is,
$|\Pi_{G_i}|=|\Pi_{G_{i-1}}|$) and
doesn't increase $b$
(that is, $b(G_i)=b(G_{i-1})$),
so doesn't increase $\deff(\Pi)$.
A step of type A decreases $b$ (and doesn't affect $|\Pi|$),
so decreases $\deff(\Pi)$.
A step of type B increases each of $|\Pi|$ and $b$ by 1, so
doesn't affect $\deff(\Pi)$.\qqqed

\mn
{\em Proof of Claim} 2.
If
$G_{i-1}$ is rigid with core $\{S_1\dots S_{r-1}\}$ and $Q$
is not in the core, then (i)
for each $U\in \{S_1\dots S_{r-1}\}$ there is some max cut
with
$V(Q)$ and $U$ contained in different blocks, and
(ii) $Q\sim U$ for at least two choices of
$U\in \{S_1\dots S_{r-1}\}$.
Say step $i$ is of type B$_j$ if there are exactly $j$
($\in [2,r-1]$)
such $U$'s.  It is enough to show that if this is the case,
then step $i+1$, if taken
(i.e. if $T\neq i$), is either of type
$B_l$ for some $l\leq j-1$, or of type A.

Suppose (w.l.o.g.) that
$G_i=G_{i-1}+e$ with
$e\in \nabla(Z_Q,S_1)\cap \ext(\Pi)$.
Then $G_i$ is rigid with core $\{S'_1\dots S'_{r-1}\}$ satisfying
(i) $S_k'\supseteq S_k$ for each $i$;
(ii) $Q\not\sim S_1'$
(in $G_i$); and
(iii) $\{k:Q\sim S'_k ~\mbox{in}~ G_i\}\sub
\{k:Q\sim S_k ~\mbox{in}~ G_{i-1}\}$.
(Because:
addition of $e$ increases $b$
(as noted in the proof of
Claim 1),
so does not increase the set of max cuts;
this gives rigidity, (i) and (iii),
and also (ii) once we observe
that addition of $e$ doesn't
increase the size of any cut with
$V(Q)$ and $S_1$ in the same block.
Actually, for $k>1$,
$Q\sim S'_k ~\mbox{in}~ G_i$ iff
$Q\sim S_k ~\mbox{in}~ G_{i-1}$, but we don't need this.)

In particular, since $G_i$ is rigid, the $(i+1)$st step,
if taken, must be of
type A or B; and if it is of type B, then (ii) and (iii)
imply that it is type B$_l$ for
some $l\leq j-1$.\qqqed

\medskip
Now to complete the proof of Lemma~\ref{AB}, just notice
that Claim 1 (with the assumption $\deff_{G_0}(\Pi)< d$)
guarantees that there are at most $\qq-1$
steps of type A,
and this together with
Claim 2 implies that all steps of types A and B are contained
in at most $\qq$ intervals of length at most $r-1$
(the extra interval corresponding to (ii) in Claim 2).
This gives the first assertion of Lemma~\ref{AB}
(actually with $r-1$ in place of $r$), and the
second follows since if $i$ violates \eqref{ordinarystep} then either step $i$
or step $i+1$ is of type $A$ or $B$.
\qqqed

\begin{lemma}\label{TLlemma}
For each $G_0\in \g$,
there is some $T\in \{0\dots L\}$ for which
the number of legal sequences $G_0\dots G_T$
is at least $L^{-1}[(1-\gb)M/(r-1)]^{T-r\qq}$.
\end{lemma}

\nin
(Recall $\gb$ was defined at \eqref{gb}.
Of course for small enough $T$ this just says that
there {\em is} a legal sequence.)

\mn
{\em Proof.}
Let $G_0\dots G_{i-1}$ be a legal initial segment
(defined in the obvious way).
If $G_{i-1}$ is as in (c) then the number of possibilities
for $G_i$ is
\begin{eqnarray*}
|G_{i-1}\cap \inte (\Pi)|
&> &|G_0\cap \inte (\Pi)|-L \\
&>& \tfrac{n^2p}{2(r-1)} - (r-1)\gd n^2p-L
\\
&>& \tfrac{M}{r-1} - \tfrac{\gd n^2p}{2(r-1)}
-(r-1)\gd n^2p-L\\
&>& \tfrac{M}{r-1} -r\gd n^2p-L
= (1-\gb)M/(r-1),
\end{eqnarray*}
where the second inequality uses
\eqref{G04} and the third uses $M< (1+\gd)n^2p/2$.

If $G_{i-1}$ is as in (b) then, as noted earlier, we just
want to say there is {\em some} legal choice for $G_i$.
Since $Q$ is not in the core, we have
$Q\sim U$ for at
least two choices of $U\in \{S_1\dots S_{r-1}\}$,
say $S_1$ and $S_2$,
and
$e$ (the edge to be added to $G_{i-1}$)
can be any member of
$\nabla(Z_Q,(S_1\cup S_2)\cap (A_2\cup \cdots\cup A_{r-1}))\sm G_{i-1}$;
so we just need to say this set is nonempty, which is true because:
\[
|\nabla(Z_Q,(S_1\cup S_2)\cap (A_2\cup \cdots\cup A_{r-1}))|
~~~~~~~~~~~~~~~~~~~~~~~~~~~~~~~~~
\]
\beq{nablaZ}
~~~~~~~~~~~~~~
> |Z_Q|\left[\tfrac{2}{r}+ (r-2)\tfrac{(1-\gd)}{r-1}-1\right]n
=
|Z_Q|\cdot\tfrac{r-2}{r-1}\left(\tfrac{1}{r}-\gd\right)n,
\enq
since $|S_1|,|S_2|> n/r
$
and $|A_j|> (1-\gd)n/(r-1)$ for each $j$
(since the $S_j$'s form a core and $\Pi$ is balanced),
while, using
Lemma~\ref{AB}
and \eqref{G01},
\begin{eqnarray*}
|G_{i-1}\cap \nabla(Z_Q,V\sm Z_Q)|
&\leq &
\sum _{x\in Z_Q}d_{G_{0}}(x) +r\qq < |Z_Q| (1+\gd)np +r\qq,
\end{eqnarray*}
which is (much) smaller than the bound in \eqref{nablaZ}
(using \eqref{po1} and
$\qq = 2r^2|Q| < 2r^2|V(Q)|\gS\leq 4r^2|Z_Q|\gS$).

\medskip
Thus, again using Lemma~\ref{AB} (to say a legal sequence
involves at most $r\qq$ steps that are not of type C),
Lemma~\ref{TLlemma}
follows from the next little
(presumably known) observation.

\begin{lemma}\label{treelemma}
Suppose $\T$ is a tree of depth at most $L>0$ and $W$ is a subset
of the internal vertices of $\T$ such that
each internal vertex not in $W$
has at least $\Delta$ children and each path from the root
contains at most s vertices of $W$.
Then there is some $T\in \{0\dots L\}$ for which the number of
leaves at depth $T$ is at least $L^{-1}\gD^{T-s}$.
\end{lemma}
\nin
{\em Proof.}
For each $T<L$ and leaf $w$ at depth $T$, add (to $\T$) a
$\gD$-branching subtree of depth $L-T$
rooted at $w$, forming a tree $\T'$.
Then
$\T'$ has at least $\gD^{L-s} $ leaves
(which are, of course, all at depth $L$), e.g. since
the natural root-leaves random walk down $\T'$
reaches no leaf with probability more than $\gD^{-(L-s)}$.
On the other hand, the number of leaves in $\T'$ is precisely
$\sum_T m_T\gD^{L-T}$, where
$m_T$ is the number of leaves at
depth $T$ in $\T$.
The lemma follows.
\qed

\bigskip
Let $\cup_{T=0}^L\g_T$
be a partition of $\g$ such that for each $T$ and $G_0\in \g_T$
the number of legal sequences $G_0\dots G_T$
is at least $L^{-1}[(1-\gb)M/(r-1)]^{T-r\qq}$.

\medskip
We next give {\em upper} bounds on the numbers of legal sequences
$G_0\dots G_T$ for $T\in \{0\dots L\}$.
For this part of the argument we think of starting with $G_T$
and moving (now mostly by {\em adding} edges) to $G_0$.
For typographical reasons we now set $\C{n}{2}=N$.

Notice that if $G_0\dots G_T$ is a legal sequence
then, by Lemma~\ref{AB} (and the fact that only steps of type $B$
add edges), $M-T\leq |G_T|< M-T+2r\qq$.
Note also---just to keep things slightly cleaner---that
(using \eqref{po1})
\beq{NGT}
\sum_{0\leq i< 2r\qq}\C{N}{M-T+i}< \C{N}{M-T+2r\qq}
\enq
(so the r.h.s. bounds the number of possibilities for $G_T$
for a given $T$).

\medskip
We first consider $T=L$.
Here we use the second assertion of Lemma~\ref{AB}.
If $i$ satisfies \eqref{ordinarystep}
then
$G_i=G_{i-1}-e$ with
$e$ contained in some $(\cee,G_i)$-component
(since adding $e$ to $G_i$ doesn't increase $b$);
so, since $G_i$ is non-rigid, the number of
possibilities for $G_{i-1}$ (given $G_i$) is at most
$(1-\ga )n^2/(2(r-1))$.
This bounds the total number of legal sequences
$G_0\dots G_L$ by
\beq{T=L}
\Cc{N}{M-L+2r\qq}L^{r\qq}n^{2r\qq }
[(1-\ga )n^2/(2(r-1))]^{L-r\qq }.
\enq
Here the first term counts choices of $G_L$, and
the term $L^{r\qq}$ is for specification of a set of at most
$r\qq$ indices
$i$ for which \eqref{ordinarystep} fails
(and for which we use the
trivial bound $n^2$ on the number of possibilities
for $G_{i-1}$ given
$G_i$).

\medskip
We next consider $T<L$.
Here, in contrast to what we did for $T=L$, our
goal is to say that the number of possibilities
for $G_T$ is small.
Suppose $G_T$ has core $\{S_1\dots S_{r-1}\}$
with $V(Q)\sub S_1$.  We show that in this case
\beq{GTsat}
\mbox{$G_T$ satisfies $F(S_2\dots S_{r-1})$.}
\enq

\mn
{\em Proof.} 
Notice that
\beq{nQST}
G_T\cap\nabla(V(Q),S_2\dots S_{r-1})\sub \ext(\Pi)
\enq
(since
\[G_T\cap\nabla(V(Q),S_2\dots S_{r-1})
\sub G_T\cap\nabla(S_1\dots S_{r-1})\sub \crit(G_T)\]
and $\crit(G_T)\cap ~\inte(\Pi)=\0$).
We consider the cases $|V(Q)|>13$ and $|V(Q)|\leq 13$
separately.

If $|V(Q)|>13$ then
\beq{GTcap}
G_T\cap \nabla(W_Q,S_2\dots S_{r-1})\sub G_0
\enq
(since edges added in moving from $G_0$ to $G_T$
meet $V(Q)$ only in $Z_Q=V(Q)\sm W_Q$).
Combining this with \eqref{nQST},
which in particular implies that
\[\nabla_{G_T}(V(Q),S_2\cup\cdots \cup S_{r-1})\cap
\nabla(V(Q),A_1\sm V(Q))=\0,\]
we have
\[
\K_{G_T}(Q',S_2\dots S_{r-1})
\sub \K_{G_0}(Q',A_2\dots A_{r-1})
\]
and thus
\[
\kappa_{G_T}(Q',S_2\dots S_{r-1})
\leq \kappa_{G_0}(Q',A_2\dots A_{r-1}).
\]
Since $Q'\sub Q\sub Q_{G_0}(\Pi)$, this gives \eqref{GTsat}.

\medskip
If $|V(Q)|\leq 13$ (in which case $Q'=Q$ and $W_Q=V(Q)$),
then
we don't quite have \eqref{GTcap}, but can (we assert) say
\beq{k-G}
\kappa_{G_T}(Q,S_2\dots S_{r-1})
\leq \kappa_{G_0}(Q,A_2\dots A_{r-1})+ o(\gL_r(n,p)),
\enq
which again
gives $F(S_2\dots S_{r-1})$ for $G_T$.

For \eqref{k-G},
notice that $|G_T\sm G_0|\leq r\qq=O(1)$
(by Lemma~\ref{AB}, since all edges of $G_T\sm G_0$ are added
in steps of type B).
On the other hand, because of \eqref{nQST}, each
member of $\K_{G_T}(Q,S_2\dots S_{r-1}))
\sm \K_{G_0}(Q,A_2\dots A_{r-1})$
uses one of the $O(1)$ pairs
from $Q$, together with
at least one of the at most $r\qq= O(1)$
edges of $G_T\sm G_0$, so uses
two vertices of $V(Q)$ plus, for some $s\in [3,r]$,
precisely $s$ other vertices incident with edges of $G_T\sm G_0$.
Thus, since the number of possibilities for these $s$ vertices is $O(1)$,
\eqref{k-G} follows from \eqref{G05}.
\qed

Lemma~\ref{basic''}
(applicable since $|G_T|\geq M-L > (1-2\gd) \C{n}{2} p$)
%
%
and \eqref{NGT} now
bound the number of choices for
$G_T$ by $\xi\Cc{N}{M-T+2r\qq}$, where
\beq{xi}
\xi=\exp[-10K\log r|Q'|\log n]
\leq \exp[-2K\log r|Q|\log n],
\enq
so we may crudely bound the number of
legal sequences of length $T$
by
\beq{finalxi}
\xi \Cc{N}{M-T+2r\qq}N^{T}.
\enq
(The $N^{T}$
could of course be improved along the
lines of the above discussion
for $T=L$.)

\medskip
Combining the bounds in \eqref{T=L} and \eqref{finalxi} with the fact that
each $G_0\in \g_T$ is the first term of at least
$L^{-1}[(1-\gb)M/(r-1)]^{T-r\qq }$ legal sequences $(G_0\dots G_T)$,
we have
\begin{eqnarray}
|\g_L| &\leq&
L\left[\tfrac{r-1}{(1-\gb)M}\right]^{L-r\qq }
\Cc{N}{M-L+2r\qq}L^{r\qq}n^{2r\qq }
\left[\tfrac{(1-\ga )n^2}{2(r-1)}\right]^{L-r\qq }.\nonumber\\
&<&
n^{4r\qq }
\Cc{N}{M-L+2r\qq}
\left[\tfrac{(1-\ga)N}{(1-\gb)M}\right]^L\label{gL}
\end{eqnarray}
and, for $T<L$ (with $\xi$ as in \eqref{xi}),
\begin{eqnarray}
|\g_T| &\leq&
L\left[\tfrac{r-1}{(1-\gb)M}\right]^{T-r\qq }
\xi\Cc{N}{M-T+2r\qq}N^T
\nonumber\\
&<&
n^{2r\qq } \xi  \Cc{N}{M-T+2r\qq}
\left[\tfrac{(r-1)N}{(1-\gb)M}\right]^T\label{gT}.
\end{eqnarray}

Thus, noting that \eqref{po1} implies, for any $-2rd\leq i\leq M$,
\beq{MNT}
\Cc{N}{M-i} < [(1+o(1))M/N]^i\Cc{N}{M},
\enq
we have
\beq{GL}
|\g_L| <
n^{4r\qq }(N/M)^{2r\qq}
(1-\ga+\gb)^L
\Cc{N}{M}
<n^{6r\qq }
(1-\ga+\gb)^L
\Cc{N}{M}
\enq
and, for $T<L$,
\beq{GT}
|\g_T| <
n^{2r\qq }(N/M)^{2r\qq}
\left[\tfrac{r-1+o(1)}{1-\gb}\right]^T\xi
\Cc{N}{M}
<n^{4r\qq }
\left[\tfrac{r-1}{1-\gb}\right]^T\xi
\Cc{N}{M}
\enq
(where, to make things a little easier to look at, we used
\eqref{gbsmall} in \eqref{GL} (to say
$(1-\ga+o(1))/(1-\gb)< 1-\ga+\gb$)
and $(1+o(1))^{L}<n^{o(\qq)}$ in \eqref{GT}).

Finally,
summing these bounds and using
\eqref{K}, \eqref{Ld},
\eqref{gbsmall} and \eqref{xi}
gives \eqref{GnMtoshow'}:
\[
|\g| < \left[
n^{6r\qq }(1-\ga+\gb)^L +
Ln^{4r\qq }(\tfrac{r-1}{1-\gb})^L\xi\right]
\Cc{N}{M}
<\exp [-3|Q|\log n]\Cc{N}{M}.
\]
(Here
$(1-\ga +\gb)^L\approx \exp[-50 r^3|Q|\log n]$
dominates $n^{6rd} = \exp[12 r^3|Q|\log n]$,
and in
\begin{eqnarray*}
Ln^{4r\qq }(\tfrac{r-1}{1-\gb})^L\xi
&<& \exp[8r^3 |Q|\log n + L\log r -2K\log r |Q|\log n]\\
&=& \exp[8r^3 |Q|\log n -K\log r |Q|\log n],
\end{eqnarray*}
the term $8r^3 |Q|\log n$ in the exponent is negligible.)
\qqqed

\section{Remarks}\label{Remarks}

{\bf A.}
An obvious (but probably formidable) challenge is to prove Theorem~\ref{MT}
with the correct $C$.  The natural guess 
is that
\[
C> [2r/(r+1)]^{\frac{2}{(r+1)(r-2)}}
\]
suffices,
this being what's needed
to guarantee that (w.h.p.) all edges lie in $K_r$'s.
Note, though, that
the even more precise
``stopping time" version---which
says that if
we choose $e_1,\ldots \in E(K_n)$,
with $e_i$ uniform from edges as yet unchosen, and stop as
soon as every $e_i$ is in a $K_r$, then w.h.p. the resulting
$G$ satisfies $t_r(G)=b_r(G)$---is
{\em not} correct.
To see this (informally),
suppose $xy$ is the last edge
added in forming $G$.  There is then some $uv\in G$
such that every $K_r$ on $uv$ also contains $xy$.  But in this case
$t_r(G)>b_r(G)$ whenever there is a maximum cut
with (for example) $u,v$ and $x$ in a single block,
and this is not a low probability event.

\mn
{\bf B.}
For a property (i.e. isomorphism-closed
collection) $\f$ of graphs on $[n]$,
set
$\mu_p(\f)=\Pr(G_{n,p}\in \f)$,
and define the {\em threshold} for $\f$ to be
\[
p_c(\f):=\min\{p:\mu_q(\f)\geq 1/2 ~~\forall q\geq p \}.
\]
If $\f$ is increasing (preserved by addition of edges) then
$\mu_p(\f)$ is increasing in $p$ and $p_c(\f)$ is that
$p$ (unique except in trivial cases)
for which $\mu_{p_c}(\f)=1/2$.  (This is not the original definition of threshold in
\cite{ER}
but is more or less equivalent.
Of course these notions make sense more generally, but for this brief discussion
we stick to graph properties.)

The property
$\f_r:=\{\ttt_r(G)=\bb_r(G)\}$ of
Theorem~\ref{MT}
is not increasing and $\mu_p(\f_r)$ is
is not increasing in $p$
(for a given $n$);
rather it is close to 1 for $p$ either
large enough or quite small,
and is easily seen to be close to 0
for some intermediate values.
Still, there is an interesting possibility
(which for $r=3$ was suggested in \cite{DKMantel}), namely,
{\em could it be that, for given $r $ and $n$,
$\mu_p(\f_r)$
has just
one local minimum?}

In fact it would seem to be interesting to prove such a statement for
{\em any} (natural) nonincreasing property,
and similarly interesting to identify some natural
situation(s) in which $\mu_p(\f)$ {\em is} increasing
although $\f$ is not;
might this, for example, be true of the property
$\{t_r(G)< (1-1/(r-1)+\eps)|G|\}$ ({\em cf.} Theorem~\ref{8.34})?

\mn
{\bf C.}
Another obvious question is, does
Theorem~\ref{MT} extend to graphs $H$ other than cliques; that is, if
$t_H(G)$ and $b_H(G)$ are the maximum
values of $|K|$ for $K$ ranging over,
respectively, $H$-free and ($\chi(H)-1$)-partite subgraphs of $G$
(where $\chi$ is chromatic number),
when is $G_{n,p}$ likely to satisfy
\[
\f_H:=\{t_H(G)=b_H(G)\}~?
\]
It is easy to see
that the question only makes sense when $H$ is {\em critical}, that is,
contains a ({\em color-critical}) edge
$e$ such that $\chi(H-e)<\chi(H)$.
As noted in Section~\ref{Intro},
the
result of \cite{BSS} mentioned there
holds in this generality,
and it is suggested by the authors of \cite{BPS} that their main result
(Theorem~\ref{Tbps} above) should as well.
Here again there is a natural guess.
Say $\g_H$ holds for $G$ if each $e\in E(G)$ is color-critical in some copy of $H$ in $G$.
%
\begin{conj}\label{Hconj}
For any H with a color-critical edge,
$p_c(\f_H)=O( p_c(\g_H))$.
\end{conj}
\nin
(An old theorem of M. Simonovits~\cite{Simonovits}
says that if $H$ is critical then $K_n$ satisfies $\f_H$
for large enough $n$.)
For $H=K_r$, Conjecture~\ref{Hconj} is Theorem~\ref{MT}.
The threshold for $\g_H$ is not a mystery, but takes some
space and is omitted here.
One may also guess ({\em cf.} {\bf A} above) that in fact $p_c(\f_H)\sim p_c(\g_H)$.

\mn
{\bf D.}
Trivially necessary for $t_r(G)=b_r(G)$
is that every max(imum) cut be maxi{\em mal} $K_r$-free;
thus Theorem~\ref{MT} implies
that for $p$ in our
range this condition holds w.h.p.
This is in fact an instance of  Lemma~\ref{L2.3}, but
{\em is there an easier way to prove it?}
For $r=3$ (where the requirement is that for each
max cut $(A,B)$ and $xy\in G[A]$, $x$ and $y$ share
a neighbor in $B$), the proof implicit in
\cite{DKMantel} is simple once found; but finding it
was the real key to that paper.

Here we are back to the difficulties associated with max cuts
({\em cf.} the beginning of Section~\ref{RandC}).
On this theme, a simple question suggested by the present work
is:
{\em for what $p$ is it true that
$G_{n,p}$ (w.h.p.) admits no max cut $(A_1\dots A_{r-1})$
such that some $x$ has all its neighbors in a single $A_i$?}

When $r\geq 4$, the proof of
Lemma~\ref{Xdefect} can be adapted to give this for
$p>C_rn^{-1/2}$,
but it should really be both
easier and true for considerably smaller $p$, perhaps requiring only
$p\gg n^{-1}\log n$.
For $r=3$ we don't even know that
$p>Cn^{-1/2}$ is enough, though
the same guess seems reasonable:

\begin{conj}\label{CutConj}
If $p\gg n^{-1}\log n$, then w.h.p. no (ordinary) max cut of $G_{n,p}$
contains all (or even 51\% of) the edges at any vertex.
\end{conj}
\nin
Thus $p$ should be large enough that a {\em typical} cut
contains
only about half the edges at any vertex; a {\em max} cut will
of course
tend to contain more, but the guess
is that this effect is relatively mild.
(It follows from \cite[Theorem~1.4]{BPS} that the conclusion holds
for $p$ at least about $n^{-1/3}\log^{2/3}n$.)

\mn
{\bf E.}
Long as the above argument is, the full proof of Theorem~\ref{MT} is
longer still, in that we {\em begin} with the already very difficult
assertion (Theorem~\ref{8.34}) that every large enough $F\sub G=G_{n,p}$
is nearly $(r-1)$-partite.
Note, though, that we really only need this when $F$ is {\em maximum} $K_r$-free
(and for $p$ slightly larger than what's specified in \eqref{p'}).
In fact both \cite{BSS} and \cite{BPS} (which of course preceded \cite{Conlon-Gowers})
begin with such limited
versions of Theorem~\ref{8.34}, and it would be interesting to see
whether a version adequate to present purposes could be proved relatively easily.

\bn
Department of Mathematics\\
Rutgers University\\
Piscataway NJ 08854\\
rvdemarco@gmail.com\\
jkahn@math.rutgers.edu

\end{document}